\documentclass[11pt,a4paper]{article}
\usepackage[utf8]{inputenc}
\usepackage{amsmath,bm}
\usepackage{amsmath}
\usepackage{amsthm}
\usepackage{amsfonts}
\usepackage[square, comma, sort&compress,numbers]{natbib}
\usepackage{amssymb}
\usepackage{setspace}
\usepackage{graphicx}
\usepackage{color}
\usepackage[marginal]{footmisc}
\usepackage[section]{placeins}

\usepackage{indentfirst}
\usepackage[left=2cm,right=2cm,top=2cm,bottom=2cm]{geometry}
\author{}
\usepackage{float}
\newtheorem{rem}{Remark} 

\newtheorem*{defi*}{Definition}
\newtheorem{cor}{Corollary}
\newtheorem{lem}{Lemma}
\newtheorem{thm}{Theorem}

\newtheorem*{lem*}{Lemma}
\newtheorem*{thm*}{Theorem}
\allowdisplaybreaks
\usepackage{authblk}
\makeatletter
\newcommand*{\rom}[1]{\expandafter\@slowromancap\romannumeral #1@}
\makeatother
\title{{\bf Persistent heterodimensional cycles in periodic perturbations of Lorenz-like attractors}}
\date{}
\author{Dongchen Li}
\author{Dmitry Turaev}
\affil{{Department of Mathematics, Imperial College London}}
\begin{document}
\def\D{\mathrm{D}}
\def\d{\mathrm{d}}
\def\tr{\mathrm{tr}\,}
\def\diff{Di\!f\!f}
\def\Poincare{Poincaré\,\,}
\bibliographystyle{plainnat}

\date{}
\maketitle

\footnote{This work is supported by the EPSRC grant EP/PO26001/1.}
\par{}
\noindent{\bf Abstract.} 
We prove that heterodimensional cycles can be created by unfolding a pair of homoclinic tangencies in a certain class of $C^r$-diffeomorphisms $(r=3, \dots, \infty,\omega)$. This implies the existence of a $C^2$-open domain in the space of dynamical systems with a certain type of symmetry where systems with heterodimensional cycles are dense in $C^r$. In particular, we describe a class of three-dimensional flows with a Lorenz-like attractor such that an arbitrarily small time-periodic perturbation of any such flow can belong to this domain - in this case the corresponding heterodimensional cycles belong to a chain-transitive attractor of the perturbed flow.
\par{}
\par{}

\section{Main results}\label{sec:intro}

A heterodimensional cycle is formed by intersections between invariant manifolds of hyperbolic periodic orbits of different indices (dimensions of unstable manifolds). By this definition, they only appear in dimension three or more for diffeomorphisms, or dimension four or higher if we consider systems of autonomous differential equations. Heterodimensional cycles in such dynamical systems create a basic mechanism that causes non-hyperbolicity and breaks structural stability. Early examples involving heterodimensional cycles were studied by Abraham and Smale \citep{as70} and Shub \citep{shub71}. Later on, a systematic study was carried out by Diaz and his collaborators in \citep{dr92,di95a,di95b,bd96}. Bonatti and Diaz built in \citep{bd08} a comprehensive theory of $C^1$ diffeomorphisms having heterodimensional cycles of co-index one (i.e., when the difference between the indices is one). They also showed the $C^1$-robustness of heterodimensional cycles - a $C^1$-small perturbation of a system with a heterodimensional cycle can always be constructed such that the perturbed system gets into a $C^1$-open domain in the space of dynamical systems where systems with heterodimensional cycles are dense (in $C^\infty$ or $C^\omega$ sense). A general higher smoothness version of this result is missing and a $C^r$ theory (with $r>1$) of perturbations of heterodimensional cycles is much less developed (see, however, \citep{dr92,di95a,di95b,ks12,ast17,ast18}).

The aim of this work is to provide more examples where heterodimensional cycles appear naturally in multidimensional systems. In particular, we show that heterodimensional cycles can be born out of a certain type of homoclinic tangencies (after a $C^r$-small perturbation, for an arbitrarily large $r$, including the case of perturbations small in the real-analytic sense). As homoclinic tangencies persist in the so-called Newhouse domains ($C^2$-open regions in the space of dynamical systems where systems with homoclinic tangencies are $C^r$-dense for every $r\geq 2$ \cite{gts93,pv94}), this gives us the persistence of heterodimensional cycles in the corresponding type of the Newhouse domain.

Our main application is the problem of a periodic perturbation of Lorenz-like attractors. There are different approaches to Lorenz Attractors, e.g. the Guckenheimer-Williams \citep{gw79} and Afraimovich-Bykov-Shilnikov \citep{abs77,abs83} geometric models, and modern generalisations in \citep{mpp04}.
Here we understand the Lorenz attractor as an object described by the Afraimovich-Bykov-Shilnikov geometrical model \citep{abs77,abs83}. This means that we take an autonomous system of ODEs that has a saddle equilibrium state with a one-dimensional unstable manifold. We take a cross-section to the stable manifold and assume that all orbits that start from the cross-section return to its inner part in a positive time (except for the orbits that start from the stable manifold - these tend to the equilibrium state). We also assume uniform hyperbolicity for the return map to the cross-section (exact conditions for that can be written as in \citep{abs77,abs83}). A small neighbourhood of the closure of the set of all orbits that start from the cross-section is a strictly forward-invariant region (an absorbing domain). The attractor inside this domain is the Lorenz attractor in the Afraimovich-Bykov-Shilnikov sense. In \citep{tuc99,tuc02}, it was checked with the use of rigorous numerics that the classical Lorenz system satisfies the conditions of \citep{abs77,abs83}. The same is true for an open set of parameter values in the Morioka-Shimizu model \citep{ctz17} and the extended Lorenz model \citep{ot17}. 

The Morioka-Shimizu model and the extended Lorenz model are important because they serve as normal forms for several codimension-3 bifurcations of equilibrium states which have three Lyapunov exponents simultaneously equal to zero, in systems with certain types of $Z_2$-symmetry \citep{sst93,pst98}. Therefore, the existence of the Lorenz-like attractor in these normal forms also implies that the Lorenz-like attractor is born at the unfolding of such ``triple instability'' bifurcations in an arbitrary system of differential equations. 

More importantly (see \citep{sst93}), the same systems serve as normal forms for some codimension-3 bifurcations of periodic orbits (with 4 zero Lyapunov exponents - one Lyapunov exponent is always zero for a periodic orbit, so having 3 more zero Lyapunov exponents is a codimension-3 bifurcation). This means that some iteration of the \Poincare map near any periodic orbit undergoing such triple instability bifurcation is close (in appropriately chosen coordinates) to the time-1 map of the flow of the corresponding normal form. 
It is the same as to say that some iteration of the Poincare map is the period map of some time-periodic perturbation of this normal form. Since these particular normal forms, as we mentioned, have a Lorenz-like attractor for a certian region of parameter values, these bifurcations give rise to attractors obtained by applying a small time-periodic perturbation to a Lorenz-like attractor. Multidimensional systems of differential equations can have an unbounded number of periodic orbits, any of which can undergo the ``triple instability'' bifurcations which we discuss here, provided there are at least three bifurcation parameters and the flow does not contract three-dimensional volumes (so there is no effective reduction to a low-dimensional case). Different scenarios where these bifurcations happen and the system acquires one or several periodically perturbed Lorenz-like attractors are presented in \cite{gst96,gst08,gst09,ggkt,go17}.

The question of a time-periodic perturbation of the Lorenz-like attractors is also interesting in its own right. A general theory proposed in \citep{ts08} asserts that after any sufficiently small time-periodic perturbation is applied to a system with a Lorenz-like attractor the period map will have a unique chain-transitive attractor
$\cal A$. The equilibrium state of the non-perturbed system becomes the saddle fixed point of the period map, and this fixed point, along with its unstable manifold, belongs to $\cal A$. The unstable manifold may have homoclinic tangencies to the stable manifold. In this paper, we give conditions, under which an arbitrarily small perturbation of such tangencies can create a heterodimensional cycle that involves the fixed point (with the one-dimensional unstable manifold) and another saddle periodic orbit with a two-dimensional unstable manifold. It follows from the results of \citep{ts08}, that when the heterodimensional cycle containing the fixed point exists, it lies in $\cal A$, and the entire unstable manifolds of both its periodic points also lie in $\cal A$. This underscores very non-trivial dynamics in the attractor. In particular, since the attractor $\cal A$ contains saddles with different numbers of positive Lyapunov exponents (1 and 2), the relevance of Lyapunov exponents computations for the understanding of chaos represented by such attractors is questionable (e.g. the shadowing property could be violated \citep{dgsy27}). 

In our analysis we do not need to be restricted to the case of periodically perturbed Lorenz-like system only, we just need to assume the existence of a particular type of homoclinic tangencies. Namely, 
denote by $\diff^r(\mathcal{M})$ the space of $C^r$-diffeomorphisms on a $D$-dimensional manifold $\mathcal{M}$, where $r= 3,\dots,\infty,\omega$ and $D \geqslant 3$ unless otherwise specified.
 Let $F\in\diff^r(\mathcal{M})$ satisfy the following conditions.

(C1) $F$ has a saddle periodic point $O$ with multipliers $\gamma$, $\lambda$, $\lambda_1,\dots, \lambda_{D-2}$ such that $\lambda$ and $\gamma$ are real,
\begin{equation}\label{eq:multipliers0}
|\mathrm{Re}\,\lambda_{D-2}|<\dots<|\mathrm{Re}\,\lambda_1|<|\lambda|<1<|\gamma|
\end{equation}
and
\begin{equation}\label{eq:multipliers}
|\lambda\gamma|>1.
\end{equation}

(C2) There exist two orbits $\Gamma$ and $\tilde{\Gamma}$ of quadratic homoclinic tangency between the unstable and stable manifolds of $O$.

In order to formulate the next condition, we recall some definitions. Denote by $W^{uE}(O)$ a two-dimensional invariant manifold tangent to the eigenspace corresponding to $\lambda$ and $\gamma$ -- the unstable and weak stable multipliers of $O$, and call it the extended unstable manifold of $O$. This manifold is not unique, but it contains $W^u(O)$ and any two of these manifolds are tangent to each other at every point of $W^u(O)$. Recall also that for any diffeomorphism satisfying (C1) there is a unique strong-stable $C^r$-foliation $\mathcal{F}_0$ in the stable manifold $W^s(O)$ which includes, as a leaf, the strong-stable manifold $W^{ss}(O)$ (tangent at $O$ to the eigenspace corresponding to the multipliers smaller than $\lambda$ in the absolute value). Detailed discussion can be found in Chapter 13 of \citep{sstc2} or in \cite{tu96}.

Assume the diffeomorphism $F$ satisfies the following non-degeneracy assumption.

(C3) The homoclinic orbits $\Gamma$ and $\tilde\Gamma$ do not lie in
$W^{ss}(O)$, and the manifold $W^{uE}(O)$ is transverse to the strong-stable foliation $\mathcal{F}_0$ at the points of $\Gamma$ and $\tilde\Gamma$ (in particular, $W^{uE}(O)$ is transverse to the stable manifold $W^s(O)$ at the points of $\Gamma$ and $\tilde\Gamma$).

Observe that if we add any $C^2$-small perturbation to $F$ without destroying the homoclinic tangencies, the tangencies will remain quadratic and also condition (C3) will remain fulfilled. 

Note that conditions (C1) and (C3) imply that the set consisting of the saddle O, and the two homoclinic orbits $\Gamma$ and $\tilde\Gamma$ is partially hyperbolic. Therefore, the foliation $\mathcal{F}_0$ can be smoothly extended to a neighbourhood of 
$O\cup \Gamma\cup\tilde \Gamma$, see \cite{tu96}.

It should be noticed that a single homoclinic tangency is not enough for creating heterodimensional cycles in diffeomorphisms of the type considered in this paper, i.e., those having a saddle with real multipliers being closest to the imaginary axis. It is shown in \cite{gst08} that periodic orbits of different indices can be obtained by unfolding a single orbit of homoclinic tangency. However, these points and $O$ cannot form heterodimensional cycles since they all lie in a certain two-dimensional invariant manifold (see \citep{tu96}) while heterodimensional cycles require at least 3-dimensional ambient space. Therefore, we must consider an interplay between two orbits of homoclinic tangency. This is similar to the results of \citep{li16,lt17} where we obtained heterodimensional cycles by perturbations of a pair of homoclinic loops to a saddle-focus equilibrium state.

A way to make homoclinic tangencies come in pairs is to assume a symmetry in the system. Note that Lorenz-like systems that motivate this work do possess symmetry, so when such system has a homoclinic loop it also has a second one. When we add a periodic perturbation that keeps the symmetry, the pair of homoclinic loops can transform to a symmetric pair of homoclinic tangencies of the type we consider here.

The diffeomorphism $F$ is $\mathbb{Z}_2$-symmetric if there exists a $C^r$-diffeomorphism $\mathcal{R}$ such that $\mathcal{R}^2=id$ and $\mathcal{R}\circ F= F \circ \mathcal{R}$. In order to describe our assumptions on the involution $\mathcal{R}$, consider a small neighbourhood $V$ of the point $O$. We assume that the orbit of $O$ is symmetric with respect to $\mathcal{R}$, so $\mathcal{R} O=O$. It is well-known that one can choose coordinates in $V$, with $O$ at the origin, such that  $\mathcal{R}$ will be linear in these coordinates (a nonlinear involution $v\mapsto \mathcal{R}(v)$ becomes linear:  $v^{new}\mapsto \mathcal{R}_0 v^{new}$, after the coordinate transformation $v^{new}=(v+\mathcal{R}_0 \mathcal{R}(v))/2$, where $\mathcal{R}_0$ is the derivative of $\mathcal{R}$ at zero). Choose such coordinates $v$. Let $\tau$ be the period of the point $O$. As the linear map $\mathcal{R}$ commutes with the derivative 
$\D F^{\tau}$ at $O$, the invariant subspaces of $\D F^{\tau}|_{O}$
are invariant with respect to $\mathcal{R}$ too. Denote $v=(x,y,z)$ where the $x$-, $y$-, and $z$- spaces are the eigenspaces of 
$\D F^{\tau}|_{O}$ corresponding to $\lambda$, $\gamma$, and the rest of the multipliers $\lambda_i$, respectively. As we mentioned, the $x$-, $y$- and $z$-spaces are invariant under $\mathcal{R}$. We assume that $\mathcal{R}:(x,y,z)\mapsto(\bar x,\bar y, \bar z)$ in $V$ acts in the following way:
\begin{equation}\label{eq:symmetryR}
\bar{x} = x, \quad
\bar{y} =-y,\quad
\bar{z} = \mathcal{S}z,
\end{equation}
where $\mathcal{S}$ is a linear involution that changes the signs of some of $z$-coordinates.

Denote by $\diff ^r_s(\mathcal{M})$ the subspace of $\diff ^r(\mathcal{M})$ consisting of $\mathcal{R}$-symmetric diffeomorphisms. Maps that are close to $F$ in $\diff ^r(\mathcal{M})$ (in particular, the maps that are close to $F$ in $\diff ^r_s(\mathcal{M})$) have a saddle periodic point, a hyperbolic continuation of $O$, that continuously depends on the map; its stable and unstable manifolds also depend on the map continuously. Those of these maps that have orbits of homoclinic tangency close to $\Gamma$  form a codimension-1 surface $\mathcal{H}$ in $\diff ^r(\mathcal{M})$. For the maps that belong to the surface  $\mathcal{H}\cap \diff ^r_s(\mathcal{M})$ we also have a symmetric to $\Gamma$ orbit of homoclinic tangency, $\tilde\Gamma$; conditions (C1)-(C3) are fulfilled for every map in this surface. One can define a functional $\mu$ in a neighbourhood of $F$ in $\diff ^r(\mathcal{M})$ such that
${d\mu(F_\varepsilon)}/{d\varepsilon}\neq 0$ for any one-parameter family $F_\varepsilon$ of maps in $\diff ^r(\mathcal{M})$, which is transverse to the surface $\mathcal{H}$, and 
$|\mu(F_\varepsilon)|$ measures the distance between the unstable and stable manifolds of $O$ near a certain point of $\Gamma$. Thus, the surface $\mathcal{H}$ is given by the equation $\mu=0$. Another functional we need is $\theta=-\ln|\left(\lambda|/|\gamma|\right)$ (it is a modulus of topological conjugacy and is known to play an important role in bifurcations of homoclinic tangencies \cite{gs93}). We consider any two-parameter family $F_{\varepsilon_1,\varepsilon_2}$ of diffeomorphisms from $\diff ^r_s(\mathcal{M})$ (so all diffeomorphisms in the family are symmetric) such that 
$F_{\varepsilon_1^*,\varepsilon_2^*}$ equals to the map $F$, and assume that
$$\det\frac{\partial(\mu(F_{\varepsilon_1,\varepsilon_2}),\theta(F_{\varepsilon_1,\varepsilon_2}))}
{\partial(\varepsilon_1,\varepsilon_2)}\neq 0.$$
This condition means that we can consider $\mu(\varepsilon_1,\varepsilon_2)$ and 
$\theta(\varepsilon_1,\varepsilon_2)$ as new parameters, so we further use the notation $F_{\mu,\theta}$
for the chosen family. Let $\theta^*$ be the value of $\theta$ for the original diffeomorphism $F$, so
$F=F_{0,\theta^*}$.

We also need one more ($C^1$-open) condition on the multipliers of $O$:

(C4) $|\lambda_1|<\lambda^2$ and $|\lambda||\gamma|^{\frac{1}{2}}<1$. 

We do not know if Theorem \ref{thm1} below holds without this condition, but our proof uses it in an essential way.

We can now state the main result of the paper.
\begin{thm}\label{thm1}
Let $\{F_{\mu,\theta}\}$ be the two-parameter family of diffeomorphisms in $\diff ^r_s(\mathcal{M})$ such that $F_{0,\theta^*}$ satisfies conditions (C1) - (C4). Then, there exists a sequence $\{(\mu_j,\theta_j)\}$ accumulating on $(0,\theta^*)$ such that for any sufficiently large $j$ the diffeomorphism $F_{\mu_j,\theta_j}$ has a symmetric pair of heterodimensional cycles, each of which includes the index-1 saddle periodic point $O$  and some index-2 saddle periodic point.
\end{thm}

\begin{figure}[!h]
\begin{center}
\includegraphics[width=0.50\textwidth]{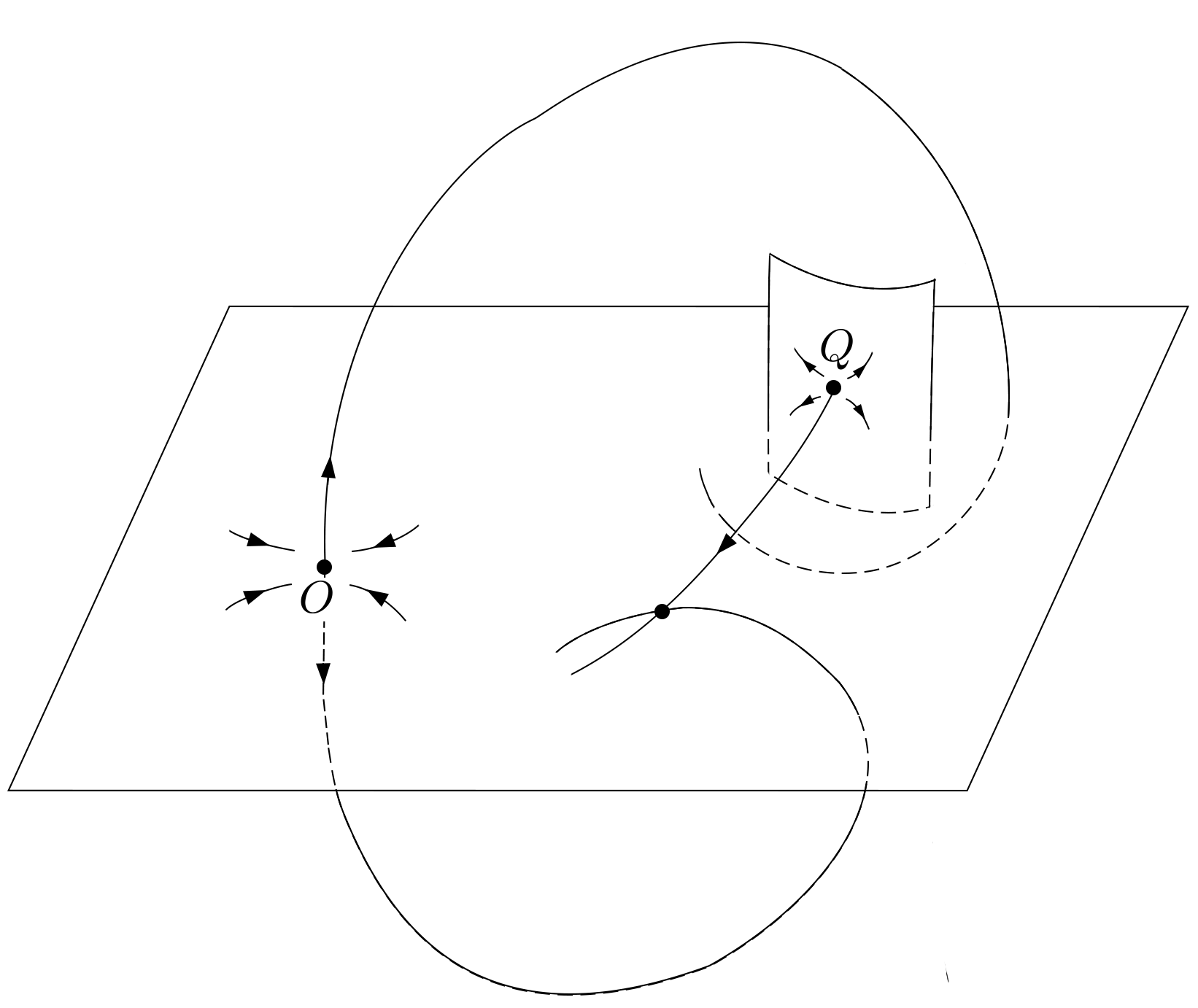}
\end{center}
\caption{A heterodimensional cycle can be obtained by splitting the homoclinic tangencies while changing $\theta$.}
\label{fig:hdc}
\end{figure}

Let us sketch the proof of this theorem. First, by changing $\mu$, we destroy the original homoclinic tangency and obtain a new one, $\hat{\Gamma}$, such that transverse homoclinics to $O$ will exist near $\hat\Gamma$ and also some additional properties are satisfied by $\hat\Gamma$ (see Lemma \ref{lem:preperturb2}). It is known (cf. \cite{gs86}) that by changing $\theta$ one can create a saddle
orbit $Q$ of index 2 near $\hat\Gamma$ (condition $|\lambda\gamma|>1$ is crucial here, as it implies expansion of areas transverse to the strongly contracting directions). By using the existence of transverse homoclinics to $O$, we prove that for any index-2 saddle periodic point near $\hat{\Gamma}$, its unstable manifold will intersect $W^s(O)$ (see Lemma \ref{lem:transverse}). Finally, we show that, by changing $\mu$ and $\theta$ together, the index-2 saddle periodic point $Q$ can be found such that that $W^s(Q)$ intersects the piece of the unstable manifold of $O$ near the orbit of homoclinic tangency which is symmetric to $\hat\Gamma$ (see Lemma \ref{lem:nontransverse}). In order to be able to do this, we need to have $W^s(Q)$ sufficiently ``straight'', which we achieve using condition (C4). The obtained existence of both intersections of $W^s(Q)$ with $W^u(O)$ and $W^u(Q)$ with $W^s(O)$ means the existence of the heterodimensional cycle involving $O$ and $Q$ (see Figure \ref{fig:hdc}).

Recall that the Newhouse region in $\diff ^r(\mathcal{M})$ is an open set comprised by diffeomorphisms having the so-called wild-hyperbolic set \citep{newhouse1}. Systems with homoclinic tangencies are dense in the Newhouse region. Moreover, any family of diffeomorphisms which is transverse to a codimension-1 surface filled by diffeomorphisms which have a saddle periodic point $O$ with a qudratic homoclinic tangency which satisfies the non-degeneracy conditions described in (C3) intersects the Newhouse region over an open set of parameter values, so parameter values corresponding to the existence of quadratic homoclinic tangencies to the hyperbolic continuation of $O$ are dense in these regions and the non-degeneracy conditions (C3) are fulfilled for these tangencies \citep{gts93}. Since our family $F_{\mu,\theta}$ is transverse to the codimension-1 surface $\mathcal{H}\cap \diff ^r(\mathcal{M})$, it follows that we have open regions in the $(\mu,\theta)$ plane
where the parameter values are dense for which the map $F_{\mu,\theta}$ has a symmetric pair of homoclinic tangencies satisfying conditions (C1)-(C4). Thus, Theorem \ref{thm1} implies the following result on the Newhouse region in $\diff ^r_s(\mathcal{M})$:

\begin{cor}\label{cor1}
There exist open sets in the plane of parameters $(\mu,\theta)$ where parameter values corresponding to the existence of a pair of symmetric homoclinic tangencies to $O$ are dense, and parameter values corresponding to the existence of heterodimensional cycles involving $O$ and an index-2 saddle periodic point are dense in these sets.
\end{cor}

Let us now consider the case without symmetry. Then, the simultaneous existence of two homoclinic tangencies given by condition (C2) is a codimension-2 phenomenon. Each of these homoclinic tangencies
can be split independently, so we can introduce two splitting parameters, $\mu_1$ and $\mu_2$, which
measure the distance between the stable and unstable manifolds near a point of $\Gamma$ and,
respectively, a point of $\tilde\Gamma$. As we have more parameters which we can perturb independently,
the result analogous to Theorem \ref{thm1} becomes easier to obtain. In particular, we do not make
assumption (C4) in the non-symmetric case. However, we need one more condition, without which the
birth of heterodimensional cycle from the pair of homoclinic tangencies satisfying (C1)-(C3) will be impossible.

Recall that a uniquely defined smooth strong-stable foliation $\mathcal{F}_0$ exists in the stable manifold of $O$. The homoclinic orbits $\Gamma$ and $\tilde\Gamma$ lie in $W^s(O)$, so for each point of these orbits there is a uniquely defined leaf of $\mathcal{F}_0$ which passes through this point. Assume that the following ``coincidence condition'' holds:

(C5) There is a leaf of $\mathcal{F}_0$ which contains, simultaneously, a point of $\Gamma$ and a point of $\tilde{\Gamma}$.

Note that if condition (C5) is not satisfied, then both orbits of homoclinic tangency will be contained in the same three-dimensional invariant manifold \citep{tu96} and, therefore, no heterodimensional cycles can be born near them. 
So, condition (C5) is necessary for the creation of heterodimensional cycles. This condition is automatically fulfilled in the symmetric case (when the involution $\mathcal{R}$ near $O$ preserves the orientation in the weak stable direction $x$, as given by (\ref{eq:symmetryR})). However, in the general case this is an additional equality-type condition, which makes the bifurcation under consideration a bifurcation of codimension 3. In principle, when we consider perturbations of systems satisfying conditions (C1)-(C3) and (C5), we may consider the distance between the nearest leaves of the foliation $\mathcal{F}_0$ passing through the points of $\Gamma$ and
$\tilde\Gamma$ as an independent bifurcation parameter. We, however, do not need this and consider an
arbitrary 2-parameter unfolding $F_\varepsilon$, with $\varepsilon=(\varepsilon_1,\varepsilon_2)$, of the map $F$ satisfying (C1)-(C3) and (C5), for which we require only that
$$\det\frac{\partial(\mu_1(F_\varepsilon),\mu_2(F_\varepsilon))}
{\partial(\varepsilon_1,\varepsilon_2)}\neq 0.$$
Thus, we can choose $(\mu_1,\mu_2)$ as new parameters.

The same strategy we used for the proof of Theorem \ref{thm1} gives us the following

\begin{thm}\label{thm2}
Let $\{F_{\mu_1,\mu_2}\}$ be a two-parameter family of diffeomorphisms in $\diff ^r(\mathcal{M})$ such that $F_{0,0}$ satisfies conditions (C1)-(C3) and (C5). Then, there exists a sequence $(\mu^1_j,\mu^2_j)\to 0$ such that for every sufficiently large $j$ the diffeomorphism $F_{\mu^1_j,\mu^2_j}$ has a heterodimensional cycle including a hyperbolic continuation of the index-1 saddle periodic point $O$ and
an index-2 saddle periodic point.
\end{thm}

Now we can return to periodically perturbed Lorenz-like systems. Examples of such systems are the classical Lorenz model \citep{lo63}
\begin{equation}\label{eq:intro:0}
\left\{\begin{array}{rcl}
\dot{x}&=&\sigma(y-x), \\
\dot{y}&=&x(\rho-z)-y, \\
\dot{z}&=&xy- \beta z,
\end{array}\right.
\end{equation}
and the Morioka-Shimizu model \cite{sm80}
\begin{equation}\label{eq:intro:1}
\left\{\begin{array}{rcl}
\dot{x}&=&y, \\
\dot{y}&=&x(1-z)-\lambda y, \\
\dot{z}&=&-\alpha z +x^2.
\end{array}\right.
\end{equation}
A computer-assisted proof for the existence of Lorenz attractor in system (\ref{eq:intro:0}) for
the values of parameters $(\sigma,\rho,\beta)$ close to $\sigma=10,\rho=28,\beta=8/3$ was given
in \citep{tuc99,tuc02} and, in \citep{ctz17}, for system (\ref{eq:intro:1}) for an open set of $(\alpha,\lambda)$ near $\alpha=0.606, \lambda=1.045$. Recall that by Lorenz attractor we mean the attractor in the sense of Afraimovich-Bykov-Shilnikov (ABS) model, see \citep{abs77,abs83}. 

\begin{figure}[!h]
\begin{center}
\includegraphics[width=0.70\textwidth]{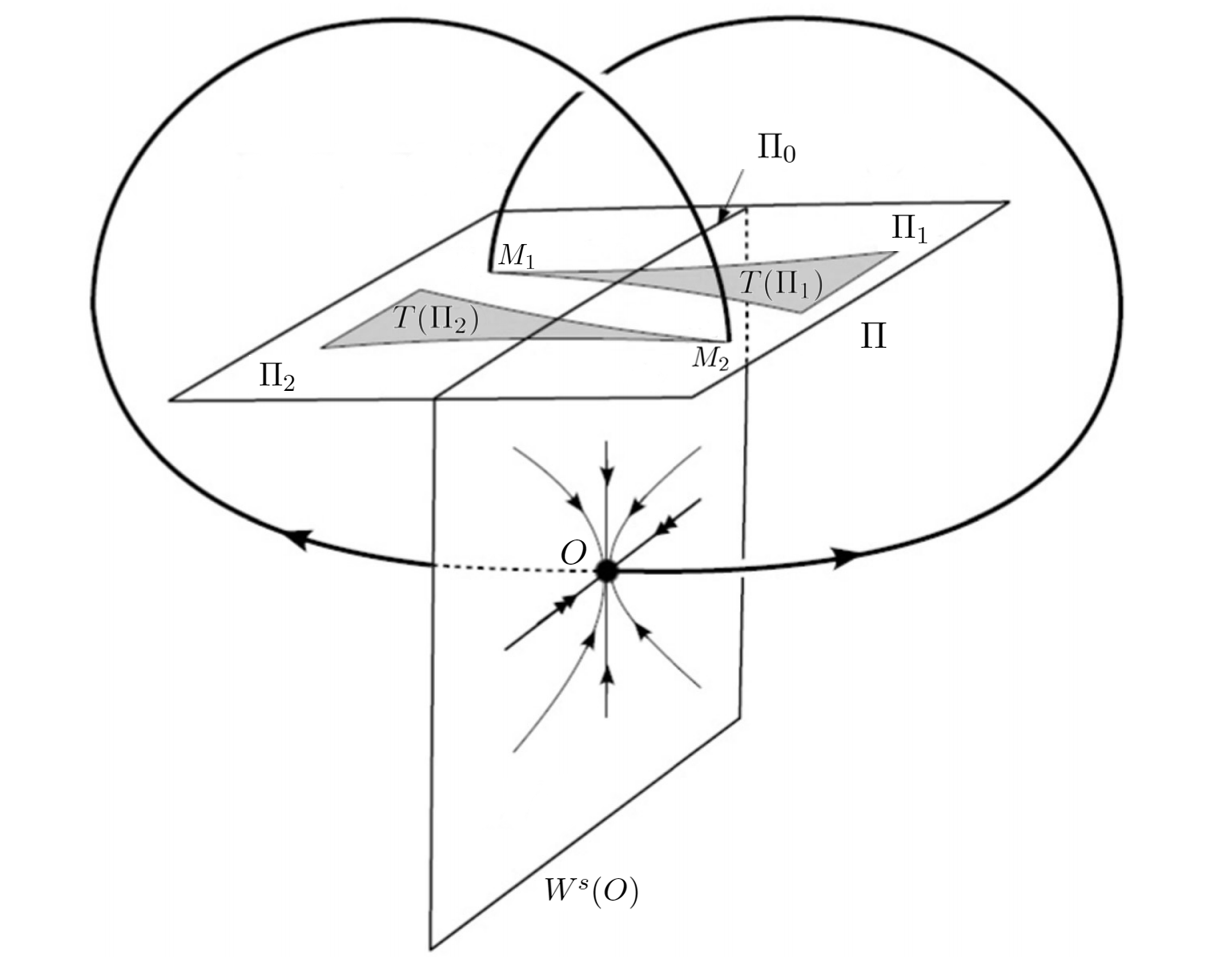}
\end{center}
\caption{The Afraimovich-Bykov-Shilnikov model.}
\label{fig:abs}
\end{figure}

Briefly, the ABS model can be described as follows. Let a smooth system of differential equations have a saddle equilibrium state $O$ with a one-dimensional unstable manifold $W^u(O)$. Assume also that the nearest to the imaginary axis characteristic exponent (an eigenvalue of the linearisation matrix) at $O$ is real and negative. Take a compact cross-section $\Pi$ (of codimension 1) transverse to a piece of the stable manifold $W^s(O)$, and let the two unstable separatrices $\Gamma_1$ and $\Gamma_2$ of $W^u(O)$ intersect $\Pi$ at some points $M_1$ and $M_2$, respectively. Denote by $\Pi_0$ the intersection of $\Pi$ with $W^s_{loc}(O)$, and by $\Pi_1$ and $\Pi_2$ the two parts separated by $\Pi_0$ so that we have $\Pi=\Pi_0\cup\Pi_1\cup\Pi_2$. Then, consider the Poincaré map $T$ on $\Pi$ induced by the orbits of the system - we assume that every orbit starting from $\Pi\backslash \Pi_0$ returns to $\Pi$, so the \Poincare map is defined everywhere on
$\Pi\backslash \Pi_0$ (the orbits that start on $\Pi_0$ tend to $O$ as $t\to+\infty$ and do not return to $\Pi$). 
Let $(u,v)$ be the coordinates on $\Pi$ such that $\{u=0\},\{u>0\}$ and $\{u<0\}$ correspond to $\Pi_0,\Pi_1$ and $\Pi_2$, respectively (see Fig. \ref{fig:abs}). 
The map $T$ is smooth outside $\Pi_0$, and for a point $M=(u,v)$ we have
\begin{equation*}
\lim_{u\to 0^+}T(M)=M_1 \quad \mbox{and}\quad \lim_{u\to 0^-}T(M)=M_2.
\end{equation*}

We assume that the image $T(\Pi)$ lies strictly in the inner part of $\Pi$, so a small neighbourhood $\cal D$ of the set formed by forward orbits starting from $\Pi$ is strictly forward-invariant, hence there is an attractor inside $\cal D$ (the Lorenz attractor). By the assumption on the characteristic exponents at $O$, the map $T$ near $\Pi_0$ is expanding in the $u$-direction and contracting in the $v$-direction. The main assumption of the ABS model is that this hyperbolicity property extends to the whole of $\Pi$. Under this assumption, there exists a smooth stable invariant foliation $\mathcal{F}$ on $\Pi$, which includes $\Pi_0$ as one of its leaves. Furthermore, the quotient map of $T$ obtained by taking quotient along the leaves of $\mathcal{F}$ is expansive. This allows for a detailed study of the structure of the attractor in $\cal D$ (see \citep{abs77,abs83} for details).

We will call the system Lorenz-like if it satisfies the above described properties of the ABS model. Note that
both models (\ref{eq:intro:0}) and (\ref{eq:intro:1}) are symmetric with respect to $(x,y)\to(-x,-y)$. In terms
of the ABS model, we will call it symmetric if the \Poincare map is symmetric with respect to an involution
that changes the sign of the expanding variable $u$.

Note that the equilibrium state $O$ is a saddle fixed point for the time-$t$ map of the system for any $t$. If we add a small $t$-periodic perturbation to a Lorenz-like system, then $O$ would continue
as a saddle fixed point of the time-$t$ map. Theorem 7 in \citep{ts08} states that for all small time-periodic perturbations of a Lorenz-like system the period map has a unique chain-transitive attractor $\mathcal{A}\subset {\cal D}$ which coincides with the set of all points attainable from $O$ by $\varepsilon$-orbits for all $\varepsilon>0$. In particular, the attractor $\mathcal{A}$ contains $O$ and its unstable manifold. Therefore, when $O$ is a part of the heterodimensional cycle, this heterodimensional cycle is in $\mathcal{A}$.

Recall that systems with homoclinic loops to $O$ are $C^\infty$-dense among Lorenz-like systems
\cite{abs77,ts08}; systems with a symmetric pair of homoclinic loops to $O$ are $C^\infty$-dense among symmetric Lorenz-like systems. For the time-$t$ map of the system (without a periodic perturbation), the homoclinic loop corresponds to a continuous family of orbits homoclinic to the fixed point $O$, i.e., to a non-transverse intersection of its stable and unstable manifolds. Thus, given any symmetric Lorenz-like system, we can add an arbitrarily small time-independent perturbation  (without destroying the symmetry) such that conditions (C1),(C2) will be satisfied. The strong-stable invariant foliation in the Lorenz-like systems \cite{abs77,abs83} also persists at small time-periodic perturbations \cite{ts08}, which implies that the non-degeneracy condition (C3) will hold automatically.

Thus, in order to apply Theorem \ref{thm1}, it remains to check condition (C4). The multipliers of $O$ for the time-1 map of an autonomous flow are the exponents of the eigenvalues of the linearisation matrix of the system at $O$. Therefore, condition (C4) will be fulfilled by the time-$t$ map of a Lorenz-like flow (and, hence, by any sufficiently small perturbation of it) if

(C4$^\prime$) ${\rm Re} \; \nu_1<2\nu_0$ and $\nu_0 +\dfrac{1}{2}\nu <0$,

\noindent where $\nu_j$ and $\nu$ are the characteristic exponents of $O$ such that
$$\dots \leq {\rm Re}\; \nu_2 \leq {\rm Re}\; \nu_1 < \nu_0 < 0 < \nu.$$

We arrive at the following

\begin{thm}\label{thm3}
Let the equilibrium state of a symmetric Lorenz-like system satisfy condition ({\rm C4}$^\prime$). Then, there exists an arbitrarily small time-periodic perturbation (which keeps the symmetry of the system) such that the attractor $\mathcal{A}$ of the period map of the perturbed system contains a symmetric pair of heterodimensional cycles, each of which involves $O$
and an index-2 saddle periodic point. Moreover, in an open neighbourhood of this map in 
$\diff ^r_s(\mathcal{D})$, these heterodimensional cycles are a part of the attractor $\mathcal{A}$ for a $C^r$-dense subset
of this neighbourhood (for any $r\leq \infty$).
\end{thm}
Note that the $C^\omega$ case is not included here because we do not know whether the perturbation for a Lorenz-like system to have a pair of homoclinic tangencies without destroying the symmetry can be made analytic. If condition (C4$^\prime$) is not fulfilled, then a weaker statement follows from Theorem \ref{thm2}.
\begin{thm}\label{thm4}
For any symmetric Lorenz-like system, there exists an arbitrarily small (in $C^r$, for any $r\leq \infty$) time-periodic perturbation such that the attractor $\mathcal{A}$ of the period map of the perturbed system contains a heterodimensional cycle involving $O$ and an index-2 saddle periodic point.
\end{thm}

Note that the Lorenz system (\ref{eq:intro:0}) does not satisfy condition (C4$^\prime$) at classical parameter values,
while the Morioka-Shimizu system (\ref{eq:intro:1}) fulfils this condition for the set of parameter values for which a proof of the existence of Lorenz attractor is obtained in \citep{ctz17}. Therefore, Theorem
\ref{thm4} is applicable to time-periodic perturbations of the Lorenz attractor in the Lorenz system, and the stronger Theorem \ref{thm3} is applicable to the periodic perturbation of the  Lorenz attractor in the Morioka-Shimizu system.

The rest of this paper is organised as follows. In Section \ref{sec:returnmap} we describe the dynamics near $O$ and define the first return map. In Section \ref{sec:preperturb} we make perturbations which give us a homoclinic tangency with some special properties required to create heterodimensional cycles. Next, we give in Section \ref{sec:index2} the condition for having a periodic point of index 2. A formula for leaves of the strong-stable foliation $\mathcal{F}^s$ is derived in Section \ref{sec:cones}. Finally, with all the preparation, we prove Theorems \ref{thm1} and \ref{thm2} in Section \ref{sec:proofofthm}.

\section{The first return map}\label{sec:returnmap}

Let a $C^r$-diffeomorphism $F$ fulfil conditions (C1)-(C3). We embed it into a parametric family $F_{\varepsilon}$ such that $F=F_{\varepsilon^*}$, where $\varepsilon$ is the set of parameters defined in the previous Section. Observe that this family is transverse to the surface of diffeomorphisms satisfying (C1)-(C3). 

Let $V$ be a small neighbourhood of $O$, and take two points $M^+,M^-\in\Gamma\cap V$ such that $M^+\in W^s_{loc}(O)$, $M^-\in W^u_{loc}(O)$, $F^{-\tau}(M^+)\notin V$ and $F^\tau(M^-)\notin V$, where $\tau$ is the period of the point $O$. Let $\Pi_0,\Pi_1\subset V$ be two small open sets containing $M^+$ and $M^-$, respectively. In what follows we consider the local map $T_0\equiv F_{\varepsilon}^{\tau}|_{V}: V\to \mathcal{M}$ and the global map $T_1\equiv F_{\varepsilon}^l|_{\Pi_1}:$ $\Pi_1\to \mathcal{M}$ where $l$ is the positive integer such that $F^l(M^-)=M^+$ (it exists, because $M^+$ and $M^-$ belong to the same orbit $\Gamma$).

Let $C^r$-coordinates $(x,y,z)\in \mathbb{R}^D$ be introduced in $V$ such that the map $T_0$ takes the form
\begin{equation}\label{eq:setting:T_0}
\begin{array}{rcl}
\bar{x} &=& \lambda(\varepsilon) x +  f_1(x,y,z,\varepsilon), \\
\bar{y} &=& \gamma(\varepsilon) y +f_2(x,y,z,\varepsilon),\\
\bar{z} &=& A(\varepsilon)z + f_3(x,y,z,\varepsilon),
\end{array}
\end{equation}
where the eigenvalues of the $(D-2)\times (D-2)$ matrix $A$ are the multipliers $\lambda_1\dots\lambda_{D-2}$; the functions $f_i$ ($i=1,2,3$) and their first derivatives vanish at the origin, and, furthermore,
\begin{equation}\label{eq:nonlinearterms}
\begin{array}{l}
f_{1,3}(0,y,0,\varepsilon)=0,\quad\quad f_2(x,0,z,\varepsilon)=0, \quad\quad f_1(x,0,z,\varepsilon)=0,  \quad\quad  f_2(0,y,0,\varepsilon)=0,\\[10pt]
\dfrac{\partial f_{1,3}}{\partial (x,y)}(0,y,0,\varepsilon)=0, \quad \quad \dfrac{\partial f_2}{\partial y}(x,0,z,\varepsilon)=0
\end{array}
\end{equation}
for all sufficiently small $x$, $y$ and $z$. The existence of such coordinate transformation is shown in \cite{gst08}. In the appendix we show that in the symmetric case (i.e., when $F\in \diff^r_s$) this transformation can be done in such a way that the involution $\mathcal{R}$ is still locally linear and satisfies \eqref{eq:symmetryR} in the new coordinates. Note that this coordinate transformation, and its first and second derivatives with respect to $(x,y,z)$, are $C^{r-2}$-smooth functions of both the parameters $\varepsilon$ and $(x,y,z)$ \cite{gst08}. Therefore, $\lambda$, $\gamma$, and $A$
in \eqref{eq:setting:T_0} are $C^{r-2}$-smooth functions of $\varepsilon$, and the functions $f_{1,2,3}$, as well as the derivatives of $f_{1,2,3}$ with respect to $(x,y,z)$  up to order 2, are $C^{r-2}$-smooth functions of $(x,y,z,\varepsilon)$.

The first two identities in \eqref{eq:nonlinearterms} mean that the local manifolds $W^s_{loc}(O)$ and $W^u_{loc}(O)$ are straightened, i.e., we have $W^s_{loc}(O)=\{y=0\}$ and $W^u_{loc}(O)=\{x=0,z=0\}$. The third identity implies that the leaves of the strong-stable foliation $\mathcal{F}_0$ in $W^s_{loc}(O)$ have the form $\{x=c,y=0\}$ and the quotient map on $W^s_{loc}(O)$ obtained by factorising over the leaves of $\mathcal{F}_0$ is linear. The forth identity corresponds to the linearisation of the map restricted to $W^u_{loc}(O): \{x=0,z=0\}$. 

In order to obtain necessary formulas for the first return map to $\Pi_0$, we need, first, to consider iterates of $T_0$. Take any point $(x_0,y_0,z_0)\in V$, and let $(x_k,y_k,z_k)=T^k_0(x_0,y_0,z_0)$. The triple $(x_k,y_0,z_k)$ is a uniquely defined function of $x_0,y_k$ and $z_0$ on a small neighbourhood of ($x^+,y^-,z^+$) for any $k\geqslant 0$ (see e.g. \citep{sh67,gs90}). It follows from Lemma 7 of \cite{gst08} that if the map $T_0$ satisfies conditions (\ref{eq:nonlinearterms}), then the following relations hold for all sufficiently large $k$:
\begin{equation}\label{eq:setting:T^k}
\begin{array}{rcl}
x_k &=& \lambda(\varepsilon)^k x_0 +  \phi_k(x_0,y_k,z_0,\varepsilon), \\
y_0 &=& \gamma(\varepsilon)^{-k} y_k + \psi_k(x_0,y_k,z_0,\varepsilon), \\
z_k &=& \hat{\phi}_k(x_0,y_k,z_0,\varepsilon),
\end{array}
\end{equation}
where $\phi_k, \psi_k, \hat\psi_k$ are smooth functions such that
\begin{equation}\label{fpnorm}
\|\phi_k, \hat\phi_k \|_{_2} = o(|\lambda(\varepsilon)|^k), \qquad \|\psi_k\|_{_2}=o(|\gamma(\varepsilon)|^{-k}),
\end{equation}
and also
\begin{equation}\label{fpnorm1}
\|\hat\phi_k \|_{_1} = o(\hat\lambda^k)
\end{equation}
where $\hat{\lambda}$ is any number such that 
$\max\{\lambda^2,|\lambda_1|\}<\hat{\lambda}<|\lambda|$. We use the following notation in formulas \eqref{fpnorm} and \eqref{fpnorm1}: $\|\cdot\|_1$ stands for the maximum of the $C^0$-norms of the function and its first derivative with respect to $(x_0,y_k,z_0)$, while $\|\cdot\|_2$ denotes  the maximum of the $C^0$-norms of the function, its first derivative with respect to $(x_0,y_k,z_k,\varepsilon)$, and all its second derivatives except for the second derivative with respect to $\varepsilon$ alone.

In the case where condition (C4) is fulfilled, we obtain stronger estimates. In Appendix \ref{app:t4} we show that when $|\lambda_1|<\lambda^2$ and
$|\lambda\gamma|>1$ there exists a $C^2$-smooth extended unstable invariant manifold $W_{loc}^{uE}(O)$ which contains the local unstable manifold $W_{loc}^{u}(O)$ and is tangent to $z=0$ at the points of $W_{loc}^{u}(O)$, i.e., $W_{loc}^{uE}(O)$ is given by the equation 
$z=\eta(x,y,\varepsilon)$ where $\eta(0,y,\varepsilon)\equiv 0$, $\frac{\partial}{\partial x} \eta(0,y,\varepsilon)\equiv 0$. 
Furthermore, in $W_{loc}^{uE}(O)$ there is an invariant foliation $\mathcal{F}^{uE}$ with the leaves of the form $h(x,y,\varepsilon)=const$ where
$h(x,0,\varepsilon)\equiv x$ and $h(0,y,\varepsilon)\equiv 0$. The functions $\eta$ and $h$ are $C^2$, but if the coordinates are introduced where
the map $T_0$ gets into the form \eqref{eq:setting:T_0},\eqref{eq:nonlinearterms}, the second derivative with respect to $\varepsilon$ alone may not exist. It is also shown in the Appendix that in the symmetric case the manifold $W_{loc}^{uE}(O)$ and the invariant foliation $\mathcal{F}^{uE}$ on it are invariant with respect to the involution $\mathcal{R}$, i.e., $\eta(x,-y,\varepsilon) \equiv \mathcal{S}\eta(x,y,\varepsilon)$ and $h(x,-y,\varepsilon)\equiv h(x,y,\varepsilon)$. From now on, we will omit $\varepsilon$ in all expressions for simplicity.  

We can now choose new coordinates $z^{new}=z-\eta(x,y)$ and $x^{new}=h(x,y)$. It is easy to see that the map keeps its form
\eqref{eq:setting:T_0},\eqref{eq:nonlinearterms} in the new coordinates, and estimates \eqref{eq:setting:T^k},\eqref{fpnorm} and \eqref{fpnorm1} hold.
In the symmetric case, we also have that formula \eqref{eq:symmetryR} for the involution $\mathcal{R}$ remains unchanged.

In the new coordinates the invariant manifold $W_{loc}^{uE}(O)$ and foliation $\mathcal{F}^{uE}$ get straightened: 
$W_{loc}^{uE}(O)$ is given by $\{z=0\}$ and the leaves of $\mathcal{F}^{uE}$ are $\{x=const,z=0\}$. This implies that in the new coordinates
\begin{equation}\label{eq:nonlinearterms1}
f_3(x,y,0)=0, \quad\quad  f_1(x,y,0)=0
\end{equation}
(the first equation follows from the invariance of $W_{loc}^{uE}(O)$; the invariance of $\mathcal{F}^{uE}$ implies that $f_1(x,y,0)=f_1(x,0,0)$, which
gives the second equation of \eqref{eq:nonlinearterms1} by virtue of the third equation of \eqref{eq:nonlinearterms}).

\begin{lem}\label{lem:firstd} Once identities \eqref{eq:nonlinearterms} and \eqref{eq:nonlinearterms1} are fulfilled, one can find positive constant 
$\lambda_0<\lambda^2$ such that, for all $k\geq 0$,
\begin{equation}\label{eq:firstd}
\left\|\dfrac{\partial x_k}{\partial z_0}\right\|\leqslant \lambda_0^k,
\quad\quad
\left\|\dfrac{\partial z_k}{\partial z_0}\right\|\leqslant \lambda_0^k.
\end{equation}
\end{lem}

\noindent{\it Proof.}
We can rewrite formula \eqref{eq:setting:T_0} for $T_0$ as
\begin{equation*}
\begin{array}{rcl}
\bar{x} &=& \lambda x +  f_1(x,y,z), \\
y &=& \gamma^{-1} \bar{y} -\gamma^{-1}f_2(x,y,z),\\
\bar{z} &=& Az + f_3(x,y,z),
\end{array}
\end{equation*}
from which one deduces the following relation between $(x_0,y_k,z_0)$ and its $j$-th iterate $(x_j,y_j,z_j)$ ($1\leqslant j \leqslant k$):
\begin{equation}\label{eq:ssleaves:9}
\begin{array}{rcl}
x_j&=&\lambda^j x_0+\sum^j_{s=1}\lambda^{s-1}f_1(x_{j-s},y_{j-s},z_{j-s}),\\[5pt]
y_j&=&\gamma^{j-k} y_k - \sum^k_{s=j+1}\gamma^{-s+j} f_2(x_{k-s+j},y_{k-s+j},z_{k-s+j}),\\[5pt]
z_j&=&A^j z_0+\sum^j_{s=1}A^{s-1}f_3(x_{j-s},y_{j-s},z_{j-s}).
\end{array}
\end{equation}
By formulas A.18, A.20 and A.34 in \citep{gst08}, we have
\begin{equation}\label{eq:ssleaves:10}
\|y_j\|\leqslant C|y_k|\cdot |\gamma|^{j-k},\quad 
\left\|\dfrac{\partial y_j}{\partial z_0}\right\|\leqslant C |\gamma|^{j-k}
\end{equation}
for some constant $C$.
Since $f_3$ vanishes at $z=0$ (see \eqref{eq:nonlinearterms1}), and its derivative vanishes at the origin, it follows that
$$\|f_3\|\leq \delta\|z\|$$
where $\delta$ can be made as small as we need by taking the neighbourhood $V$ of the otigin sufficiently small. Therefore,
$$\|\bar z\|\leq (\|A\|+\delta) \|z\| \leq \lambda_0 \|z\|$$
(we can always choose such $\lambda_0$ satisfying $\lambda_0<\lambda^2$ because $|\lambda_1|<\lambda^2$ by the assumption of this lemma). This gives
\begin{equation}\label{eq:ssleaves:10a}
\|z_j\|\leqslant \|z_0\| \lambda_0^j.
\end{equation}

Now assume  that the inequalities 
\begin{equation}\label{eq:ssleaves:11}
\left\|\dfrac{\partial (x_s,z_s)}{\partial z_0}\right\|\leqslant \lambda_0^s
\end{equation}
hold for all $s=0,\dots, j-1$ (they are, obviously true for $s=0$) and prove that they remain true for $s=j$. By induction, this will prove the lemma.

By differentiating equations \eqref{eq:ssleaves:9}, we find
\begin{equation}\label{eq:ssleaves:12}
\begin{array}{rcl}
\dfrac{\partial x_j}{\partial z_0} &=&
\sum^{j}_{s=1}\lambda^{s-1}\Big(
\dfrac{\partial f_1}{\partial x}\dfrac{\partial x_{j-s}}{\partial z_0} 
+ \dfrac{\partial f_1}{\partial y}\dfrac{\partial y_{j-s}}{\partial z_0}
+\dfrac{\partial f_1}{\partial z}\dfrac{\partial z_{j-s}}{\partial z_0}  
\Big),\\[5pt]
\dfrac{\partial z_{j}}{\partial z_0} &=& 
A^j + \sum^{j}_{s=1}A^{s-1}\Big(
\dfrac{\partial f_3}{\partial x}\dfrac{\partial x_{j-s}}{\partial z_0} 
+ \dfrac{\partial f_3}{\partial y}\dfrac{\partial y_{j-s}}{\partial z_0}
+\dfrac{\partial f_3}{\partial z}\dfrac{\partial z_{j-s}}{\partial z_0}  
\Big).
\end{array}
\end{equation}
Recall that the $C^2$ function $f_1$ vanishes both at $z=0$ and $y=0$ while the $C^2$ function $f_3$ vanishes at $z=0$ (see \eqref{eq:nonlinearterms},\eqref{eq:nonlinearterms1}) and its derivative is zero at the origin. Therefore,
$$\left\|\dfrac{\partial f_1}{\partial (x,z)}\right\|\leq K \|y\|, \qquad \left\|\dfrac{\partial f_1}{\partial y}\right\|\leq K \|z\|,$$
$$\left\|\dfrac{\partial f_3}{\partial (x,y)}\right\|\leq K \|z\|, \qquad \left\|\dfrac{\partial f_3}{\partial z}\right\|\leq \delta,$$
where $K$ and $\delta$ are some constants and $\delta$ can be chosen as small as we want (by choosing the neighbourhood $V$ small enough).
By plugging these inequalities into \eqref{eq:ssleaves:12}, we obtain
$$
\begin{array}{rcl}
\left\|\dfrac{\partial x_{j}}{\partial z_0}\right\| &\leq &
K \sum^{j}_{s=1}|\lambda|^{s-1}\Big(\|y_{j-s}\| \cdot \left\|\dfrac{\partial (x_{j-s},z_{j-s})}{\partial z_0}\right\|+
\|z_{j-s}\| \cdot \left\|\dfrac{\partial y_{j-s}}{\partial z_0}\right\|\Big),\\[10pt]
\left\|\dfrac{\partial z_{j}}{\partial z_0}\right\| &\leq &
\|A\|^j+\sum^{j}_{s=1}\|A\|^{s-1}\Big(K\|z_{j-s}\|\cdot\left\|\dfrac{\partial x_{j-s}}{\partial z_0}\right\|
+ K\|z_{j-s}\| \cdot\left\|\dfrac{\partial y_{j-s}}{\partial z_0}\right\|
+\delta \left\|\dfrac{\partial z_{j-s}}{\partial z_0}\right\|\Big). 
\end{array}
$$
Now, using estimates \eqref{eq:ssleaves:10},\eqref{eq:ssleaves:10a} (where one should replace $j$ by $(j-s)$) and \eqref{eq:ssleaves:11} (where one should change $s$ to $(j-s)$), we obtain
\begin{equation}\label{eq:ssleaves:14}
\begin{array}{rcl}
\left\|\dfrac{\partial x_{j}}{\partial z_0}\right\| &\leq &
K \sum^{j}_{s=1}|\lambda|^{s-1}\Big(C |y_k| \cdot |\gamma|^{j-s-k} \cdot \lambda_0^{j-s} +
\|z_0\|\lambda_0^{j-s} \cdot C |\gamma|^{j-s-k}\|\Big)\leq\\
 &\leq &  \dfrac{KC}{|\lambda|} (|y_k| + |z_0|) \;\lambda_0^j \;\sum^{j}_{s=1} \left(\dfrac{|\lambda|}{|\gamma| \lambda_0}\right)^s,\\~\\
\left\|\dfrac{\partial z_{j}}{\partial z_0}\right\| &\leq &
\|A\|^j+\sum^{j}_{s=1}\|A\|^{s-1}\Big(K\|z_0\|\lambda_0^{j-s} \cdot \lambda_0^{j-s}
+ K\|z_0\|\lambda_0^{j-s} \cdot C |\gamma|^{j-s-k}
+\delta \lambda_0^{j-s}\Big)\leq \\
 &\leq & \lambda_0^j+ \dfrac{K \|z_0\| (C+1) +\delta }{\|A\|} \;\lambda_0^j\; \sum^{j}_{s=1} \left(\dfrac{\|A\|}{\lambda_0}\right)^s.
\end{array}
\end{equation}
Recall that we assume $|\lambda\gamma|>1$. Hence, if $\lambda_0<\lambda^2$ is chosen close enough to $\lambda^2$, we have
$\frac{|\lambda|}{|\gamma| \lambda_0}<1$. Also, since $|\lambda_1|<\lambda^2$, where $\lambda_1$ is the largest, in the absolute value, eigenvalue of $A$, we have that $\lambda_0<\lambda^2$ can be chosen such that $\frac{\|A\|}{\lambda_0}<1$. This means that the sums
$\sum^{j}_{s=1} \left(\frac{|\lambda|}{|\gamma| \lambda_0}\right)^s$ and $\sum^{j}_{s=1} \left(\frac{\|A\|}{\lambda_0}\right)^s$
in \eqref{eq:ssleaves:14} are uniformly bounded for all $j$. Therefore, since  $|y_k|$, $\|z_0\|$ and $\delta$ 
can be taken as small as we need by choosing the neighbourhood $V$ small enough, the estimates \eqref{eq:ssleaves:14}
imply that the inequalities \eqref{eq:ssleaves:11} hold for $s=j$. Therefore, by induction, they hold for all $s$. At $s=k$ we obtain the lemma.
\qed

We now proceed to obtain necessary formulas for the global map $T_1$. Let us write its Taylor expansion near the point $M^-$. 
At $\varepsilon=\varepsilon^*$, the point $M^-$ is homoclinic, so its image $M^+=T_1 M^-$ belongs to the local stable manifold,
and the curve $T_1 W^u_{loc}$ has a quadratic tangency to $W^s_{loc}$. In the coordinate system where the local stable and 
unstable manifolds are straightened, i.e., they are given by the equations $\{y=0\}$ and, respectively, $\{x=0,z=0\}$, we have $M^-=(0,y^-,0)$
and $M^+=(x^+,0,z^+)$ and the Taylor expansion for $T_1: (x,y,z) \mapsto (x',y',z')$ is given by 
\begin{equation}\label{eq:setting:T_1}
\begin{array}{rcl}
x'-x^+ &=& ax+b(y - y^-) + a_{13} z +h_1(x,y-y^-,z), \\ [5pt]
y' &=& y^+(\varepsilon) + cx + d(y-y^-)^2 + a_{23} z + h_2(x,y-y^-,z),\\[5pt]
z'-z^+ &=& a_{31}x+a_{32}(y - y^-) + a_{33} z +h_3(x,y-y^-,z), \\ [5pt]
\end{array}
\end{equation}
where $d\neq 0$ and the Taylor expansions for functions $h_{1,2,3}$ start with quadratic terms (the term $d(y-y^-)^2$ is taken out of $h_2$,
so $h_2$ does not contain it). We will use the coordinate system where the map $T_0$ is in the form \eqref{eq:setting:T_0} and the identities
\eqref{eq:nonlinearterms} hold. 

When we vary $\varepsilon$, the map $T_1$ can be kept in the form \eqref{eq:setting:T_1} where the coefficients and the functions $h_{1,2,3}$ now depend on $\varepsilon$ (e.g. we choose $y^-(\varepsilon)$ in such a way that there is no linear term in $(y-y^-(\varepsilon))$ in the equation for $y'$ in \eqref{eq:setting:T_1}). We however take $d$ independent of $\varepsilon$, so $h_2$ is allowed to include the $(y-y^-(\varepsilon))^2$-term with the coefficient which vanishes at $\varepsilon=\varepsilon^*$. Recall that the coordinates we use are of class $C^2$, but the second derivative with respect to $\varepsilon$ alone may not exist. Thus, we have that all the coefficients, as well as the functions $h_{1,2,3}$ and their first derivatives with respect to $(x,y,z)$ are at least $C^1$ functions of $\varepsilon$. So, we can write
\begin{equation}\label{eq:setting:3}
h_{1,3}=O(x^2+(y-y^-)^2+z^2),
\qquad
h_2=O(x^2+z^2+|x| \cdot |y-y^-| +\|z\| \cdot |y-y^-|) + o((y-y^-)^2)_{\varepsilon\to \varepsilon^*},
\end{equation}
and
\begin{equation}\label{eq:setting:3.1}
\frac{\partial h_{1,2,3}}{\partial \varepsilon}=o(|x|+\|z\|+|y-y^-|), \qquad 
\frac{\partial^2 h_{1,2,3}}{\partial \varepsilon\partial(x,y,z)}=o(1)_{(x,y-y^-(\varepsilon),z)\to 0}.
\end{equation}

By construction, the value of $y^+(\varepsilon)$ measures the magnitude of splitting between the curve $T_1 W^u_{loc}$ and the local stable manifold. Thus, $\mu(F_\varepsilon)=y^+(\varepsilon)$ can be taken as the parameter governing the splitting of the homoclinic tangency at the point $M^+$. It is our standing condition that $\partial\mu/\partial\varepsilon\neq 0$, so we simply assume that $\mu$ is one of the parameters $\varepsilon$ (see the explanation before Theorem 1).

Note that our conditions in Section \ref{sec:intro} imply that
\begin{equation}\label{eq:nc}
d\neq0,\,x^+\neq 0\,\mbox{ and  }\, bc\neq0
\end{equation}
in formula (\ref{eq:setting:T_1}). The first two inequalities come, respectively, from the facts that the tangency is quadratic and it is not in the strong-stable manifold of $O$. The third one follows from the transversality of the extended unstable manifold $W^{uE}(O)$ to the strong-stable foliation $\cal F_0$,
see Condition (C3). 

Indeed, the first identity in the second line of (\ref{eq:nonlinearterms}) implies that $W^{uE}_{loc}$ is tangent to the plane $z=0$ at the points of $W^{u}_{loc}$ (see \cite{gst08}); in particular, it is tangent to
$z=0$ at the homoclinic point $M^-$. So, the tangent plane to the image $T_1 W^{uE}_{loc}$ is given by
$$x'-x^+ = a x + b (y-y^-), \qquad y'= c x, \qquad z'=a_{31}x+a_{32}(y-y^-).$$
The transversality of $T_1W^{uE}$ to $\mathcal{F}_0$ just means that this tangent plane intersects the strong-stable leaf $\{x'=x^+,y'=0\}$ at a single point (the point $M^+$). This is equivalent to the requirement that the equation
$$0= ax + b(y-y^-), \qquad 0 = cx$$
has only one solution ($x=0,y=y^-$), which implies $bc\neq 0$.

\begin{figure}[!h]
\begin{center}
\includegraphics[width=0.40\textwidth]{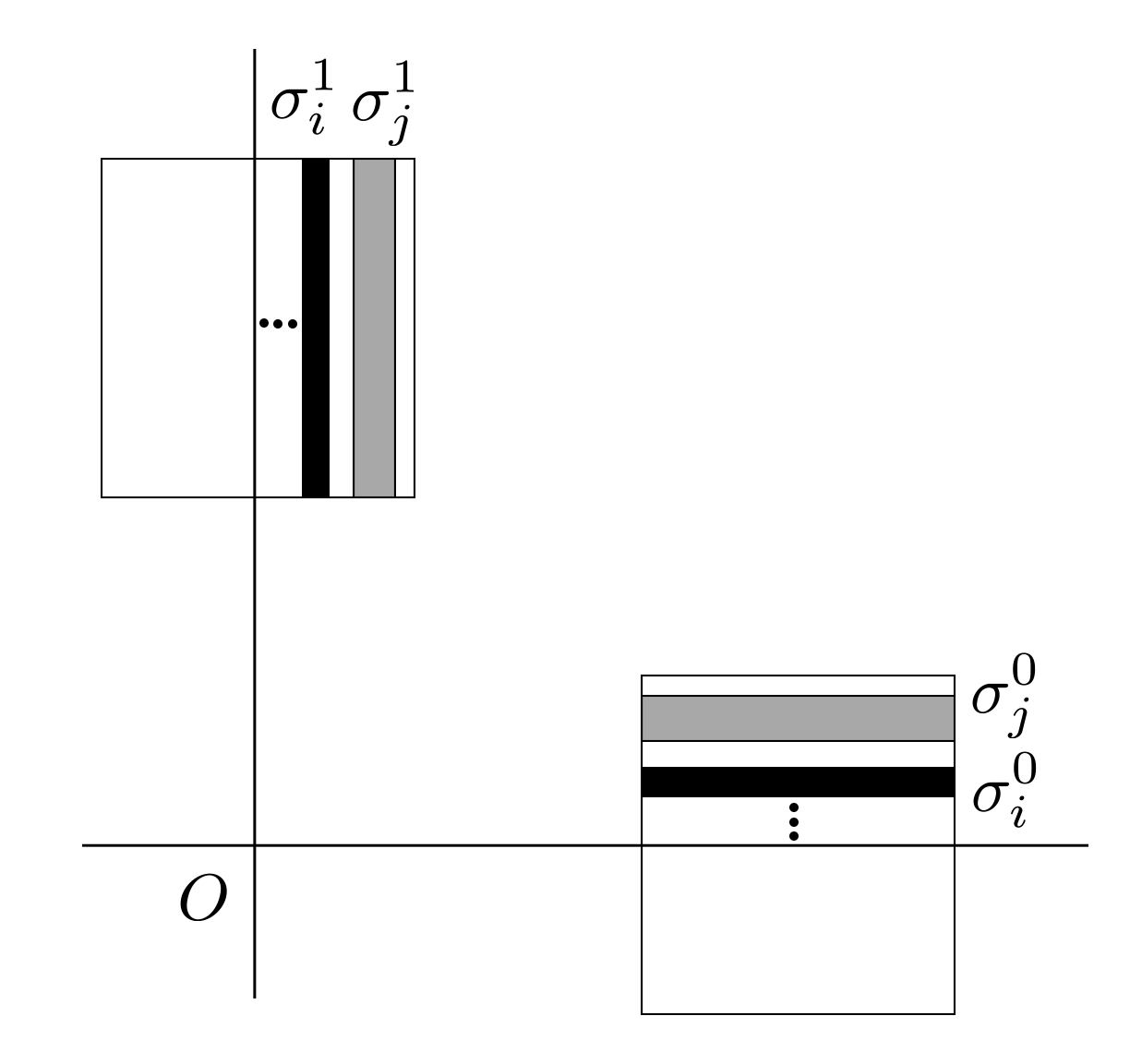}
\end{center}
\caption{The projections of the countable sequences of disjoint sets $\sigma_k^0$ along the leaves of $\mathcal{F}^s$ onto $\{z=0\}$.}
\label{fig:T_0}
\end{figure}

We can now define the maps $T_1T_0^k$ of the first return to $\Pi_0$. We fix the choice of the neighbourhoods $\Pi_0$ and $\Pi_1$ as follows: $\Pi_0=\{(x,y,z)\mid |x-x^+|<\delta/2,|y|<\delta,\|z-z^+\|<\delta/2\}$ and $\Pi_1=\{(x,y,z)\mid |x|<\delta,|y-y^-|<\delta/2,\|z\|<\delta\}$, where $\delta>0$ is small such that $T_0(\Pi_0)\cap \Pi_0=\emptyset$ and $T_0^{-1}(\Pi_1)\cap \Pi_1=\emptyset$. Let $k^*$ be the smallest number such that $T_0(\Pi_0)\cap\Pi_1\neq\emptyset$. There are two countable sequences of disjoint subsets $\sigma^0_k\subset \Pi_0$ and $\sigma^1_k:=T^k_0(\sigma^0_k)\subset\Pi_1$ such that $k\geqslant k^*$, and $\sigma^0_k\to W_{loc}^s(O)$ and $\sigma^1_k\to W^u_{loc}(O)$ as $k\to +\infty$ (see Figure \ref{fig:T_0}). Therefore, the first-return map $T:\Sigma^0:=\bigcup^{+\infty}_{k_0}\sigma^0_k \to \Pi_0$ is defined as
\begin{equation}\label{eq:setting:T}
T(M)=T_1\circ T_0^k(M) \quad \mbox{if} \quad M\in\sigma^0_k.
\end{equation}
For a point $M\in\Sigma^0$ we call the corresponding $k$ in \eqref{eq:setting:T} {\it the stay number} of $M$. 
The image of $\Sigma^0$ under $T$ may not be entirely contained in $\Pi_0$. However, throughout this paper, we only consider points sufficiently close to $M^+$ such that their images lie in $\Pi_0$.

In the same way, a global map $\tilde T_1$ and a first-return map $\tilde T$ are defined near the second orbit of homoclinic tangency, $\tilde\Gamma$. In the symmetric case, i.e., when $F_\varepsilon\in\diff ^r_s(\mathcal{M})$, the maps $T_1$ and $\tilde T_1$ are related by the symmetry $\mathcal{R}$. Namely, we denote by $\tilde{M}^+$ and $\tilde{M}^-$ the points that are $\mathcal{R}$-symmetric to $M^+$ and $M^-$. These two points satisfy $\tilde{M}^+\in W^s_{loc}(O)\cap\tilde{\Gamma}$, $\tilde{M}^-\in W^u_{loc}(O)\cap\tilde{\Gamma}$, and have coordinates $(x^+,0,\mathcal{S} z^+)$ and $(0,-y^-,0)$. We can choose the neighbourhood $\Pi_0$ such that it will contain both points $M^+$ and $\tilde{M}^+$. In order to achieve this, note that the directions corresponding to coordinates $z$ are strongly contracting, so we can just let $\Pi_0$ be the set $\{(x,y,z)\mid |x-x^+|<\delta/2,|y|<\delta,\|z\|<\delta\}$ and choose $x^+$ sufficiently small. When $\delta$ is small, the property $T_0(\Pi_0)\cap \Pi_0=\emptyset$ and $T_0^{-1}(\Pi_1)\cap \Pi_1=\emptyset$ holds. The neighbourhood $\tilde\Pi_1$ is defined as $\tilde{\Pi}_1= {\mathcal R}\Pi_1=\{(x,y,z)\mid |x|<\delta,|y+y^-|<\delta/2,\|\mathcal{S}z\|<\delta\}$, which implies $T_0^{-1}(\tilde{\Pi}_1)\cap\tilde{\Pi}_1=\emptyset$.

The second global map $\tilde{T}_1\equiv F^l|_{\tilde{\Pi}_1}:(x,y,z)\mapsto(x',y',z')$ takes the form 
\begin{equation}\label{eq:setting:T_11}
\begin{array}{rcl}
x'-x^+ &=& ax-b(y + y^-) + a_{13}\mathcal{S}  z +h_1(x,-y-y^-,\mathcal{S} z), \\ [5pt]
y' &=& -\mu - cx - d(y+y^-)^2  - a_{23}\mathcal{S}  z - h_2(x,-y-y^-,\mathcal{S} z),\\[5pt]
z'-\mathcal{S}z^+ &=& \mathcal{S}a_{31}x-\mathcal{S}a_{32}(y + y^-) + a_{33} z +\mathcal{S} h_3(x,-y-y^-,\mathcal{S} z), \\ [5pt]
\end{array}
\end{equation}
with the same coefficients and functions $h_{1,2,3}$ as in \eqref{eq:setting:T_1}.

There is a countable sequence of disjoint subsets $\tilde{\sigma}^0_k\subset \Pi_0$ such that $\tilde{\sigma}^1_k=T^k_0(\tilde{\sigma}^0_k)\subset\tilde{\Pi}_1$, where $k\geqslant k^*$, and $\tilde{\sigma}^0_k\to W_{loc}^s(O)$ and $\tilde{\sigma}^1_k\to W^u_{loc}(O)$ as $k\to +\infty$. The first return map $\tilde{T}:\tilde{\Sigma}^0=\bigcup^{+\infty}_{k_0}\tilde{\sigma}^0_k \to \tilde{\Pi}_0$ is defined as
\begin{equation}\label{eq:setting:T'}
\tilde{T}(M)=\tilde{T}_1\circ T_0^k(M) \quad \mbox{if} \quad M\in\tilde{\Sigma}^0_k.
\end{equation}

\section{An adjustment to the homoclinic tangency}\label{sec:preperturb}
In order to create a heterodimensional cycle in the small neighbourhood $U$ of $O\cup\Gamma\cup\tilde{\Gamma}$, we need the homoclinic tangency to satisfy the following conditions:\\[5pt]
\indent (a) the signs of $cdx^+$ and $cx^+y^-$ are positive, where $c$ and $d$ are the coefficients in the global map \eqref{eq:setting:T_1}; and\\[5pt]
\indent (b) there are two transverse homoclinic points in $W^u_{loc}(O)$ close to $M^-$ such that $M^-$ lies between these two points.\\[5pt]

In Section \ref{sec:proofofthm1}, conditions (a) and (b) are used to show the existence of the non-transverse and, respectively, transverse intersections between the invariant manifolds of two periodic orbits of different indices. In this section we prove that unfolding the original homoclinic tangency produces new homoclinic tangencies satisfying the above conditions. Depending on the signs of $c$ and $d$, the original homoclinic tangency falls into one of the four classes: (1) $cdx^+<0,dy^-<0$, (2) $cdx^+<0,dy^->0$, (3) $cdx^+>0,dy^-<0$, and (4) $cdx^+>0,dy^->0$. We start with showing that tangencies of classes (1), (3), and (4) can be replaced by tangencies of class (2).

\begin{lem}\label{lem:preperturb1}
Take any smooth one-parameter family $F_\mu$ of diffeomorphisms, where $\mu$ is the splitting parameter for the homoclinic tangency $\Gamma$, and $F_0$ fulfils conditions (C1)-(C3). Then, there exists a sequence $\{\mu_k\}$ accumulating on $\mu=0$ such that the saddle $O$ of $F_{\mu_k}$ has a class (2) homoclinic tangency and a tangency point $M^-_k\in W^u_{loc}(O)$ satisfying $M^-_k\to M^-$ as $k\to +\infty$.
\end{lem}

\noindent{\it Proof.} We will assume $x^+>0$ and $y^->0$ throughout this section since this can be always achieved by changing signs of $x$ and/or $y$ at the very beginning. There is nothing to prove if the original tangency already belongs to class (2). For the remaining three cases, we first construct new tangencies, and then show that some of those tangencies belong to class (2). 
\par{}
Let us create a secondary homoclinic tangency by making the curve $T_1\circ T^k_0\circ T_1(W^u_{loc}(O))$ intersect $W^s_{loc}(O)$ non-transversely. By formula \eqref{eq:setting:T_1} for $T_1$ (where one should take $y^+(\varepsilon)=\mu$), the image $(x_0,y_0,z_0)=T_1(x,y,z)$ of a point $(x,y,z)\in \Pi_1$ is given by
\begin{equation}\label{eq:preperturb:0}
\begin{array}{rcl}
x_0-x^+ &=& ax+b(y - y^-) + a_{13} z +h_1(x,y-y^-,z), \\ [5pt]
y_0 &=& \mu + cx + d(y-y^-)^2 + a_{23} z + h_2(x,y-y^-,z),\\[5pt]
z_0-z^+ &=& a_{31}x+a_{32}(y - y^-) + a_{33} z +h_3(x,y-y^-,z). \\ [5pt]
\end{array}
\end{equation}
Consequently, the image $T_1(W^u_{loc}(O))$ has the form
\begin{eqnarray}
y_0&=&\mu+\dfrac{d}{b^2}(x_0-x^+)^2+h_2(0,\frac{x_0-x^+}{b},0),\label{eq:preperturb:1.1}\\
z_0-z^+&=&\dfrac{a_{32}}{b}(x_0-x^+)+h_3(0,\frac{x_0-x^+}{b},0),\label{eq:preperturb:1.2}
\end{eqnarray}
where $h_2(0,(x_0-x^+)/{b},0)=o((x_0-x^+)^2)$ and $h_3(0,(x_0-x^+)/{b},0)=o(|x_0-x^+|)$. For any point $(x_0,y_0,z_0)\in T_1(W^u_{loc}(O))\cap \sigma^0_k$, we can find its $k$-th iterate $(x_k,y_k,z_k)=T^k_0(x_0,y_0,z_0)$ by formula \eqref{eq:setting:T^k}:
\begin{eqnarray}
x_k &=& \lambda^k x_0 + o(\lambda^k), \label{eq:preperturb:2.1}\\[5pt]
y_0 &=& \gamma^{-k} y_k + o(\gamma^{-k}), \label{eq:preperturb:2.2}\\[5pt]
z_k &=& O(\hat{\lambda}^k).\label{eq:preperturb:2.3}
\end{eqnarray}
The point $(x_0,y_0,z_0)$ is a homoclinic point if $T_1(x_k,y_k,z_k)=(\bar{x},\bar{y},\bar{z})\in W^s(O)$, namely, the coordinate $\bar{y}$ equals zero. From the second equation in \eqref{eq:setting:T_1}, we have
\begin{equation}\label{eq:preperturb:3}
\bar{y} = \mu + cx_k + d(y_k-y^-)^2 + a_{23} z_k + h_2(x_k,y_k-y^-,z_k)=0.
\end{equation}
By plugging \eqref{eq:preperturb:1.1} and \eqref{eq:preperturb:2.3} into \eqref{eq:preperturb:2.2}, and plugging \eqref{eq:preperturb:2.1} and \eqref{eq:preperturb:2.3} into \eqref{eq:preperturb:3}, we obtain the following system whose solutions correspond to homoclinic points $(x_0,y_0,z_0)\in T_1(W^u_{loc}(O))$:
\begin{equation}\label{eq:preperturb:4}
\begin{array}{rcl}
0&=&\mu-\gamma^{-k}y^- -\gamma^{-k}(y_k-y^-)+\dfrac{d}{b^2}(x_0-x^+)^2+u_1(x_0,y_k,\mu)+u_2(x_0,\mu),\\[5pt]
0&=&\mu+c\lambda^k x^+ +c\lambda^k(x_0-x^+)+d(y_k-y^-)^2+u_3(x_0,y_k,\mu)+u_4(x_0,y_k,\mu),
\end{array}
\end{equation}
where $u_1=o(\gamma^{-k}), u_2=o(x_0^2), u_3=o(\lambda^k)$, and $u_4=o(|\lambda|^k+y_k^2)$. After letting ${X}=x_0-x^+$ and ${Y}=y_k-y^-$, system \eqref{eq:preperturb:4} recasts as
\begin{equation}\label{eq:preperturb:5}
\begin{array}{rcl}
0&=&\mu-\gamma^{-k}y^- -\gamma^{-k}Y+\dfrac{d}{b^2}X^2+\hat{u}_1(X,Y,\mu)+\hat{u}_2(X,\mu),\\[5pt]
0&=&\mu+c\lambda^k x^+ +c\lambda^kX+dY^2+\hat{u}_3(X,Y,\mu)+\hat{u}_4(Y,\mu),
\end{array}
\end{equation}
where $\hat{u}_1=o({\gamma}^{-k}),\hat{u}_2=o({X}^2),\hat{u}_3=o({\lambda}^k)$ and $\hat{u}_4=o(|\lambda|^k+{Y}^2)$.

A non-degenerate homoclinic tangency corresponds to a solution to system \eqref{eq:preperturb:5} with multiplicity two. This corresponds to the vanishing determinant of the Jacobian matrix. Now, by letting the Jacoby matrix of system \eqref{eq:preperturb:5} have determinant zero, expressing $\mu$ as a function of $X$ and $Y$ from the first equation of \eqref{eq:preperturb:5}, and plugging this expression for $\mu$ into the second one, we arrive at the following system:
\begin{equation}\label{eq:preperturb:6.1}
\begin{array}{rcl}
0&=&c\lambda^k\gamma^{-k}+4\dfrac{d^2}{b^2}(X+v_1(X,Y))(Y+v_2(X,Y))+o(\lambda^k\gamma^{-k}),\\[5pt]
0&=&c\lambda^k x^+ +\gamma^{-k}y^- +c\lambda^k X+\gamma^{-k}Y+dY^2-\dfrac{d}{b^2}X^2+o(\lambda^k\gamma^{-k}),
\end{array}
\end{equation}
where $v_1=o(|\gamma|^{-k}+|{X}|)$ and $v_2=o(|\lambda|^k+|{Y}|)$. With the further coordinate transformation
\begin{equation}\label{eq:preperturb:6.2}
(\hat{X},\hat{Y})=(X+v_1(X,Y),Y+v_2(X,Y)),
\end{equation}
we obtain
\begin{equation}\label{eq:preperturb:6}
\begin{array}{rcl}
0&=&c\lambda^k\gamma^{-k}+4\dfrac{d^2}{b^2}\hat{X}\hat{Y}+o(\lambda^k\gamma^{-k}) ,\\[5pt]
0&=&c\lambda^k x^+ +\gamma^{-k}y^- +c\lambda^k \hat{X}+\gamma^{-k}\hat{Y}+d\hat{Y}^2-\dfrac{d}{b^2}\hat{X}^2+o(|\lambda|^k+|\gamma|^{-k}).
\end{array}
\end{equation}
Quadratic tangencies of the original system correspond to non-degenerate solutions to \eqref{eq:preperturb:6}, and the value of $\mu=\mu_k$ corresponding to such tangency can be found from either of the equations in \eqref{eq:preperturb:5}.

\begin{figure}[!h]
\begin{center}
\includegraphics[width=1\textwidth]{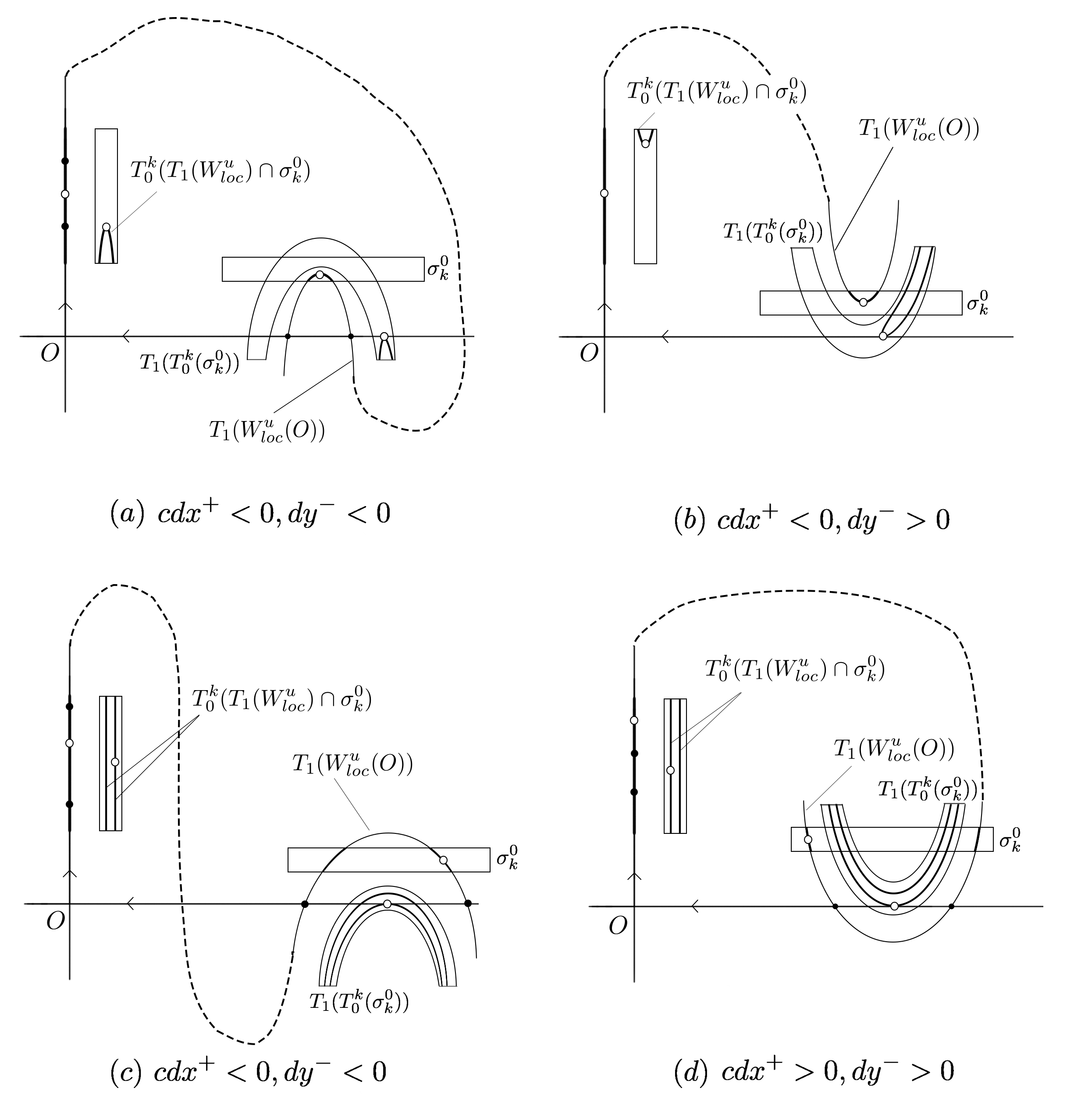}
\end{center}
\caption{Creation of secondary homoclinic tangencies for $x^+,y^->0$. Here we project the iterates of $W^u_{loc}(O)$ and $\sigma^0_k$ onto the two-dimensional plane $\{z=0\}$ along the leaves of $\mathcal{F}$ (note that such projection is well-defined by the non-degeneracy condition (C2)), and take $\mu=\mu^i_k$ for some $i\in\{1,2\}$. The horizontal and the vertical strips are the projections of $\sigma^0_k$ and $T^k_0(\sigma^0_k)$, and the hollowed dots denote the points in the orbit of the homoclinic tangency while the solid dots denote those in the transverse homoclinic orbits.}
\label{fig:Lemma2}
\end{figure}

In what follows, we find solutions to \eqref{eq:preperturb:6}. Let $k$ be even so that $\lambda^k$ and $\gamma^{-k}$ are always positive. Consider first class (1), where $cdx^+<0$ and $dy^-<0$. We do the following scaling:
\begin{equation*}\label{eq:preperturb:11}
(\hat{X},\hat{Y})\mapsto |\lambda|^\frac{k}{2}\sqrt{\left|\dfrac{cx^+}{d}\right|}\bigg(-\dfrac{b^2\gamma^{-k}}{4dx^+}U,V\bigg).
\end{equation*}
In the new variables system \eqref{eq:preperturb:6} takes the form 
\begin{equation}\label{eq:preperturb:12}
\begin{array}{rcl}
1&=&UV+o(1)_{k\to +\infty},\\
1&=&V^2+o(1)_{k\to +\infty}.
\end{array}
\end{equation}
For any sufficiently large $k$ the above system has two non-degenerate solutions $(1+o(1),1+o(1))$ and $(-1+o(1),-1+o(1))$, corresponding to two solutions to system \eqref{eq:preperturb:6}:
\begin{equation}\label{eq:preperturb:12.1}
\begin{array}{rcl}
(\hat{X}_k^1,\hat{Y}_k^1)&=&\left(-\dfrac{b^2|\lambda|^\frac{k}{2}\gamma^{-k}}{4dx^+}\sqrt{\left|\dfrac{cx^+}{d}\right|}+o(|\lambda|^\frac{k}{2}\gamma^{-k}),|\lambda|^\frac{k}{2}\sqrt{\left|\dfrac{cx^+}{d}\right|}+o(|\lambda|^\frac{k}{2})\right),\\[15pt]
(\hat{X}_k^2,\hat{Y}_k^2)&=&\left(\dfrac{b^2|\lambda|^\frac{k}{2}\gamma^{-k}}{4dx^+}\sqrt{\left|\dfrac{cx^+}{d}\right|}+o(|\lambda|^\frac{k}{2}\gamma^{-k}),-|\lambda|^\frac{k}{2}\sqrt{\left|\dfrac{cx^+}{d}\right|}+o(|\lambda|^\frac{k}{2})\right).
\end{array}
\end{equation}
These two solutions give us two homoclinic tangency points ${M}_k^1,{M}_k^2\in T_1(W^u_{loc}(O)$ for two different $\mu$ values $\mu_k^1$ and $\mu_k^2$ (see Figure \ref{fig:Lemma2}(a)). From equations \eqref{eq:preperturb:1.2}, \eqref{eq:preperturb:2.2} and \eqref{eq:preperturb:6.2}, we find the coordinates of these tangency points as
\begin{equation}\label{eq:preperturb:9}
{M}_k^1=(\hat{X}^1+x^++o(\gamma^{-k}),\gamma^{-k}(\hat{Y}^1+y^-+o(1)),z_1)\quad \mbox{and}\quad {M}_k^2=(\hat{X}^2+x^+,\gamma^{-k}(\hat{Y}^2+y^-+o(1)),z_2),
\end{equation}
where we do not write the $z$-coordinates explicitly. Let $M^-_k=(0,y^-_k,0)\in W^u_{loc}(O)$ be the pre-image of any of the points ${M}_k^1$ and ${M}_k^2$. By \eqref{eq:preperturb:0} and \eqref{eq:preperturb:9}, we have $y^-_k - y^-=(\hat{X}^i_k+o(\hat{X}^i_k)+o(\gamma^{-k}))/b$. This immediately shows that $M^-_k \to M^-$ as $k\to +\infty$. The first equation in \eqref{eq:preperturb:5} yields the corresponding $\mu$ values, which are $\mu_k^{i}=\gamma^{-k}y^-(1+o(1))$ ($i=1,2$). 
\begin{rem}\label{class2}
Note that the condition $dy^-<0$ has not been used in the above computation. In fact, we can also create new tangencies for class (2) in the same way (see Figure \ref{fig:Lemma2}(b)).
\end{rem}
Now consider classes (3) and (4), where we have $cdx^+>0$. By using the scaling
\begin{equation*}\label{eq:preperturb:7}
(\hat{X},\hat{Y})\mapsto b|\lambda|^\frac{k}{2}\sqrt{\dfrac{cx^+}{d}}\bigg(U,-\dfrac{\gamma^{-k}}{4dx^+}V\bigg),
\end{equation*}
and dividing the first and second equation of \eqref{eq:preperturb:6} to $c\lambda^k\gamma^{-k}$ and $c\lambda^k x^+$, respectively, we arrive at the following system
\begin{equation}\label{eq:preperturb:8}
\begin{array}{rcl}
1&=&UV+o(1)_{k\to +\infty},\\
1&=&U^2+o(1)_{k\to +\infty}.
\end{array}
\end{equation}

For any sufficiently large $k$, system \eqref{eq:preperturb:8} has non-degenerate solutions $(1+o(1),1+o(1))$ and $(-1+o(1),-1+o(1))$, which lead to two solutions to system \eqref{eq:preperturb:6} as
\begin{equation}\label{eq:preperturb:8.1}
\begin{array}{rcl}
(\hat{X}_k^1,\hat{Y}_k^1)&=&\left(b|\lambda|^\frac{k}{2}\sqrt{\dfrac{cx^+}{d}}+o(|\lambda|^\frac{k}{2}),-
\dfrac{b|\lambda|^\frac{k}{2}\gamma^{-k}}{4dx^+}\sqrt{\dfrac{cx^+}{d}}+o(|\lambda|^\frac{k}{2}\gamma^{-k})\right),\\[15pt]
(\hat{X}_k^2,\hat{Y}_k^2)&=&\left(-b|\lambda|^\frac{k}{2}\sqrt{\dfrac{cx^+}{d}}+o(|\lambda|^\frac{k}{2}),
\dfrac{b|\lambda|^\frac{k}{2}\gamma^{-k}}{4dx^+}\sqrt{\dfrac{cx^+}{d}}+o(|\lambda|^\frac{k}{2}\gamma^{-k})\right).
\end{array}
\end{equation}
For each sufficiently large $k$, these two solutions give us two points of homoclinic tangency ${M}_k^1,{M}_k^2\in T_1(W^u_{loc}(O))$ (see Figure \ref{fig:Lemma2}(c) and (d)). Similar to the discussion for class (1), for the pre-image $M^-_k$ of any of the points ${M}_k^1$ and ${M}_k^2$, we have $M^-_k \to M^-$ as $k\to +\infty$. The corresponding $\mu$ values can be found from the second equation in \eqref{eq:preperturb:5}, which gives $\mu_k^i=-cx^+\lambda^k(1+o(1))$ ($i=1,2$). 

We proceed to compute the signs of the coefficients $c$ and $d$ corresponding to the new homoclinic tangencies. We have shown that for each sufficiently large $k$ there exist two values of $\mu=\mu^i_k(i=1,2)$ that correspond to a homoclinic tangency. The associated global map for this tangency is
$$\hat{T}:=T_1\circ T^k_0\circ T_1: (x,y,z)\mapsto(\bar{x},\bar{y},\bar{z}).$$
By denoting $T^{-1}_1(M^i_k)=(0,y^i_k,0)$, the coefficients $c^i_k$ and $d^i_k$ of $\hat{T}$ are given by

\begin{equation}\label{eq:preperturb9}
c^i_k=\dfrac{\partial \bar{y}(0,y^i_k,0)}{\partial x}\quad\mbox{and}\quad
d^i_k=\frac{1}{2}\dfrac{\partial^2 \bar{y}(0,y^i_k,0)}{\partial y^2},
\end{equation}
where $\bar{y}$ is related to $(x_k,y_k,z_k)=T^k_0(x_0,y_0,z_0)=T^k_0\circ T_1(x,y,z)$ by \eqref{eq:preperturb:3}. We note from \eqref{eq:preperturb:4} - \eqref{eq:preperturb:6.1} that 
\begin{equation}\label{eq:preperturb:10a}
\hat{Y}=Y+v_2=Y+\dfrac{1}{2d}\left(\dfrac{\partial \hat{u}_3}{\partial Y}+\dfrac{\partial \hat{u}_4}{\partial Y}\right)
=\dfrac{1}{2d}\left(2d(y_k-y^-)+\dfrac{\partial (cx_k+a_{23}z_k)}{\partial y_k}+\dfrac{\partial h_2}{\partial y_k}\right)=\dfrac{1}{2d}\dfrac{\partial \bar{y}}{\partial y_k}.
\end{equation}
This fact along with equations \eqref{eq:preperturb:0} and \eqref{eq:preperturb:2.1} - \eqref{eq:preperturb:2.3} yields
\begin{equation}\label{eq:preperturb10}
\begin{array}{rcl}
c^i_k
&=&\left(\dfrac{\partial \bar{y}}{\partial x_k}\dfrac{\partial x_k}{\partial x}+\dfrac{\partial \bar{y}}{\partial y_k}\dfrac{\partial y_k}{\partial x}+\dfrac{\partial \bar{y}}{\partial z_k}\dfrac{\partial z_k}{\partial x}\right)\left|_{(x,y,z)=(0,y^i_k,0)}\right.\\[15pt]
&=&ac\lambda^k+2cd\gamma^{k}\hat{Y}^i_k+o(\hat{Y}^i_k)+o({\lambda}^k),
\end{array}
\end{equation}
where $\hat{Y}^i_k$ is given by \eqref{eq:preperturb:12.1} or \eqref{eq:preperturb:8.1}.

Let us now compute $d^i_k$ which is given by
\begin{equation}\label{eq:preperturb:d}
d^i_k=\frac{1}{2}\dfrac{\partial}{\partial y}\left(
\dfrac{\partial \bar{y}}{\partial x_k}\dfrac{\partial x_k}{\partial y}+\dfrac{\partial \bar{y}}{\partial y_k}\dfrac{\partial y_k}{\partial y}+\dfrac{\partial \bar{y}}{\partial z_k}\dfrac{\partial z_k}{\partial y}
\right)\left|_{(x,y,z)=(0,y^i_k,0)}\right..
\end{equation}
It can be easily seen from \eqref{eq:preperturb:0} and \eqref{eq:preperturb:2.1} - \eqref{eq:preperturb:2.3} that 
\begin{equation}
\dfrac{\partial}{\partial y}\left(
\dfrac{\partial \bar{y}}{\partial x_k}\dfrac{\partial x_k}{\partial y}+\dfrac{\partial \bar{y}}{\partial z_k}\dfrac{\partial z_k}{\partial y}
\right)\left|_{(x,y,z)=(0,y^i_k,0)}\right.=o(\lambda^k).
\end{equation}
Regarding the rest of the derivatives in \eqref{eq:preperturb:d}, we note from the first equation of \eqref{eq:preperturb:0} that
\begin{equation*}
y-y^-=\dfrac{(x_0-x^+)}{b}+o(x_0-x^+)=\dfrac{X}{b}+o(X)=\dfrac{\hat{X}+o(\gamma^{-k})}{b}(1+o(1)),
\end{equation*}
see \eqref{eq:preperturb:6.1}. Together with equations \eqref{eq:preperturb:0} and \eqref{eq:preperturb:2.2}, this leads to
\begin{equation}\label{eq:preperturb:10b}
\dfrac{\partial y_k}{\partial y}\left|_{(x,y,z)=(0,y^i_k,0)}\right.
=2d\gamma^k(y^i_k-y^-)+o(y^i_k-y^-)=\dfrac{2d\gamma^k}{b}(\hat{X}^i_k+o(\gamma^{-k}))(1+o(1)),
\end{equation}
where $\hat{X}^i_k$ is given by \eqref{eq:preperturb:12.1} or \eqref{eq:preperturb:8.1}. Now, with the help of \eqref{eq:preperturb:10a} and \eqref{eq:preperturb:10b}, we obtain
\begin{equation}\label{eq:preperturb11}
\begin{array}{rcl}
d^i_k
&=&\frac{1}{2}\left(
\dfrac{\partial^2 \bar{y}}{\partial y_k^2}\left(\dfrac{\partial y_k}{\partial y}\right)^2+ \dfrac{\partial \bar{y}}{\partial y_k}\dfrac{\partial^2 y_k}{\partial y^2}  
\right)\left|_{(x,y,z)=(0,y^i_k,0)}\right.+ o(\lambda^k)\\[15pt]
&=&\dfrac{4d^3\gamma^{2k}}{b^2}(\hat{X}^i_k+o(\gamma^{-k}))^2(1+o(1))+2d^2\gamma^k\hat{Y^i_k}+o(\lambda^k).
\end{array}
\end{equation}

For class (1), where $cdx^+<0$ and $dy^-<0$, we plug the solutions \eqref{eq:preperturb:12.1} into the above equations and get 
\begin{equation}\label{eq:preperturb12}
\begin{array}{rcl}
c^i_k&=&(-1)^{(i+1)}2cd|\lambda|^\frac{k}{2}\gamma^k  \sqrt{\left|\dfrac{cx^+}{d}\right|}+o(|\lambda|^\frac{k}{2}\gamma^k),\\[15pt]
d^i_k&=&(-1)^{(i+1)}2d^2|\lambda|^\frac{k}{2}\gamma^k \sqrt{\left|\dfrac{cx^+}{d}\right|}+o(|\lambda|^\frac{k}{2}\gamma^k),
\end{array}
\end{equation}
which implies $c^1_kd^1_k x^+<0$ and $d^1_ky^->0$. Therefore, by taking $\mu_k=\mu^1_k$ and $M^-_k=T^{-1}_1(M^1_k)$, we obtain a homoclinic tangency that belongs to class (2), as required.

Let now $cdx^+>0$. With the corresponding solutions \eqref{eq:preperturb:8.1}, equations \eqref{eq:preperturb10} and \eqref{eq:preperturb11} yield
\begin{equation}\label{eq:preperturb13}
\begin{array}{rcl}
c^i_k&=&(-1)^i \dfrac{bc|\lambda|^\frac{k}{2}}{2x^+}\sqrt{\dfrac{cx^+}{d}}\left(1+o(1)\right),\\[15pt]
d^i_k&=&4cd^2\lambda^k\gamma^{2k}x^+ +o(\lambda^k\gamma^{2k}).
\end{array}
\end{equation}
Observe that $c^i_k(i=1,2)$ have different signs and $d^i_k$ always have the same sign as $d$. It follows that for class (4) where $cdx^+>0$ and $dy^->0$ one can obtain the desired class (2) homoclinic tangency by picking $i$ such that $c^i_k<0$. If the original tangency belongs to class (3) where $cdx^+>0$ and $dy^-<0$, then we can first obtain a class (1) tangency by choosing $i$ such that $c^i_k>0$. After this, repeat what we did for class (1) tangency. \qed

\begin{figure}[!h]
\begin{center}
\includegraphics[width=0.6\textwidth]{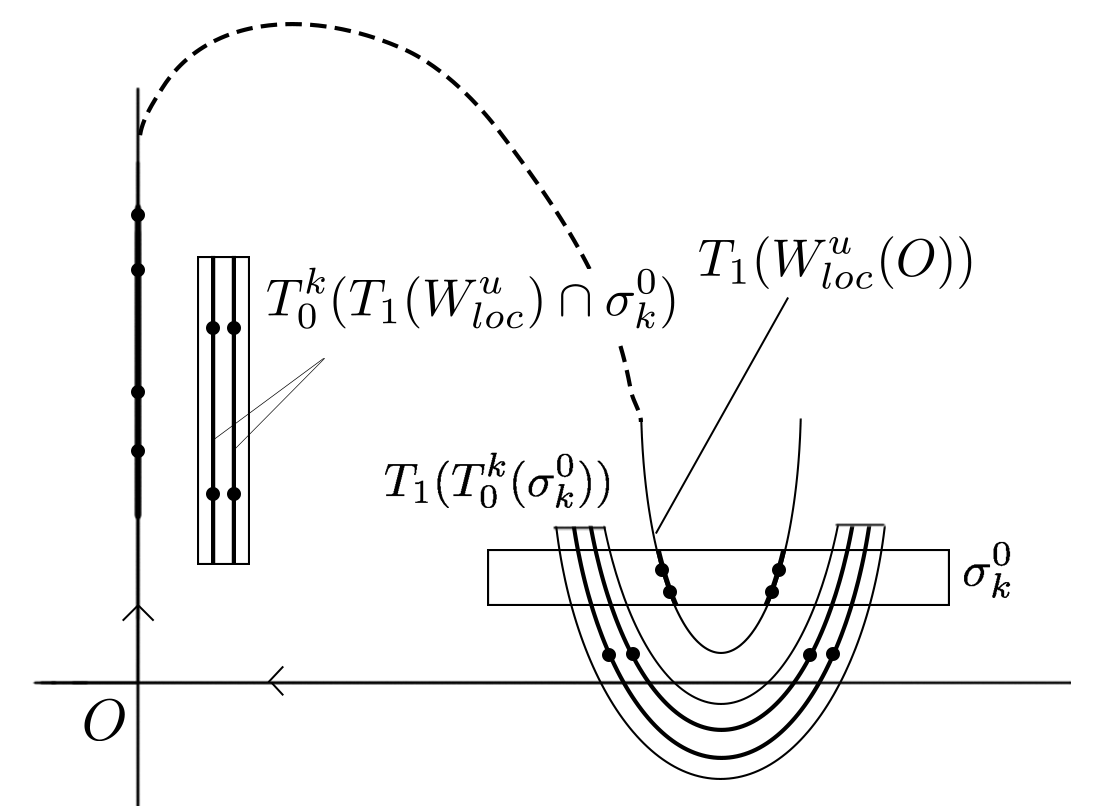}
\end{center}
\caption{Transverse homoclinic points at $\mu=0$.}
\label{fig:Lemma2a}
\end{figure}

We are now in the position to show that a homoclinic tangency satisfying conditions (a) and (b) can be recovered from any kind of the original tangency.
\begin{lem}\label{lem:preperturb2}
For any smooth one-parameter family $F_\mu$ of diffeomorphisms with the diffeomorphism $F_0$ satisfying conditions (C1)-(C3), there exists a sequence $\{\mu_k\}$ accumulating on $\mu=0$ such that the saddle $O$ of $F_{\mu_k}$ has a new homoclinic tangency point $M^-_k$ for which
$cdx^+>0$ and $cx^+y^->0$, and in $W^u_{loc}(O)\cap\Pi_1$ there exist two transverse homoclinic points $N^1_k$ and $N_k^2$ such that the $y$-coordinate of $M^-_k$ lies between those of $N^1_k$ and $N^2_k$. The distance between the points $N^{1,2}_k$ and $M_k^-$ 
tends to zero as $k\to +\infty$.
\end{lem}
\noindent{\it Proof.} By Lemma \ref{lem:preperturb1}, it is sufficient to prove Lemma \ref{lem:preperturb2} only for the case where the homoclinic tangency of $F_0$ belongs to class (2), namely, we may assume that $cdx^+<0$ and $dy^->0$. We start with showing that in this case there exist infinitely many transverse homoclinic points at $\mu=0$. Indeed, non-degenerate solutions of system \eqref{eq:preperturb:5} correspond to transverse homoclinic points. By using the scaling
$$
(X,Y)\mapsto \left( b|\gamma|^{-\frac{k}{2}}\sqrt{\dfrac{y^-}{d}}U,|\lambda|^{\frac{k}{2}}\sqrt{\dfrac{cx^+}{d}}V\right),
$$
we rewrite system \eqref{eq:preperturb:5} at $\mu=0$ as
\begin{equation}
\begin{array}{rcl}
1&=&U^2+o(1)_{k\to +\infty},\\[5pt]
1&=&V^2+o(1)_{k\to +\infty}.
\end{array}
\end{equation}
This gives four non-degenerate solutions to \eqref{eq:preperturb:5} at $\mu=0$
\begin{equation}\label{eq:preperturb2:1}
(X,Y)=\left(\pm b|\gamma|^{-\frac{k}{2}}\sqrt{\dfrac{y^-}{d}}+o(|\gamma|^{-\frac{k}{2}}),\pm|\lambda|^{\frac{k}{2}}\sqrt{\dfrac{cx^+}{d}}+o(|\lambda|^{\frac{k}{2}})\right)
=:\left(\pm \tilde{X}+o(|\gamma|^{-\frac{k}{2}}),\pm\tilde{Y}+o(|\lambda|^{\frac{k}{2}})\right)
\end{equation}
for any sufficiently large $k$. These solutions correspond to four transverse homoclinic points in $T_1(W^u_{loc}(O))$:
\begin{equation}\label{eqtra}
\begin{array}{rcl}
N^1_k=\left(x^++\tilde{X}+o(|\gamma|^{-\frac{k}{2}}),\gamma^{-k}(y^-+\tilde{Y})+o(\gamma^{-k}),z^1\right),\\[5pt]
N^2_k=\left(x^++\tilde{X}+o(|\gamma|^{-\frac{k}{2}}),\gamma^{-k}(y^--\tilde{Y})+o(\gamma^{-k}),z^2\right),\\[5pt]
N^3_k=\left(x^+-\tilde{X}+o(|\gamma|^{-\frac{k}{2}}),\gamma^{-k}(y^-+\tilde{Y})+o(\gamma^{-k}),z^3\right),\\[5pt]
N^4_k=\left(x^+-\tilde{X}+o(|\gamma|^{-\frac{k}{2}}),\gamma^{-k}(y^--\tilde{Y})+o(\gamma^{-k}),z^4\right).
\end{array}
\end{equation}
Denote $T^{-1}_1(N^i_k)$ by $\hat{N}^i_k=(0,\hat{y}^i_k,0)$. It follows from the first equation of \eqref{eq:setting:T_1} that $\hat{y}^{1,2}_k>y^-$ and $\hat{y}^{3,4}_k<y^-$, which means that the tangency point $M^-$ is bounded by the four transverse homoclinic points $\hat{N}^i_k$ (see Figure \ref{fig:Lemma2a}). Moreover, we have from the second equation of \eqref{eq:setting:T_1} that $\hat{y}^1_k>\hat{y}^2_k$ and $\hat{y}^3_k>\hat{y}^4_k$. By transversality, for each fixed k, all four homoclinic intersections persist for all sufficiently small $\mu$.

In what follows we prove that there exists a sequence $\{\mu_m\}$ accumulating on $\mu=0$ such that for each sufficiently large $m$ the diffeomorphism $F_{\mu_m}$ has a non-transverse homoclinic point $M^-_m\in W^u_{loc}(O)$ that belongs to class (4) and satisfies either $M^-_m\to\hat{N}^2_k$ or $M^-_m\to\hat{N}^3_k$ as $m\to +\infty$. This will complete the proof of the lemma after noting that class (4) tangencies satisfy condition (a), both $\hat{N}^2_k$ and $\hat{N}^3_k$ are bounded by the two transverse homoclinic points $\hat{N}^1_k$ and $\hat{N}^4_k$, and these points all tend to $M^-$ as $k\to +\infty$.  

\begin{figure}[h]
\begin{center}
\includegraphics[width=0.5\textwidth]{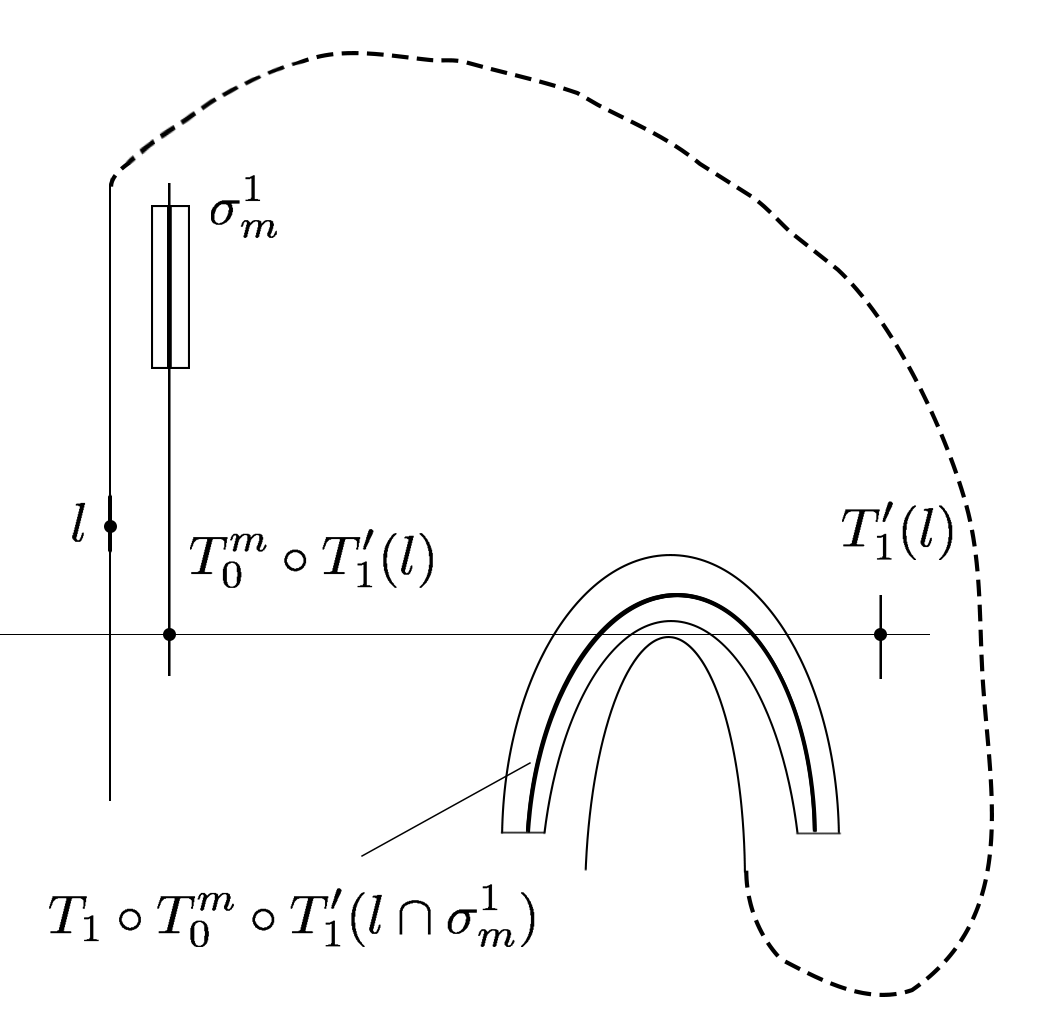}
\end{center}
\caption{By changing $\mu$, one can make $T_1\circ T^m_0 \circ T'_1(l)$ intersect $W^s_{loc}(O)$ non-transversely. Here $l\in W^u_{loc}(O)$ is a small piece containing the transverse homoclinic point.}
\label{fig:model}
\end{figure}

We denote as $T_1'$ the restriction of the global map $T_1$ to a small neighbourhood of the transverse homoclinic point $\hat{N}^2_k=(0,\hat y_k,0)$.
We denote $T_1'(\hat{N}_k^2)=N^2_k=(\hat{x}^+,0,\hat{z}^+)$ and write the Taylor expansion of $T_1'$ about the point $\hat{N}^2_k$
as 
\begin{equation}\label{eq:model2}
\begin{array}{rcl}
\bar{x}-\hat{x}^+ &=& a'x+b'(y - \hat{y}_k) + a'_{13} z +h'_1(x,y,z), \\ [5pt]
\bar{y} &=& c'x + d'(y-\hat{y}_k) + a'_{23} z + h'_2(x,y,z),\\[5pt]
\bar{z}-\hat{z}^+ &=& a'_{31}x+a'_{32}(y - \hat{y}_k) + a'_{33} z +h'_3(x,y,z), \\ [5pt]
\end{array}
\end{equation}
where $h'_{1,2,3}=O(x^2+y^2+z^2)$. The coefficients in these formula are obtained by evaluating, at $(0,\hat{y}_k,0)$, the first derivatives of the map $T_1$ given by \eqref{eq:setting:T_1}. Obviously,
\begin{equation}\label{eq:model2aa}
a'=a+\dots,\quad b'=b+\dots,\quad c'=c+\dots,\quad d'=2d(\hat{y}_k-y^-)(1+\dots),
\end{equation}
where the dots denote terms that tend to zero as $k\to +\infty$. We now create a homoclinic tangency by finding a point $M^-_m\in W^u_{loc}(O)$ close to $\hat{N}^2_k$ such that $M_m^+:=T_1\circ T^m_0 \circ T'_1(M)\in W^s_{loc}(O)$ for some $m$, and the curve 
$T_1\circ T^m_0 \circ T'_1 (W^u_{loc}(O))$ is tangent to $W^s_{loc}(O)$ at the point $M_m^+$ as shown in Figure \ref{fig:model}. 

Let $m$ be even, so that $\lambda^m$ and $\gamma^{-m}$ are positive. The image $T'_1(W^u_{loc}(O))$ is given by
\begin{eqnarray*}
y_0&=&\dfrac{d'}{b'}(x_0-\hat{x}^+)+o(x_0-\hat{x}^+),\\
z_0-\hat{z}^+&=&\dfrac{a'_{32}}{b'}(x_0-\hat{x}^+)+o(x_0-\hat{x}^+).
\end{eqnarray*}
For any point $(x_0,y_0,z_0)\in T'_1(W^u_{loc}(O))\cap \sigma^0_m$, we can find its $m$-th iterate $(x_m,y_m,z_m)=T^m_0(x_0,y_0,z_0)$ by using formula \eqref{eq:setting:T^k}:
\begin{eqnarray}\label{eq:model2a}
\begin{array}{rcl}
x_m &=& \lambda^m x_0 + o(\lambda^m), \\[5pt]
y_0 &=& \gamma^{-m} y_m + o(\gamma^{-m}), \\[5pt]
z_m &=& O(\hat{\lambda}^m),
\end{array}
\end{eqnarray}
The point $(x,y,z)$ is a homoclinic point if and only if $T_1(x_m,y_m,z_m)\in W^s(O)$, namely,
$$0 = \mu + cx_m + d(y_m-y^-)^2 + a_{23} z_m + h_2(x_m,y_m,z_m).$$
Then, by repeating the same procedure as was used to find equation \eqref{eq:preperturb:5}, we obtain  
\begin{equation}\label{eq:model3}
\begin{array}{rcl}
0&=&-\gamma^{-m}y^- -\gamma^{-m}Y+\dfrac{d'}{b'}X+{u}_1(X,Y,\mu)+{u}_2(X,\mu),\\[5pt]
0&=&\mu+c\lambda^m \hat{x}^+ +c\lambda^mX+dY^2+{u}_3(X,Y,\mu)+{u}_4(Y,\mu),
\end{array}
\end{equation}
where $X=x-\hat{x}^+$, $Y=y_m-y^-$, ${u}_1=o({\gamma}^{-m}),{u}_2=o(X),{u}_3=o({\lambda}^m)$, and ${u}_4=o(\lambda^m+Y^2)$.
\par{}
In order to have a homoclinic tangency, we need the Jacobian matrix of the right-hand side of \eqref{eq:model3} to have zero determinant, namely,
\begin{equation}\label{eq:model4a}
c\lambda^m\gamma^{-m}+\dfrac{2dd'}{b'}(Y+v(X,Y))+o(Y)+o(\lambda^m\gamma^{-m})=0,
\end{equation}
where $v=o(\lambda^m+|Y|)$. After the coordinate transformation
\begin{equation}\label{eq:model4b}
(\hat{X},\hat{Y})=(X,Y+v(X,Y)), 
\end{equation}
equations \eqref{eq:model3} keep their form, and equation \eqref{eq:model4a} is recast as
\begin{equation}\label{eq:model4}
c\lambda^m\gamma^{-m}+\dfrac{2dd'}{b'}\hat{Y}+o(\hat{Y})+o(\lambda^m\gamma^{-m})=0.
\end{equation}
The quadratic tangencies correspond to non-degenerate solutions to the system consisting of \eqref{eq:model3} and \eqref{eq:model4}. With a straightforward computation one can find the solutions as
\begin{equation}\label{eq:model5}
\begin{array}{rcl}
\hat{X}_m&=&\dfrac{b'\gamma^{-m}y^-}{d'}+o(\gamma^{-m}),\\[5pt]
\hat{Y}_m&=&-\dfrac{b'c\lambda^m\gamma^{-m}}{2dd'}+o(\lambda^m\gamma^{-m}),\\[5pt]
\mu_m&=&-c\lambda^m \hat{x}^+ +o(\lambda^m),
\end{array}
\end{equation}
where $m$ is sufficiently large, and each solution gives a non-transverse homoclinic point 
$M^-_m \in W^u_{loc}(O)$ corresponding to a quadratic tangency at $\mu={\mu}_m$.

The global map associated to $M^-_m$ is $\hat{T}:=T_1\circ T^m_0\circ T'_1:(x,y,z)\mapsto(\bar{x},\bar{y},\bar{z})$, and the corresponding coefficients $c_m$ and $d_m$ are given by
\begin{equation}\label{eq:model6}
\hat{c}=\dfrac{\partial \bar{y}}{\partial x}_{M^-_m} \quad\mbox{and}\quad \hat{d}=\frac{1}{2}\dfrac{\partial^2 \bar{y}}{\partial y^2}|_{M^-_m}.
\end{equation}
Similar to the computation of such coefficients in the proof of Lemma \ref{lem:preperturb1}, by applying the chain rule to equations \eqref{eq:setting:T_1}, \eqref{eq:model2} and \eqref{eq:model2a}, and using the formulas \eqref{eq:model4a} and \eqref{eq:model5}, we have
\begin{eqnarray}
c_m&=&a'c\lambda^m + 2c'd\gamma^m\hat{Y}_m + o(\lambda^m)= c\lambda^m\left(\dfrac{a'd'-b'c'}{d'}\right)+o(\lambda^m),\label{eq:model7}\\
d_m&=&dd^{'2}\gamma^{2m}+o(\gamma^m).\label{eq:model8}
\end{eqnarray}
Equation \eqref{eq:model8} means that $d_m$ has the same sign as $d$, which is positive.
Equation \eqref{eq:model7} for $c_m$ can be recast as 
\begin{equation}
c_m=c\lambda^m\left(\dfrac{2ad(\hat{y}_k-y^-)-bc+\dots}{2d(\hat{y}_k-y^-)(1+\dots)}\right)+o(\lambda^m).
\end{equation}
Since $\hat{y}_k-y^-$ can be sufficiently small, the estimates in \eqref{eq:model2aa} imply that the sign of $c_m$ is the same as $-b/(d(\hat{y}_k-y^-))$. It follows from $d>0$ and $\hat{y}_k-y^->0$ that if $b<0$, then we have $c_m>0$, and this gives us the class (4) homoclinic tangency; if $b>0$, then we just need to consider the point $\hat{N}^3_k$, for which $\hat y_k-y^-<0$, instead of $\hat{N}^2_k$ in \eqref{eq:model2}.\qed

\section{Invariant cone fields}\label{sec:cones}
In this section, we prove the existence of certain invariant cone fields in $\Pi_0$. These cone fields will help in two ways. First, estimates for the strong-stable leaves are obtained from stable invariant cones in Lemmas \ref{lem:ssleavesweak} and \ref{lem:ssleaves}. Second, we use the cones to obtain estimates for the multipliers of periodic orbits.

Recall that $\sigma^0_k\subset\Pi_0$ $(k\geqslant k^*)$ are the sets of points whose images under $T^k_0$ belong to $\Pi_1$, where $k^*$ is the smallest integer such that $T^{k^*}_0(\Pi_0)\cap \Pi_1 \neq \emptyset$. Denote by $\Sigma^0$ the union of all $\sigma^0_k$ with $k\geqslant k^*$. For any $X\in\Sigma^0$, we have $T(X)=T_1\circ T^k_0(X)$ where $k$ is such that $X\in\sigma^0_k$.

\begin{lem}\label{lem:ucone}
If $k^*$ is sufficiently large, then there exist constants $K>0$ and $M>0$ such that the
cone field $\mathcal{C}^{cu}$ over $\Sigma^0$ (the center unstable cone filed)
defined as
\begin{equation}\label{eq:ucone:1}
\mathcal{C}^{cu}(X)=\{(\Delta x, \Delta y, \Delta z)\mid \|\Delta z\|\leqslant K(|\Delta x|+|\Delta y|)\}
\end{equation}
is strictly forward-invariant under the derivative $\D T$ of the first-return map $T$
(here, $(\Delta x, \Delta y, \Delta z)$ are coordinates in the tangent space to $\Sigma^0$). 
Moreover,
\begin{equation}\label{dtvxy}
\|\D T(X) V\| \geq M |\lambda|^k \|V\| 
\end{equation}
for any $V\in \mathcal{C}^{cu}(X)$.
\end{lem}

\noindent{\it Proof.}
Take any $X\in\sigma^0_k$ and let $V_0=(\Delta x_0, \Delta y_0, \Delta z_0)$ be a vector in the tangent space at the point $X$ such that 
\begin{equation}\label{eq:ucone:0}
\|\Delta z_0\|\leqslant K (|\Delta x_0|+|\Delta y_0|),
\end{equation}
where $K>0$ is some constant. Denote $\D T^{k}_0(X)V_0=(\Delta x_1, \Delta y_1, \Delta z_1)$ and $\D T_1\D T^{k}_0(X)V_0=(\Delta x_2,$ $\Delta y_2, \Delta z_2)$. By formula \eqref{eq:setting:T^k} and noting that the first derivatives of the functions $\phi,\hat{\phi}$ and $\psi$ in \eqref{eq:setting:T^k} are bounded, we have the following relations:
\begin{eqnarray}
&\Delta x_1 = \lambda^{k} \Delta x_0 + o(\lambda^{k})(\Delta x_0+\Delta y_1 + \Delta z_0), \label{eq:ucone:3}\\[10pt]
&\Delta y_0 = \gamma^{-{k}} \Delta y_1 + o({\gamma^{-{k}}})(\Delta x_0+\Delta y_1 + \Delta z_0), \label{eq:ucone:4}\\[10pt]
&\Delta z_1 = O(\hat{\lambda}^{k})(\Delta x_0+\Delta y_1 + \Delta z_0).\label{eq:ucone:5}
\end{eqnarray} 
Equations \eqref{eq:ucone:3} and \eqref{eq:ucone:4} can be recast as
\begin{eqnarray}
&\Delta x_0 = \lambda^{-{k}}\Delta x_1 (1+o(1))+ o(1) (\Delta y_1+\Delta z_0),\label{eq:ucone:6a}\\
&\Delta y_0=\gamma^{-k}\Delta y_1(1+o(1)) + o(\gamma^{-{k}})(\Delta x_0+\Delta z_0).\label{eq:ucone:6b}
\end{eqnarray}
By plugging these two equations into \eqref{eq:ucone:0}, we obtain
$$\|\Delta z_0\|\leqslant K|\lambda|^{-{k}}|\Delta x_1| (1+o(1))+ o(1) |\Delta y_1|,$$
where we denote by $o(1)$ the terms that go to zero as $k\to +\infty$. The above equation together with \eqref{eq:ucone:5} and \eqref{eq:ucone:6a} implies
\begin{equation}\label{eq:ucone:8}
\|\Delta z_1 \|\leqslant O(\hat{\lambda}^{k}{\lambda^{-{k}}})|\Delta x_1| + O(\hat{\lambda}^{k})|\Delta y_1|
\end{equation}
and
\begin{equation}\label{eqsexp}
\|\Delta x_0\| + \|\Delta y_0\| = O(\lambda^{-{k}}) (\|\Delta x_1\| + \|\Delta y_1\|).
\end{equation}

The derivative $\D T_1$ is uniformly bounded in a small neighbourhood $\Pi_1$, so we have 
$$\|\Delta z_2\| \leqslant \sup \|\D T_1\|(|\Delta x_1| +|\Delta y_1|+\|\Delta z_1\|).$$
Hence, when $k^*$ is large enough, the above inequality together with \eqref{eq:ucone:8} gives
\begin{equation}\label{eq:ucone:10}
\|\Delta z_2\| \leqslant (1+\sup \|\D T_1\|)(|\Delta x_1| +|\Delta y_1|).
\end{equation}
Note that by \eqref{eq:setting:T_1} we have 
$$\begin{pmatrix}
\Delta x_2\\
\Delta y_2
\end{pmatrix}
=B_1
\begin{pmatrix}
\Delta x_1\\
\Delta y_1
\end{pmatrix}
+
B_2\Delta z_1,
$$
for some matrices $B_1$ and $B_2$, whose norm is uniformly bounded. In fact, $B_1$ is close to $\left(\begin{array}{cc} a & b \\ c & 0\end{array}\right)$, so, by \eqref{eq:nc}, 
$\det(B_1)\neq 0$, i.e., $B_1$ is invertible. Thus, we have
\begin{equation}\label{eq:ucone:12}
\begin{pmatrix}
\Delta x_1\\
\Delta y_1
\end{pmatrix}
=B_1^{-1}
\begin{pmatrix}
\Delta x_2\\
\Delta y_2
\end{pmatrix}
-
B_1^{-1}B_2\Delta z_1.
\end{equation}
By taking $k^*$ sufficiently large, equations \eqref{eq:ucone:8} and \eqref{eq:ucone:12} imply
\begin{equation}\label{eq:ucone:13}
\left\|\begin{pmatrix}
\Delta x_1\\
\Delta y_1
\end{pmatrix}\right\|
\leqslant 2\|B_1^{-1}\|
\left\|\begin{pmatrix}
\Delta x_2\\
\Delta y_2
\end{pmatrix}\right\|.
\end{equation}
We now combine the two inequalities \eqref{eq:ucone:10} and \eqref{eq:ucone:13}. It follows that, by taking $k^*$ sufficiently large, we have 
\begin{equation}
\|\Delta z_2\|< 4\|B_1^{-1}\|(1+\sup \|\D T_1\|)(|\Delta x_2|+|\Delta y_2|),
\end{equation}
which implies the lemma after letting $K=4\|B_1^{-1}\|(1+\sup \|\D T_1\|)$; estimate
\eqref{dtvxy} follows from \eqref{eq:ucone:12} and \eqref{eqsexp}.\qed

The existence of the center-unstable cone field $\mathcal{C}^{cu}$ implies that the areas of certain surfaces are expanded by the map $T$. We denote by $\mathcal{A}(S)$ the area of a surface $S$.
\begin{lem}\label{lem:areaexpansion} There exists $L>0$ such that
for any surface $S\subset\sigma^0_k$ such that its tangent space at every point lies in the cone field $\mathcal{C}^{cu}$, we have
\begin{equation}\label{eq:expansion:1}
\mathcal{A}(T(S))>L|\lambda\gamma|^k\mathcal{A}(S).
\end{equation}
\end{lem}

\noindent{\it Proof.} Since the cone field $\mathcal{C}^{cu}$ has the same form \eqref{eq:ucone:1} at all points, we have that, for any surface $S$ whose tangent space lies in $\mathcal{C}^{cu}$, its equation takes the form $z=S(x,y)$ and the derivatives $\partial S/\partial x$ and $\partial S/\partial x$ are uniformly bounded away from zero and infinity. Thus, there exist positive constants $L_1$ and $L_2$ such that
\begin{equation}\label{eq:expansion:2}
L_2\mathcal{A}(\pi_0(S))<\mathcal{A}(S) < L_1\mathcal{A}(\pi_0(S)),
\end{equation} 
where $\pi_0$ is the projection onto the $(x,y)$-plane. Since $C^{cu}$ is invariant under $\D T$, the tangent space of $T(S)$ also lies in $C^{cu}$. Therefore, 
\begin{equation}\label{eq:expansion:3}
\mathcal{A}(T(S)) > L_2\mathcal{A}(\pi_0(T(S))).
\end{equation} 
Let $G=\pi_0\circ T|_{z=0}:(x,y)\mapsto (\bar{x},\bar{y})$. We note that
$$\mathcal{A}(\pi_0(T(S))) =  \int_{\pi_0(T(S))} dxdy =  \int_{\pi_0(S)}|\det\D G|dudv $$
and
$$\mathcal{A}(\pi_0(S))=\int_{\pi_0(S)} dudv.$$
Therefore, in order to prove the lemma, it is sufficient to show that there exists $L_3>0$ such that
\begin{equation}\label{eq:expansion:4}
\left|\det\D G\right|>L_3|\lambda\gamma|^k.
\end{equation}

In what follows we prove inequality \eqref{eq:expansion:4}. By \eqref{eq:setting:T^k}, the map $T_0^k|_{z=0}$
is given by
\begin{equation*}\label{eq:expansion:5}
\begin{array}{rcl}
x_k &=& \lambda^k x +  \phi_k(x,y_k,0), \\
y &=& \gamma^{-k} y_k + \psi_k(x,y_k,0),\\
z_k&=&\hat\phi_k(x,y_k,0).
\end{array}
\end{equation*}
By \eqref{fpnorm},\eqref{fpnorm1}, this map can be rewritten as
\begin{equation}\label{eq:expansion:5b0}
\begin{array}{rcl}
x_k &=& \lambda^k x +  \tilde \phi_k(x,y), \\
y_k &=& \gamma^{k} y + \tilde \psi_k(x,y),\\
\end{array}
\end{equation}
where
\begin{equation}\label{5bexp0}\begin{array}{l}
\tilde \phi_k=o(|\lambda|^k), \qquad\partial_x \tilde\phi_k= o(|\lambda|^k), \qquad \partial_y\tilde\phi_k=o(|\lambda\gamma|^k),\\
\tilde \psi_k=o(1), \qquad\partial_x \tilde\psi_k= o(1), \qquad \partial_y\tilde\psi_k=o(|\gamma|^k).\end{array}
\end{equation}
One can also express $z_k$ as a function of $x_k$ and $y_k$ and see that this function satisfies
\begin{equation}\label{zkestz0}
z_k=O(\hat\lambda^k), \qquad \partial_{x_k} z_k=O(\hat\lambda^k{\lambda}^{-k}), \qquad \partial_{y_k} z_k=O(\hat\lambda^k).
\end{equation}

The map $G$ can be written as the composition of the map \eqref{eq:expansion:5b0} and the map
$T_1|_{T^k_0(\{z=0\})}$ which, by \eqref{eq:setting:T_1}, is given by 
\begin{equation*}
\begin{array}{l}
\bar x =x^+ +ax_k+b(y_k - y^-) + a_{13} z_k +h_1(x_k,y_k,z_k), \\ [5pt]
\bar y = \mu + cx_k + d(y_k-y^-)^2 + a_{23} z_k + h_2(x_k,y_k,z_k).
\end{array}
\end{equation*}
The above formulas yield
\begin{equation}\label{eq:expansion:7}
\dfrac{\partial (x_k,y_k)}{(x,y)}=
\begin{pmatrix}
\lambda^k+o(\lambda^k) & o({\lambda}^k\gamma^k)\\
o(1)_{k\to +\infty} & \gamma^k+o(\gamma^k)
\end{pmatrix},
\end{equation}
and
\begin{equation}\label{eq:expansion:8}
\dfrac{\partial (\bar{x},\bar{y})}{(x_k,y_k)}=
\begin{pmatrix}
a+O(|\hat{\lambda}^k|\lambda|^{-k}|+|y_k-y^-|) & b+O(|\lambda|^k+|y_k-y^-|)\\
c+O(|\hat{\lambda}^k|\lambda|^{-k}|+|y_k-y^-|) & 2d(y_k-y^-)+O(|\lambda|^k+(y_k-y^-)^2)
\end{pmatrix}.
\end{equation}
A straightforward computation gives 
$$\left|\det\D G\right| = \left|\det \dfrac{\partial (\bar{x},\bar{y})}{(x_k,y_k)} \det \dfrac{\partial (x_k,y_k)}{(x,y)}\right|=
 |bc + O(y_k-y^-)|\; |\lambda\gamma|^k+o(|\lambda\gamma|^k).$$
The term $y_k-y^-$ is bounded by the small number $\delta$ (the size of $\Pi_0$ and $\Pi_1$), and $bc\neq 0$ by \eqref{eq:nc}. It follows that \eqref{eq:expansion:4} holds indeed for some $L_3>0$. \qed

We proceed to find a stable cone field.

\begin{lem}\label{lem:scone}
There exists a stable cone field $\mathcal{C}^{s}$ over $\Sigma^0\cap T(\Sigma^0)$ which is strictly backward-invariant under $\D T$. 
The cone at the point $X\in\sigma^0_k\cap T(\Sigma^0)$ is given by
\begin{equation}\label{eq:scone:4}
\mathcal{C}^s(X)=\{(\Delta x, \Delta y, \Delta z)\mid |\Delta x|\leqslant K_1 \hat{\lambda}^{k}|\lambda|^{-{k}} \|\Delta z\|\ ,|\Delta y| 
\leqslant K_2 \hat{\lambda}^{k}|\gamma|^{-{k}} \|\Delta z\|\},
\end{equation}
where $K_1$ and $K_2$ are some positive constants, independent of $X$. The restriction of 
$\D T$ to $\mathcal{C}^s$ is contracting, i.e., there exists $M>0$ such that
\begin{equation}\label{dtvxyst}
\|\D T(X) V\| \leq M \hat\lambda^k \|V\| 
\end{equation}
for any $V\in \mathcal{C}^{s}(X)$.

\end{lem}

\noindent{\it Proof.}
Let $Y=T_1\circ T^k_0(X)$ and let $V_2=(\Delta x_2, \Delta y_2, \Delta z_2)$ be a vector in the tangent space at $Y$ such that
\begin{equation}\label{eq:scone:4a}
|\Delta x_2|\leqslant S \|\Delta z_2\|\,\,\mathrm{and}\,\,|\Delta y_2| \leqslant S \|\Delta z_2\|,
\end{equation}
where $S>0$ is a constant. Denote $\D T^{-1}_1(Y)V_2=(\Delta x_1, \Delta y_1, \Delta z_1)$. From \eqref{eq:setting:T_1}, we have
\begin{equation}\label{eq:scone:5}
\begin{array}{rcl}
(\Delta x_2, \Delta y_2)&=& B_1 (\Delta x_1, \Delta y_1)+B_2 \Delta z_1, \\[5pt]
\Delta z_2&=& B_3 (\Delta x_1, \Delta y_1)+B_4 \Delta z_1,
\end{array}
\end{equation}
where $B_i$ $(i=1\dots 4)$ are some matrices whose norms are uniformly bounded. Note that $B_1$ is close to
$\left(\begin{array}{cc} a & b \\ c & 0\end{array}\right)$ and $bc\neq 0$ by \eqref{eq:nc}, so the matrix $B_1$ is invertible.

Now, equation \eqref{eq:scone:5} can be rewritten as
\begin{equation}\label{eq:scone:7}
\begin{array}{rcl}
(\Delta x_1, \Delta y_1)&=& B_1^{-1}(\Delta x_2, \Delta y_2)-B_1^{-1}B_2 \Delta z_1, \\[5pt]
\Delta z_2& =& B_3 B_1^{-1}  (\Delta x_2,\Delta y_2) + (B_4-B_1^{-1} B_2) \Delta z_1.
\end{array}
\end{equation}
By choosing $S$ such that $S\|B_2 B_1^{-1}\|<1$, we obtain 
$$ \| B_3 B_1^{-1}(\Delta x_2,\Delta y_2)\| \leqslant S\|B_3 B_1^{-1}\| \, \|\Delta z_2\|< \hat S \|\Delta z_2\|,$$
where $\hat S<1$ is a constant, independent of the choice of the point $Y$ and the vector $V_2$.

Hence, the second equation in \eqref{eq:scone:7} implies 
$\|\Delta z_2\|=O(\|\Delta z_1\|)$, which by \eqref{eq:scone:4a} further implies
\begin{equation}\label{eq:scone:7a}
|\Delta x_2|=O(\|\Delta z_1\|) \quad\mathrm{and}\quad  |\Delta y_2|=O(\|\Delta z_1\|).
\end{equation}
Finally, from the first equation in \eqref{eq:scone:7} we find
\begin{equation}\label{eq:scone:8a}
|\Delta x_1| +|\Delta y_1| \leqslant l\|\Delta z_1\|,
\end{equation}
where $l$ is some positive constant, independent of the choice of $Y$ and $V_2$.

Denote $\D T^{-k}_0 \D T^{-1}_1(Y)V_2=(\Delta x_0, \Delta y_0, \Delta z_0)$. By formula \eqref{eq:setting:T^k}, noting that the first derivatives
of $\phi,\hat{\phi}$ and $\psi$ are bounded, we have the following relations:
\begin{eqnarray}
\Delta x_1& =& \lambda^{k} \Delta x_0 + o(\lambda^{k})(\Delta x_0+\Delta y_1 + \Delta z_0), \label{eq:scone:9}\\[5pt]
\Delta y_0& =& \gamma^{-{k}} \Delta y_1 + o({\gamma^{-{k}}})(\Delta x_0+\Delta y_1 + \Delta z_0), \label{eq:scone:10}\\[5pt]
\Delta z_1 &= &O(\hat{\lambda}^{k})(\Delta x_0+\Delta y_1 + \Delta z_0).\label{eq:scone:11}
\end{eqnarray} 
Estimates \eqref{eq:scone:11} and \eqref{eq:scone:8a} give
\begin{eqnarray*}\label{eq:scone:12}
\|\Delta z_1\| &=& O(\hat{\lambda}^{k})(|\Delta x_0|+\|\Delta z_0\|),\\
\|\Delta y_1\| &=& O(\hat{\lambda}^{k})(|\Delta x_0|+\|\Delta z_0\|).
\end{eqnarray*}
With these estimate and \eqref{eq:scone:8a}, equation \eqref{eq:scone:9} yields
$$|\Delta x_0| = O(\lambda^{-k}) \|\Delta z_1\|+o(\|\Delta z_0\|).$$
By plugging the above equation into \eqref{eq:scone:11}, we obtain
\begin{equation}\label{eq:scone:14}
|\Delta z_1| = O(\hat{\lambda}^{k})  \|\Delta z_0\|,
\end{equation}
which, along with \eqref{eq:scone:8a}, further implies
\begin{equation}\label{eq:scone:15}
|\Delta x_1| + |\Delta y_1| = O(\hat{\lambda}^{k})  \|\Delta z_0\|.
\end{equation}
Finally, the above equation together with \eqref{eq:scone:9} and \eqref{eq:scone:10} leads to
\begin{equation*}\label{eq:scone:15a}
\begin{array}{rcl}
|\Delta x_0|= O(\hat{\lambda}^{k}\lambda^{-k})  \|\Delta z_0\|, \\
|\Delta y_0|= O(\hat{\lambda}^{k}\gamma^{-k})  \|\Delta z_0\|.
\end{array}
\end{equation*}

This formula shows that the image by $\D T^{-1}$ of a vector satisfying \eqref{eq:scone:4a} lies in the cone \eqref{eq:scone:4}.
If $k^*$ was taken sufficently large, then every vector from the cone \eqref{eq:scone:4} satisfies \eqref{eq:scone:4a}, i.e.,
we have proven the required invariance of the cone field \eqref{eq:scone:4}. Estimate
\eqref{dtvxyst} follows from \eqref{eq:scone:14}, \eqref{eq:scone:15} and the uniform
boundedness of $\D T$. \qed

The strong-stable foliation $\mathcal{F}_0$ which exists in the stable manifold $W^s(O)$ extends to an invariant foliation $\mathcal{F}^s$
in a small neighborhood of the homoclinic cycle $O\cup \Gamma\cup \tilde \Gamma$ we consider here (see \cite{tu96}). As the tangents
to the leaves of the invariant foliation $\mathcal{F}^s$ must lie in the stable invariant cone $\mathcal{C}^{s}$, Lemma \ref{lem:scone} 
immediately implies the following formula for the leaves of $\mathcal{F}^s$. 

\begin{lem}\label{lem:ssleavesweak}
The leaf of the strong-stable foliation $\mathcal{F}^s$ through a point $(x^*,y^*,z^*)\in\Sigma^0$ with a stay number $k$ takes the form
\begin{equation}\label{eq:ssleavesweak}
\begin{array}{rcl}
x&=&x^*+\varphi_1(z; x^*,y^*,z^*), \\[5pt]
y&=&y^*+\varphi_2(z; x^*,y^*,z^*),
\end{array}
\end{equation}
where 
$$\begin{array}{l}\displaystyle
\varphi_1=O({\hat{\lambda}}^k\lambda^{-k}), \qquad \frac{\partial\varphi_1}{\partial z}=O({\hat{\lambda}}^k\lambda^{-k}),\\[10pt]
\displaystyle
\varphi_2=O(\hat{\lambda}^{k}\gamma^{-k}), \qquad \frac{\partial\varphi_1}{\partial z}= O(\hat{\lambda}^{k}\gamma^{-k}).
\end{array}$$
\end{lem}
Note that we do not estimate the derivatives of $\varphi_{1,2}$ with respect to $(x^*,y^*,z^*)$ here.\\

In the proof of  Lemma \ref{lem:scone}, we have not used condition (C4) on the multipliers of $O$. Formula \eqref{eq:ssleavesweak} will only be helpful in the non-symmetric case (Theorem \ref{thm2}) where we have more parameters to do the bifurcation. When it comes to the symmetric case (Theorem \ref{thm1}),
we need a better estimate, which will be obtained by taking into account condition (C4).
\begin{lem}\label{lem:ssleaves}
If condition (C4) is satisfied, then the strong-stable leaf through a point $(x^*,y^*,z^*)\in\Sigma^0$ with a stay number $k$ assumes the same form as in \eqref{eq:ssleavesweak}, but the function $\varphi_1$ now satisfies
\begin{equation}
\varphi_1=O({\lambda}_0^k\lambda^{-k}), \qquad \frac{\partial\varphi_1}{\partial z}= O({\lambda}_0^k\lambda^{-k}),
\end{equation}
where $\lambda_0$ can be taken arbitrarily close to $|\lambda_1|$.
\end{lem}

\noindent {\it Proof.}
Take any point $X \in\sigma_k^0$ and consider a vector $(\Delta x_0,\Delta y_0,\Delta z_0)$ in the tangent space, at $X$, to 
the leaf of the invariant foliation $\mathcal{F}^s$ through $X$. We need to show
that 
\begin{equation}\label{sameestk}
|\Delta x| \leq K {\lambda}_0^k\lambda^{-k} \|\Delta z\|
\end{equation}
for some constant $K$, independent of $X$.

Let $(x_k,y_k,z_k)=T^k_0 X$ and $(\Delta x_1,\Delta y_1,\Delta z_1)=\D T_0^k(\Delta x_0,\Delta y_0,\Delta z_0)$. 
By formula \eqref{eq:setting:T^k}, we have
\begin{equation}\label{eq:ssleaves:4}
\begin{array}{rcl}
\Delta x_1 &=& \lambda^k(1+\dots) \Delta x_0 + o(\lambda^k) \Delta y_1 + \dfrac{\partial x_k}{\partial z_0} \Delta z_0,\\[10pt]
\Delta z_1 &=& \dfrac{\partial z_k}{\partial z_0} \Delta z_0 + o(\lambda^k) \Delta x_0 + o(\lambda^k) \Delta y_1,
\end{array}
\end{equation}
where the dots denote terms that tend to zero as $k\to+\infty$. Since the vector $(\Delta x_0,\Delta y_0,\Delta z_0)$ is in the stable cone $C^{s}$, its image
$V$ by $\D (T_1 T_0^k)$ is also in $C^s$. So, as we have shown in the proof of Lemma \ref{lem:scone}, the vector
$(\Delta x_1, \Delta y_1, \Delta z_1)= \D T_1^{-1} V$ must satisfy
$$|\Delta x_1|+|\Delta y_1|= O(\|\Delta z_1\|),$$
see (\ref{eq:scone:8a}). Plugging this into \eqref{eq:ssleaves:4} gives
$$\lambda^k(1+\dots) \Delta x_0 = O(\|\dfrac{\partial x_k}{\partial z_0}\| + \|\dfrac{\partial z_k}{\partial z_0}\|) \Delta z_0.$$
By Lemma \ref{lem:firstd}, this inequality implies \eqref{sameestk}. \qed

\section{The index-2 condition}\label{sec:index2}
In this section we find a condition which ensures that a period-2 point of $T$ is a saddle of index 2. We will start with a result describing the multipliers 
of a periodic point (in a more general case where period-$n$ orbits are considered).

Let $X\in\Sigma^0$ be a period-$n$ point such that $X=T^n(X)=T_1\circ T_0^{k_n}\circ T_1\circ T_0^{k_{n-1}}\circ\dots\circ T_1\circ T^{k_1}(X)$, where $k_1,\dots,k_n$ are the corresponding stay numbers. We sort the eigenvalues of $\D T^{n}$ the multipliers of $X$ in decreasing order by their absolute values
and denote them as $\nu_1,\dots,\nu_{D}$. By Lemmas \eqref{lem:ucone} and \eqref{lem:scone}, the derivative  $\D T^n$ at $X$ has a pair of invariant cones,
which implies the existence of a two-dimensional invariant subspace $E^{cu}$ (in the center-unstable cone) and a $(D-2)$-dimensional invariant subspace $E^s$
in the stable cone. Estimates \eqref{dtvxy} and \eqref{dtvxyst} for $\D T^n$ restricted to $E^{cu}$ and, respectively, $E^s$ immediately give the the following estimate on the multipliers of $X$.

\begin{lem}\label{lem:eigenvalues}
The eigenvalues of $\D T^{n}|_{E^{cu}}$ are $\nu_1$ and $\nu_2$, and the eigenvalues of $\D T^{n}|_{E^s}$ are $\nu_3,\dots,\nu_D$. Moreover, we have
\begin{equation}\label{eq:eigenvalues:1}
|\nu_i|^{-1}=O(|\lambda|^{k_1+\dots+k_n}),\quad i=1,2,
\end{equation}
and
\begin{equation}\label{eq:eigenvalues:1a}
|\nu_i|=O(\hat{\lambda}^{k_1+\dots+k_n}),\quad i=3,4,\dots,D.
\end{equation}
\end{lem}

We now consider orbits of period 2, and find the condition under which such point is an index-2 saddle, i.e., $|\nu_1|>1$ and $|\nu_2|>1$.
Let $Q\in\Pi_0$ be a period-2 point of $T$ with stay numbers $k$ and $m$. Denote $Q_{01}=Q=(x_{01},y_{01},z_{01})$, $Q_{11}=T_0^k(Q)=(x_{11},y_{11},z_{11})$, $Q_{02}=T_1\circ T_0^k(Q)=(x_{02},y_{02},z_{02})$ and
$Q_{12}=T_0^m\circ T_1\circ T_0^k(Q)=(x_{12},y_{12},z_{12})$.
\begin{lem}\label{lem:index2} There exist functions $r_{1,2,3,4}$, which depends on the integers $m$ and $k$, parameters and
the coordinates of the points $Q_{ij}$, such that the point $Q$ is a saddle of index 2 if and only if there exists some number $s\in(-1,1)$ such that
\begin{equation}\label{eq:index2:1}
(y_{11}-y^-+r_1)(y_{12}-y^-+r_2)=r_3 + r_4 s.
\end{equation}
The functions $r_{1,2,3,4}$ satisfy
\begin{equation}\label{r123emk+}\begin{array}{l}
r_1=O((y_{11}-y^-)^2+|\lambda|^k+|\gamma|^{-m}), \qquad r_2=O((y_{12}-y^-)^2+|\lambda|^m+|\gamma|^{-k}),\\ 
r_3=O(|\lambda|^k|\gamma|^{-k}|+|\lambda|^m|\gamma|^{-m}+|\lambda|^{(k+m)}), \qquad r_4=O(\lambda^{(k+m)}).
\end{array}
\end{equation}
\end{lem}

\noindent{\it Proof.}
One can check that the condition $|\nu_1|,|\nu_2|>1$ is equivalent to 
$$|\nu_1\nu_2|>1\quad\mbox{and}\quad\dfrac{\nu_1+\nu_2}{\nu_1\nu_2+1}=s, \qquad -1<s<1.$$
This can be written as
\begin{equation}\label{eq:index2:4}
|\det \D T^2|_{E^{cu}}| > 1 \quad\mbox{and}\quad \dfrac{\tr \D T^2|_{E^{cu}}}{\det \D T^2|_{E^{cu}}+1}=s, \qquad -1<s<1,
\end{equation}
where $E^{cu}$ is the two-dimensional invariant subspace introduced before Lemma \ref{lem:eigenvalues}. In what follows, we use $(\Delta x,\Delta y, \Delta z)$
to denote 
a vector in $E^{cu}$. Note that $\Delta z$ is a function of $\Delta x$ and $\Delta y$. We, thus, need to compute the trace and the determinant of 
$\D T^2|_{E^{cu}}:(\Delta x,\Delta y)\mapsto (\Delta \bar{x},\Delta \bar{y})$.

Denote
\begin{equation}\label{eq:index2:5}
\eta_1=y_{11}-y^- \quad \mbox{and} \quad \eta_2=y_{12}-y^-.
\end{equation}
Take a vector $V=(\Delta x_1,\Delta y_1,\Delta z_1)\in E^{cu}$. Formula \eqref{eq:setting:T^k} implies that 
\begin{equation}\label{eq:index2:6a}
\D T_0^k|_{E^{cu}}V=A_1
\begin{pmatrix}
\Delta x_1\\
\Delta y_1
\end{pmatrix} 
=
\begin{pmatrix}
\lambda^k+o(\lambda^k) & o({\lambda}^k\gamma^k)\\
o(1)_{k\to +\infty} & \gamma^k+o(\gamma^k)
\end{pmatrix}
\begin{pmatrix}
\Delta x_1\\
\Delta y_1
\end{pmatrix} 
=:
\begin{pmatrix}
\Delta x_2\\
\Delta y_2
\end{pmatrix}.
\end{equation}
Note that the $\Delta z_1$ component is a bounded function of $(\Delta x_1,\Delta y_1)$ and its contribution to $\Delta x_2$ and $\Delta y_2$
goes into the small terms in $A_1$.

After noting $x_{11}=O(\lambda^k)$ and $z_{11}=O(\hat{\lambda}^k)$ from \eqref{eq:setting:T^k},
we can write the matrix $\D T_1(Q_{11})$  as
\begin{equation}\label{eq:index2:6}
\left(
\begin{array}{ccc}
a+O(|\lambda|^k+|\eta_1|) & b+O(|\lambda|^k+|\eta_1|) & a_{13}+O(|\lambda|^k+|\eta_1|)  \\
c+O(|\lambda|^k+|\eta_1|) & 2d\eta_1+O(|\lambda|^k+\eta_1^2) & a_{23}+O(|\lambda|^k+|\eta_1|)  \\
a_{31}+O(|\lambda|^k+|\eta_1|) & a_{32}+O(|\lambda|^k+|\eta_1|) & a_{33}+O(|\lambda|^k+|\eta_1|)
\end{array}
\right).
\end{equation}
Since the vector $\D T^k_0 V_1=(\Delta x_2, \Delta y_2,\Delta z_2)$ belongs to $\D T^k_0 \mathcal{C}^{cu}$, we 
have from equation \eqref{eq:ucone:8} that
$$\Delta z_1 = O(\hat{\lambda}^{k}\lambda^{-k})\Delta x_1 + O(\hat{\lambda}^{k})\Delta y_1.$$
Along with \eqref{eq:index2:6}, this leads to
\begin{equation}\label{eq:index2:8}
\D (T_1\circ T^k_0)|_{E^{cu}} V_1=
A_2
\begin{pmatrix}
\Delta x_2\\
\Delta y_2
\end{pmatrix} 
=
\begin{pmatrix}
a+O(|\hat{\lambda}^k\lambda^{-k}|+|\eta_1|) & b+O(|\lambda|^k+|\eta_1|)\\
c+O(|\hat{\lambda}^k\lambda^{-k}|+|\eta_1|) & 2d\eta_1+O(|\lambda|^k+\eta_1^2)
\end{pmatrix}
\begin{pmatrix}
\Delta x_2\\
\Delta y_2
\end{pmatrix} 
=:
\begin{pmatrix}
\Delta x_3\\
\Delta y_3
\end{pmatrix},
\end{equation}
where the contribution of $\Delta z_2$ goes into the $O(\cdot)$ terms.

By repeating the same procedure, we also obtain the following formulas 
for $\D T_0^m|_{\D(T_1\circ T^k_0){E^{cu}}}$ and $\D T_1|_{\D(T_0^m\circ T_1\circ T^k_0){E^{cu}}}$:
\begin{equation}\label{eq:index2:9}
\D T_0^m|_{\D(T_1\circ T^k_0){E^{cu}}}=A_3=
\begin{pmatrix}
\lambda^m+o(\lambda^m) & o({\lambda}^m\gamma^m)\\
o(1)_{m\to +\infty} & \gamma^m+o(\gamma^m)
\end{pmatrix},
\end{equation}
and
\begin{equation}\label{eq:index2:10}
\D T_1|_{\D(T_0^m\circ T_1\circ T^k_0){E^{cu}}}=A_4=
\begin{pmatrix}
a+O(|\hat{\lambda}^m\lambda^{-m}|+|\eta_2|) & b+O(|\lambda|^m+|\eta_2|)\\
c+O(|\hat{\lambda}^m\lambda^{-m}|+|\eta_2|) & 2d\eta_2+O(|\lambda|^m+\eta_2^2)
\end{pmatrix}.
\end{equation}

Now we can write the map $\D T^2_{E^{cu}}$ as the product $A_4A_3A_2A_1$. By equations \eqref{eq:index2:6a} and 
\eqref{eq:index2:8} - \eqref{eq:index2:10}, we have 
\begin{equation}\label{eq:index2:11}
A_2A_1=
\left(
\begin{array}{cc}
o(1)_{k\to+\infty}  & b\gamma^k+o(\gamma^k)\\
c\lambda^k+o(|\lambda|^k+|\eta_1|)   &  \gamma^k (2d \eta_1+o(1)_{k\to+\infty} \eta_1 +O(|\lambda|^k+\eta_1^2))
\end{array}
\right),
\end{equation}
\begin{equation}\label{eq:index2:12}
A_4A_3=
\left(
\begin{array}{cc}
o(1)_{m\to+\infty}  & b\gamma^m+o(\gamma^m)\\
c\lambda^m+o(|\lambda|^m+|\eta_2|)   &  \gamma^m (2d\eta_2+o(1)_{m\to+\infty} \eta_2 + O(|\lambda|^m+\eta_2^2))
\end{array}
\right),
\end{equation}
which yields
$$\begin{array}{l}\tr \D T^2_{E^{cu}}= \tr (A_4A_3A_2A_1)=\\
 \qquad =\gamma^{k+m}(4d^2(\eta_1+O(\eta_1^2)(\eta_2+O(\eta_2^2))(1+\dots) +
\eta_1 O(|\lambda|^m+|\gamma|^{-k}) + \eta_2 O(|\lambda|^k+|\gamma|^{-m})+\\
\qquad\qquad\qquad\qquad\qquad\qquad  +
bc\lambda^m\gamma^{-m}(1+\dots)+bc\lambda^k\gamma^{-k}(1+\dots)),
\end{array}
$$
where the dots stand for terms that tend to zero as $m,k\to+\infty$. This equation can be rewritten as
\begin{equation}\label{eq:index2:14}
\begin{array}{rcl}
\tr \D T^2_{E^{cu}}&=&\gamma^{k+m}(1+\dots)(4d^2(\eta_1+O(\eta_1^2+|\lambda|^k+|\gamma|^{-m}))\\
&&\cdot
(\eta_2+O(\eta_2^2+|\lambda|^m+|\gamma|^{-k}) + O(|\lambda|^m|\gamma|^{-m}+|\lambda|^k|\gamma|^{-k}+|\lambda|^{k+m})).
\end{array}
\end{equation}

It follows immediately from \eqref{eq:index2:6a}, \eqref{eq:index2:8} and \eqref{eq:index2:9},\eqref{eq:index2:10} that
\begin{equation}\label{eq:index2:15}
\begin{array}{rcl}
\det{A_2 A_1}&=&-\lambda^k\gamma^k(bc+O(\eta_1)+o(1)_{k\to+\infty}),\\
\det{A_4 A_3}&=&-\lambda^m\gamma^m(bc+O(\eta_2)+o(1)_{m\to+\infty}).
\end{array}
\end{equation}
Consequently, with the fact $bc\neq 0$ by \eqref{eq:nc}, we obtain
\begin{equation}\label{eq:index2:16}
\det \D T^2_{E^{cu}}=  (\lambda\gamma)^{k+m} (bc)^2 (1+ O(|\eta_1|+|\eta_2|) + o(1)_{k,m\to+\infty}),
\end{equation}
and, since $|\lambda\gamma|>1,$
$$|\det \D T^2_{E^{cu}}|>1.$$
Therefore, by \eqref{eq:index2:14} and \eqref{eq:index2:16}, condition \eqref{eq:index2:4} is indeed equivalent to
\eqref{eq:index2:1} and \eqref{r123emk+}.
\qed

\section{Proofs of Theorems \ref{thm1} and \ref{thm2}}\label{sec:proofofthm}
We first prove Theorem \ref{thm1}. It will be proved in two steps corresponding to finding the orbits of transverse and non-transverse heteroclinic
intersections in a heterodimensional cycle. 
The proof of Theorem \ref{thm2} will be a modification of that of Theorem \ref{thm1}.

\subsection{Proof of Theorem \ref{thm1}}\label{sec:proofofthm1}
Theorem \ref{thm1} is a consequence of the following two lemmas. Recall that $\delta$ is the size of the neighbourhood $\Pi_1$ of $M^{-}$. 
\begin{lem}\label{lem:transverse}
Let $F$ satisfy conditions (C1)-(C3). If there exists two transverse homoclinic points $N_1,N_2\in W_{loc}^u(O)$ of $O$ satisfying $0<y^- - y_{_{N_1}}<\delta/2$ and $0<y_{_{N_2}} - y^- < \delta/2$, then we can find an integer $K$ such that, for any index-2 periodic point $Q$ of $F$ whose orbit lies in $\begin{matrix}\bigcup^{+\infty}_K \sigma^0_k \end{matrix}$, the intersection $W^u(Q)\cap W^s(O)$ is non-empty. The result also holds for all 
diffeomorphisms sufficiently $C^2$-close to $F$.
\end{lem}

\begin{lem}\label{lem:nontransverse}
Consider a two-parameter family $\{F_{\mu,\theta}\}$ of diffeomorphisms in $\mbox{Diff}\,^r_s(\mathcal{M^D})$ where $F_{0,\theta^*}$ 
satisfies conditions (C1) - (C4). If $cx^+y^->0$ and $cdx^+>0$, then, for any sequence $\{(k_j,m_j)\}$ of pairs of even natural numbers satisfying $k_j,m_j\to +\infty$ and $m_j/k_j \to \theta^*$ as $j\to +\infty$, there exists a sequence $\{(\mu_j,\theta_j)\}$ accumulating on $(0,\theta^*)$ such that, for any sufficiently large $j$, the diffeomorphism $F_{\mu_j,\theta_j}$ has an index-2 periodic orbit $Q_j$ satisfying $T_1\circ T_0^{m_j}\circ T_1\circ T_0^{k_j}(Q_j)=Q_j$ and $W^s(Q_j)\cap W^u(O)\neq\emptyset$.
\end{lem}

Theorem \ref{thm1} follows from these lemmas.

\noindent{\it Proof of Theorem \ref{thm1}.} Lemma \ref{lem:preperturb2} gives us a sequence $\{\mu_i\}$ accumulating on $\mu=0$ such that 
$F_{\mu_i,\theta^*}$ has a new orbit $\Gamma_i$ of homoclinic tangency to $O$. This orbit $\Gamma_i$ has a point 
$M_i=(0,y_i,0)\in W^{u}_{loc}(O)\cap\Pi_1$ accompanied by two transverse homoclinic points $N^1_i=(0,y^1_i,0)$ and $N^2_i=(0,y^2_i,0)$ such that 
$0<y_i-y^1_i<\delta/2$ and $0<y^2_i-y_i<\delta/2$. It follows that $F_{\mu_i,\theta^*}$ has the property given by Lemma \ref{lem:transverse}.

Next, we fix a sufficiently large $i$. According to Lemma \ref{lem:preperturb2}, the global map associated to $\Gamma_j$ has $cx^+y^->0$ and $cdx^+>0$. Obviously, $F_{\mu_i,\theta^*}$ with a sufficiently large $i$ fulfils conditions (C1) - (C4). Hence, Lemma \ref{lem:nontransverse} gives a sequence 
$\{(\mu^n_i,\theta^n_i)\}_n$ accumulating on $(\mu_i,\theta^*)$ such that the system $F_{\mu^n_i,\theta^n_i}$ has an index-2 periodic point $Q^n_i$ 
satisfying $T_1\circ T_0^{m(n,i)}\circ T_1\circ T_0^{k(n,i)}(Q^n_i)=Q^n_i$ and $W^s(Q^n_i)\cap W^u(O)\neq\emptyset$, where $T_0$ and $T_1$ are the local 
and global maps of $F_{\mu^n_i,\theta^n_i}$. Since Lemma \ref{lem:transverse} holds for $F_{\mu_i,\theta^*}$ and all sufficiently $C^2$-close diffeomorphisms, 
the theorem follows by taking $(\mu_j,\theta_j)=(\mu^{n_j}_{i_j},\theta^{n_j}_{i_j})$, where $\{n_j\}$ and $\{i_j\}$ are any sequences tending to positive infinity as $j\to +\infty$. \qed

We proceed to prove Lemmas \ref{lem:transverse} and \ref{lem:nontransverse}.

\noindent{\it Proof of Lemma \ref{lem:transverse}.} 
We will prove this lemma by using the fact that the map $T$ expands two-dimensional areas, which follows from the assumption $|\lambda\gamma|>1$.

Let us first define a quotient first-return map by the leaves of the invariant foliation $\mathcal{F}^s$. 
Recall that the first return map $T:\Sigma^0\to\Pi_0$ (where $\Sigma^0=\begin{matrix}
\bigcup^{+ \infty}_{k^*} \sigma^0_k
\end{matrix}$) takes the form $T(M)=T_1\circ T^k_0(X)$ for any $M\in\sigma^0_k$ (see \eqref{eq:setting:T}). Let $\pi:U_0\to\{z=0\}$ be the projection map along the leaves of $\mathcal{F}^s$. Denote by $\hat{\Pi}_i, \hat{\sigma}^0_k$ and $\hat{\Sigma}^0$ the intersections of $\Pi_i,{\sigma}^0_k$ and $\Sigma^0$ with $\{z=0\}$. The foliation $\mathcal{F}^s$ induces the quotient map from $\hat{\Sigma}^0$ to $\hat{\Pi}_0$:
\begin{equation*}\label{eq:transverse:1}
\hat{T}(M)= \pi\circ T_1\circ T^k_0(M),
\end{equation*}
for any $M\in\hat{\sigma}_k^0$. 

Consider any surface $S_k\in\sigma^0_k$ whose tangents lie in the center-unstable cone field $\mathcal{C}^{cu}$. This surface is transverse to $\mathcal{F}^s$ and the angle between them are uniformly bounded. Therefore, by the absolute continuity of $\mathcal{F}^s$, there exist constants $q_1$ and $q_2$
 which do not depend on the surface such that
$$q_1\mathcal{A}(S_k)<\mathcal{A}(\pi(S_k))<q_2\mathcal{A}(S_k),$$
where we use $\mathcal{A}(\cdot)$ to denote the area. On the other hand, Lemma \ref{lem:areaexpansion} gives
\begin{equation}
\mathcal{A}(T(S_k))>L|\lambda\gamma|^k\mathcal{A}(S_k),
\end{equation}
where $L$ is some positive constant. It follows that
$$\mathcal{A}(\pi\circ T(S_k))>q_1\mathcal{A}(T(S_k))>q_1L|\lambda\gamma|^k\mathcal{A}(S_k)>q_1q_2^{-1}L|\lambda\gamma|^k\mathcal{A}(\pi(S_k)).$$
Thus, there exists $k'$ such that for any $k>k'$ we have
\begin{equation}\label{eq:transverse:5}
\mathcal{A}(\hat{T}(S_k))>q \mathcal{A}(\pi(S_k)),
\end{equation}
for some $q>1$.

Let $K=\max(k^*,k')$ and $Q\in\sigma^0_{k_0}$ ($k_0>K$) be any index-2 periodic point of $T$.
Take any small piece $W^u$ of the unstable manifold of $Q$. The tangent space of $W^u$ lies in the cone field $\mathcal{C}^{cu}$. Inequality \eqref{eq:transverse:5} implies that $\mathcal{A}(\pi(W^u))$ increases after every iteration by $\hat{T}$. This means that one can find $n_0$ such that $T^{n}(W^u)\in \sigma^0_{k_0}$ for all $n<n_0$ and $T^{n_0}(W^u)$ insects one of the boundaries $v_1=\{x=x^+-\delta/2\},\,\,v_2=\{x=x^++\delta/2\},\,\,h_1=\{y=\gamma^{-k_0}(y^- -\delta/2)\}$ and $h_2=\{y=\gamma^{-k_0}(y^- +\delta/2)\}$ of $\sigma^0_{k_0}$. We claim that $T^{n_0}(W^u)$ intersects either $h_1$ or $h_2$. For that, we show that $T^{n_0}(W^u)$ cannot intersects $v_1$ and $v_2$. Indeed, formula \eqref{eq:setting:T^k} for the local map implies that $x$ and $z$ in \eqref{eq:setting:T_1} are of order of $\lambda^{k_0}$. Hence, the main contribution to the $x$-coordinate in \eqref{eq:setting:T_1} is given by the term $b(y-y^+)$, which is of order $\delta_y$ (recall that we let $\Pi_1=\{(x,y,z)\mid |x|<\delta,|y-y^-|<\delta_y,\|z\|<\delta \}$). It follows that, by taking $K$ sufficiently large and $\delta_y$ sufficiently small, the image $T_1\circ T^{k_0}_0(\sigma^0_{k_0})$ intersects neither $v_1$ nor $v_2$. The claim is proven.

We now take a special choice of the boundaries $h_1$ and $h_2$. Let $y=w_1(x,z)$ and $y=w_2(x,z)$ be the equations of the two pieces of $W^s(O)$ that go through the transverse homoclinic points $N_1$ and $N_2$, respectively. We replace $\Pi_1$ by its subset $\{(x,y,z)\in \Pi_1\mid w_1(x,z)<y<w_2(x,z)\}$. 
Then, all the `horizontal' boundaries of $\sigma^0_k$ are pieces of $W^s(O)$. Lemma \ref{lem:transverse} follows by noticing that $h_1$ and $h_2$ are such boundaries. 

The above computation goes through in the coordinate system where the local map $T_0$ assumes the form \eqref{eq:setting:T_0} and satisfies the identities in \eqref{eq:nonlinearterms}. This can be achieved when $F$ has at least $C^2$-smoothness. Therefore, the above result holds for any diffeomorphism sufficiently $C^2$-close to $F$. \qed

\noindent{\it Proof of Lemma \ref{lem:nontransverse}.}
We start with finding a periodic point $Q\in\Pi_0$ of period 2 and index 2. We are searching for a point $Q$ such that $T^2(Q)=T_1\circ T_0^m\circ T_1\circ T_0^k(Q)=Q$. Let $Q_{01}=Q=(x_{01},y_{01},z_{01}),Q_{11}=T_0^k(Q)=(x_{11},y_{11},z_{11}),Q_{02}=T_1\circ T_0^k(Q)=(x_{02},y_{02},z_{02})$ and $Q_{12}=T_0^m\circ T_1\circ T_0^k(Q)=(x_{12},y_{12},z_{12})$. Recall that $|\lambda\gamma|>1$, hence $\theta=-\ln|\lambda|/\ln|\gamma| <1$. Therefore,
the condition $m/k \to \theta^*$ implies $k - m \gg 0$.

By formulas \eqref{eq:setting:T^k} and \eqref{eq:setting:T_1}, the point $Q$ is a period-2
point if
\begin{eqnarray*}
\left.
\begin{array}{rcl}
x_{11}&=&\lambda^k x_{01} + \phi_k,\\[5pt]
y_{01}&=&\gamma^{-k} y_{11} +\psi_k,\\[5pt]
z_{11}&=&\hat{\phi}_k,
\end{array}
\right.
\quad
\left.
\begin{array}{rcl}
x_{02}-x^+ &=& ax_{11}+b(y_{11} - y^-) + a_{13} z_{11} +h_1, \\ [5pt]
y_{02} &=& \mu + cx_{11} + d(y_{11}-y^-)^2 + a_{23} z_{11} + h_2(,\\[5pt]
z_{02}-z^+ &=& a_{31}x_{11}+a_{32}(y_{11} - y^-) + a_{33} z_{11} +h_3,
\end{array}
\right.\nonumber\\[20pt]
\left.
\begin{array}{rcl}
x_{12}&=&\lambda^m x_{02} + \phi_k,\\[5pt]
y_{02}&=&\gamma^{-m} y_{12} + \psi_k,\\[5pt]
z_{12}&=&\hat{\phi}_k,
\end{array}
\right.
\quad
\left.
\begin{array}{rcl}
x_{01}-x^+ &=& ax_{12}+b(y_{12} - y^-) + a_{13} z_{12} +h_1, \\ [5pt]
y_{01} &=& \mu + cx_{12} + d(y_{12}-y^-)^2 + a_{23} z_{12} + h_2,\\[5pt]
z_{01}-z^+ &=& a_{31}x_{12}+a_{32}(y_{12} - y^-) + a_{33} z_{12} +h_3,
\end{array}
\right. \label{eq:nontr:0.1}
\end{eqnarray*}
which can be rewritten as
\begin{eqnarray}
x_{01}-x^+&=&a\lambda^m x_{02}+b(y_{12}-y^-)+ o({\lambda}^m) +O((y_{12}-y^-)^2),\nonumber \\[5pt]
\gamma^{-k}y_{11} + o({\gamma}^{-k})&=&\mu+ c\lambda^m x_{02}+ d(y_{12}-y^-)^2+o({\lambda}^m)+h_2(0,y_{12}-y^-,0)\nonumber,\\[5pt]
z_{01}-z^+&=& a_{31}\lambda^m x_{02}+ a_{32}(y_{11}-y^-)+o({\lambda}^m)+O((y_{12}-y^-)^2),\nonumber\\[5pt]
x_{02}-x^+&=&a\lambda^k x_{01}+b(y_{11}-y^-)+ o({\lambda}^k)+O((y_{11}-y^-)^2), \label{eq:nontr:0.2}\\[5pt]
\gamma^{-m}y_{12} + o({\gamma}^{-m})&=&\mu+ c\lambda^k x_{01}+ d(y_{11}-y^-)^2+o({\lambda}^k)+h_2(0,y_{11}-y^-,0),\nonumber\\[5pt]
z_{02}-z^+&=& a_{31}\lambda^k x_{01}+ a_{32}(y_{11}-y^-)+o({\lambda}^k)+O((y_{11}-y^-)^2) \nonumber.
\end{eqnarray}
Note that it follows from the implicit function theorem that, at sufficiently large $k,m$, the variables $x_{01},x_{02},z_{01}$ and $z_{02}$ can be expressed as functions of $y_{11}$ and $y_{12}$. Consequently, we need to consider only the equations for $y_{11}$ and $y_{12}$. By introducing $\eta_1=y_{11}-y^-$ and $\eta_2=y_{12}-y^-$, finding a period-2 point becomes equivalent to solving the following system:
\begin{eqnarray}
\gamma^{-k}(\eta_1+y^-) + o({\gamma}^{-k})&=&\mu+ c\lambda^m x^+ + bc\lambda^m \eta_1+ d\eta_2^2+o({\lambda}^m)+h_2(0,\eta_2,0),\label{eq:nontr:1}\\[5pt]
\gamma^{-m}(\eta_2+y^-) + o({\gamma}^{-m})&=&\mu+ c\lambda^k x^+ +bc\lambda^k\eta_2 + d\eta_1^2+o({\lambda}^k)+h_2(0,\eta_1,0).\label{eq:nontr:2}
\end{eqnarray} 
We will look for solutions $\eta_{1,2}$ which tend to zero as $k,m\to +\infty$.  By Lemma \ref{lem:index2}, the corresponding periodic point is of index-2 if,
for some $s\in(-1,1)$,
\begin{equation}\label{eq:nontr:3}
(\eta_1+O(\eta_1^2+\lambda^k+\gamma^{-m}))(\eta_2+O(\eta_2^2+\lambda^m))=O(\lambda^m\gamma^{-m}+\lambda^{k+m})
\end{equation}
(recall that we assume $k - m\gg 0$, so $\lambda^k=o(\lambda^m)$ and $\gamma^{-k}=o(\gamma^{-m})$; condition $|\lambda\gamma|>1$ also implies that
$\gamma^{-m}=o(\lambda^{-m})$).

After expressing $\mu$ as a function of $\eta_1$ and $\eta_2$ from \eqref{eq:nontr:1} and plugging the result into \eqref{eq:nontr:2}, we obtain
\begin{equation}\label{eq:nontr:4}
0 = cx^+\lambda^m +d(\eta_1^2-\eta_2^2)+o(\eta_1^2+\eta_2^2) + o(\lambda^m).
\end{equation}
Let
\begin{equation}\label{eq:nontr:4a}
\hat{\eta}_1=\eta_1+O(\eta_1^2+\lambda^k+\gamma^{-m}) \quad\mathrm{and}\quad \hat{\eta}_2=\eta_2+O(\eta_2^2+\lambda^m),
\end{equation}
where the $O(\cdot)$ terms are exactly those in the left-hand side of \eqref{eq:nontr:3}. Consequently, equations \eqref{eq:nontr:4} and \eqref{eq:nontr:3} become
\begin{eqnarray}
0&=& cx^+\lambda^m+d(\hat{\eta}_1^2-\hat{\eta}_2^2)+o(\lambda^m+\hat\eta_1^2+\hat\eta_2^2),\label{eq:nontr:5}\\
\hat{\eta}_1\hat{\eta}_2&=&(C_{k,m}  + O(|\hat\eta_1|+|\hat\eta_2|))\lambda^{k+m},\label{eq:nontr:6}
\end{eqnarray}
where $C_{k,m}$ is independent of $\hat\eta_{1,2}$ and uniformly bounded for all $k$ and $m$.

After we rescale the variables as follows:
\begin{equation}\label{eq:nontr:7}
(\hat{\eta}_1,\hat{\eta}_2)= (\lambda^{k+\frac{m}{2}}\xi_1,\lambda^{\frac{m}{2}}\xi_2),
\end{equation}
equations \eqref{eq:nontr:5} and \eqref{eq:nontr:6} transform to
\begin{equation}\label{eq:nontr:8}
\begin{array}{rcl}
0&=&cx^+ -d\xi^2_2+\dots, \\[5pt]
\xi_1\xi_2&=&C_{k,m}+\dots,
\end{array}
\end{equation}
where the dots denote terms that tend to zero as $k,m\to+\infty$.
By noting $cdx^+>0$ from the assumption of the lemma, we find, for all sufficiently large $k$ and $m$, two solutions
\begin{equation}\label{eq:nontr:9}
(\xi_1^*,\xi_2^*)=\pm\bigg(C_{k,m}\sqrt{\dfrac{d}{cx^+}} +
o(1)_{k,m\to +\infty},\sqrt{\dfrac{cx^+}{d}}+o(1)_{k,m\to +\infty}\bigg).
\end{equation}
Then, the corresponding values of $(\eta_1,\eta_2)$ can be found from \eqref{eq:nontr:7}, \eqref{eq:nontr:4a} and the corresponding 
values of $\mu$ can be found from either of the equations \eqref{eq:nontr:1}, \eqref{eq:nontr:2}.

We proceed to seek for the intersection $W^s(Q)\cap W^u(O)$. For an index-2 point, its local stable manifold is a leaf of $\mathcal{F}^s$.
In particular, a formula for the leaf through $Q_{02}$ is given by Lemma \ref{lem:ssleaves} as
\begin{equation*}\label{eq:nontr:10}
\begin{array}{rcl}
x&=&x_{02}+\varphi_1(z; x_{02}, y_{02},z_{02}), \\[5pt]
y&=&y_{02}+\varphi_2(z; x_{02}, y_{02},z_{02}),
\end{array}
\end{equation*}
where $\varphi_1=O(\lambda_0^{m}\lambda^{-m})$ and $\varphi_2=O(\hat{\lambda}^{m}\gamma^{-m})$. Here $\lambda_0$ is a value close to
$|\lambda_1|$ such that $|\lambda_1|<\lambda_0$.

Now let $\mathcal{W}=\{(0,y,0)\mid |y+y^-|<\varepsilon\}$ with $\varepsilon>0$ be a small piece of $W^u_{loc}(O)$ containing
the point $\tilde{M}^-=(0,-y^-,0)$. 
By formula \eqref{eq:setting:T_11}, the image $\tilde{T}_1(\mathcal{W})$ is given by
\begin{equation*}
\begin{array}{rcl}
x-x^+ &=& bt+h_1(0,t,0),\\[5pt]
-y &=& \mu + dt^2+ h_2(0,t,0),\\[5pt]
\mathcal{S}^{-1} z-z^+&=&a_{32}t+h_3(0,t,0),
\end{array}
\end{equation*}
where $t\in(-\varepsilon,\varepsilon)$. Hence, we can write the condition for the intersection of $W^u(O)\cap W^s(Q)$ as
\begin{eqnarray}
b\eta_1 + O(\lambda_0^m\lambda^{-m})&=& bt+h_1(0,t,0),\nonumber \\[5pt]
-\gamma^{-m}(\eta_2+y^-)+o({\gamma}^{-m})+O(\hat{\lambda}^{m}\gamma^{-m})) &=& \mu + dt^2+ h_2(0,t,0),\nonumber\\[5pt]
\mathcal{S}^{-1} z-z^+&=&a_{32}t+h_3(0,t,0),\nonumber
\end{eqnarray} 
which can be rewritten as
\begin{equation}\label{eq:nontr:11}
-\gamma^{-m}(\eta_2+y^-) + o({\gamma}^{-m})=\mu+d(\eta_1+O(\lambda_0^m\lambda^{-m}))^2+h_2(0,\eta_1+O(\lambda_0^m\lambda^{-m}),0).
\end{equation}
Since the $x$- and $z$-coordinates of the points in the orbit of $Q$ can be expressed as functions of $\eta_1$, $\eta_2$ and $\mu$, the right-hand side of the above equation is just a function of $\eta_1$, $\eta_2$ and $\mu$. We now express $\mu$ from \eqref{eq:nontr:11} as a function of $\eta_1$ and $\eta_2$, and obtain
\begin{equation}\label{eq:nontr:12}
\mu+d\eta_1^2+h_2(0,\eta_1,0)=-\gamma^{-m}y^-+O(\eta_1\lambda_0^m\lambda^{-m})+O(\lambda_0^{2m}\lambda^{-2m})+o({\gamma}^{-m}),
\end{equation}
which, along with \eqref{eq:nontr:2}, yields
\begin{equation}\label{eq:nontr:13}
\gamma^{-m}y^-=\dfrac{c}{2}\lambda^k x^+ +O(\eta_1\lambda_0^m\lambda^{-m})+O(\lambda_0^{2m}\lambda^{-2m})+
o({\gamma}^{-m})+o({\lambda}^k).
\end{equation}

Recall that condition (C4) gives $|\lambda|\;|\gamma|^{\frac{1}{2}}<1$. This, together with the fact $\lambda_0<\lambda^2$ given by 
Lemma \ref{lem:firstd}, implies $O(\lambda_0^{2m}\lambda^{-2m})=O(\lambda^{2m})=o(\gamma^{-m})$. By equations \eqref{eq:nontr:4a}, \eqref{eq:nontr:7} and \eqref{eq:nontr:9}, we have $\eta_1=O(\lambda^k+\gamma^{-m})$, which implies 
$O(\eta_1\lambda_0^m\lambda^{-m})=O(\eta_1\lambda^m)=o(\lambda^k)+o(\gamma^{-m})$.
With these observations, equation \eqref{eq:nontr:13} can be rewritten as
\begin{equation*}\label{eq:nontr:14}
\gamma^{-m}y^-=\dfrac{c}{2}\lambda^k x^+ +o({\lambda}^k) + o(\gamma^{-m}),
\end{equation*}
or
\begin{equation}\label{eq:nontr:15}
\lambda^k\gamma^m=\dfrac{2y^-}{cx^+} +o(\lambda^k\gamma^m)+o(1).
\end{equation}
Recall the assumption $2y^-/cx^+>0$; we also have taken $k$ and $m$ even, so both sides of equation (\ref{eq:nontr:15}) are positive.
We, therefore, may take logarithm on both sides, which gives
\begin{equation}\label{eq:nontr:16}
\theta=-\dfrac{\ln |\lambda|}{\ln|\gamma|}=\dfrac{m}{k} - \dfrac{C^*_{k,m}}{k},
\end{equation}
where $C^*_{k,m}=\ln{(2y^-/cx^+ +o({\lambda}^k\gamma^m))}/\ln|\gamma|$ is uniformly bounded, for all sufficiently large $k$ and $m$. Note that
$C^*$ is a function of $\theta$ -- it depends, for example, on the coefficients of the global and local maps, which depend on $\theta$ as a parameter.
It is important for us that $C^*$ is continuous and bounded function of $\theta$, so a value of $\theta$ that solves \eqref{eq:nontr:16} can be found 
for each sufficiently large $(k,m)$. It is also obvious, that if the sequence $\{(k_j,m_j)\}$ satisfies $k_j,m_j\to +\infty$ and $m_j/k_j\to \theta^*$ as 
$j\to +\infty$, then the values of $\theta_j$ we obtain from \eqref{eq:nontr:16} accumulate on $\theta=\theta^*$. The corresponding $\mu$ values are obtained from \eqref{eq:nontr:12} as $\mu_j=-\gamma^{m_j}y^-+o(\gamma^{m_j})$ and they tend to $0$ as $j\to +\infty$. Therefore, for each sufficiently large $j$, the map $F_{\mu_j,\theta_j}$ has an index-2 point of period 2 such that $W^s(Q_j)\cap W^u(O)\neq\emptyset$.\qed

\subsection{Proof of Theorem \ref{thm2}}\label{proofofthm2}
In the non-symmetric case, we use a simpler construction than in Theorem \ref{thm1}. In particular, we use
a different version of Lemma \ref{lem:nontransverse}. Here we have two splitting parameters 
$\mu_1$ and $\mu_2$, which correspond to two different orbits $\Gamma$ and $\tilde{\Gamma}$ of homoclinic tangency, respectively.
\begin{lem}\label{lem:nontransverse2}
Consider a two-parameter family $\{F_{\mu_1,\mu_2}\}$ of diffeomorphisms in $\mbox{Diff}\,^r(\mathcal{M})$, where $F_{0,0}$
satisfy conditions (C1)-(C3) and (C5). For every sufficiently large $k$ there exist parameter values
$\{(\mu^1,\mu^2)\}$, accumulating on 0 as $k\to+\infty$, such that the diffeomorphism $F_{\mu^1,\mu^2}$ 
has an index-2 periodic point $Q$ satisfying $T_1\circ T_0^{k}(Q)=Q$ and $W^s(Q)\cap W^u(O)\neq\emptyset$.
\end{lem}

Since Lemma \ref{lem:transverse} remains true in the general case, Theorem \ref{thm2} follows immediately by replacing Lemma \ref{lem:nontransverse} with Lemma \ref{lem:nontransverse2} in the proof of Theorem \ref{thm1}.

In what follows we prove Lemma \ref{lem:nontransverse2}. Here we consider the local map in the form \eqref{eq:setting:T_0} for which only identities in \eqref{eq:nonlinearterms} satisfied (as we do not have condition (C4) here, we cannot assume identities \eqref{eq:nonlinearterms1}). Therefore, without the identities in \eqref{eq:nonlinearterms1}, we do not have Lemmas \ref{lem:firstd} and \ref{lem:ssleaves}. We still
can use Lemma \ref{lem:ssleavesweak} which gives the equation for strong-stable leaves of a pointy $(x^*,y^*,z^*)$ in the form
\begin{equation}\label{eq:newssleaves}
\begin{array}{rcl}
x&=&x^*+\varphi_1(z; x^*,y^*,z^*), \\[5pt]
y&=&y^*+\varphi_2(z; x^*,y^*,z^*),
\end{array}
\end{equation}
where $\varphi_1=O(\hat{\lambda}^k\lambda^{-k})$ and $\varphi_2=O(\hat{\lambda}^{k}\gamma^{-k})$.

\noindent{\it Proof of Lemma \ref{lem:nontransverse2}.}
The coincidence condition (C5) implies that the small neighbourhoods $\tilde{\Pi}_1,\Pi_1$ and $\Pi_0$, and the local and global maps associated to the two homoclinic tangency orbits $\Gamma$ and $\tilde{\Gamma}$ can be defined in the same way as those in Section \ref{sec:returnmap}. The local map for $\Gamma$ and $\tilde{\Gamma}$ have the form of \eqref{eq:setting:T^k}. The two global maps $T_1$ and $\tilde{T}_1$ are given by
\begin{equation}\label{eq:nontr2:0}
\begin{array}{rcl}
x_0-x_i^+ &=& a_ix_1+b_i(y_1 - y_i^-) + a_{13}^i z_1 +h_{1}^i, \\ [5pt]
y_0 &=& \mu_i + c_ix_1 + d_i(y_1-y_i^-)^2 + a_{23}^i z_1 + h_{2}^i,\\[5pt]
z_0-z^+_i &=& a_{31}^ix_1+a_{32}^i(y_1 - y_i^-) + a_{33}^i z_1 +h_{3}^i, \\ [5pt]
\end{array}
\end{equation}
where $T_1$ corresponds to $i=1$ and $\tilde{T}_1$ corresponds to $i=2$. 

Let $Q$ be a periodic point such that $T_1\circ T_0^k(Q)=Q$. Denote
$Q_{0}=Q=(x_{0},y_{0},z_{0})$, $Q_{1}=T_0^k(Q)=(x_{1},y_{1},z_{1})$. 
By formulas \eqref{eq:setting:T^k} and \eqref{eq:nontr2:0}, the condition $Q=T_1\circ T_0^k(Q)$ is written as
\begin{eqnarray*}
\left.
\begin{array}{rcl}
x_{1}&=&\lambda^k x_{0} + \phi_k,\\[5pt]
y_{0}&=&\gamma^{-k} y_{1} +\psi_k,\\[5pt]
z_{1}&=&\hat{\phi}_k,
\end{array}
\right.
\quad
\left.
\begin{array}{rcl}
x_{0}-x^+ &=& a_1x_{1}+b_1(y_{1} - y^-) + a^1_{13} z_{1} +h^1_1, \\ [5pt]
y_{0} &=& \mu + c_1x_{1} + d_1(y_{1}-y^-)^2 + a^1_{23} z_{1} + h^1_2,\\[5pt]
z_{0}-z^+ &=& a^1_{31}x_{1}+a^1_{32}(y_{1} - y^-) + a^1_{33} z_{1} +h^1_3.
\end{array}
\right.
\end{eqnarray*}
We can express all variables here as functions of $y_1$, so the above equations reduce to
\begin{equation}\label{eq:nontr2:1a}
\gamma^{-k}(\eta+y^-) + o({\gamma}^{-k})=\mu_1+ c\lambda^k x^+ + bc\lambda^k \eta+ d\eta_2^2+o({\lambda}^k)+h_2(0,\eta_2,0),
\end{equation}
where we denote $\eta=y_1-y^-$. 

Like in the proof of Lemma \ref{lem:index2}, this fixed point is a saddle of index-2 if, for some $s\in(-1,1)$,
\begin{equation}\label{eq:index2:4aaa}
\dfrac{\tr \D (T_1\circ T_0^k)|_{E^{cu}}}{\det \D (T_1\circ T_0^k)|_{E^{cu}}+1}=s,
\end{equation}
where $E^{cu}$ is the two-dimensional invariant subspace given in Lemma \ref{lem:eigenvalues}. By formulas \eqref{eq:index2:6a}
and \eqref{eq:index2:8}, 
$$\D (T_1\circ T_0^k)|_{E^{cu}}=\left(\begin{array}{cc} o(1)_{k\to+\infty}  & b\gamma^k+o(\gamma^k)\\
c\lambda^k+o(|\lambda|^k+|\eta|)   &  \gamma^k (2d \eta+o(1)_{k\to+\infty} \eta +O(|\lambda|^k+\eta^2))
\end{array}\right)$$
(see \eqref{eq:index2:11}). Therefore, equation \eqref{eq:index2:4aaa} gives us the value of
\begin{equation}\label{etaeq}
\eta = O(\lambda^k).
\end{equation}
After that, we find $\mu$ from equation \eqref{eq:nontr2:1a} as
\begin{equation}\label{eq:nontr2:4}
\mu_1 = -c\lambda^k x^+ +o(\lambda^k).
\end{equation}

Let us now construct the intersection $W^s(Q)\cap W^u(O)$. Like in the proof of Theorem \ref{thm1}, this intersection is given by
\begin{equation}\label{eq:nontr2:5}
\begin{array}{rcl}
x_{0}-x^+ + O(\hat{\lambda}^k\lambda^{-k})&=& b_2t+h^2_{1}(0,t,0), \\[5pt]
-y_{0} +O(\hat{\lambda}^k\gamma^{-k}) &=& \mu_2 + d_2t^2+ h^2_{2}(0,t,0),\\[5pt]
 z-z_0^+&=&a_{32}^2t+h^2_{3}(0,t,0),
\end{array}
\end{equation}
which transforms into
\begin{equation}\label{eq:nontr2:6}
-\gamma^{-k}(\eta+y_1^-) + o({\gamma}^{-k})
=\mu_2+d_2\eta^2+O(\eta\hat{\lambda}^k\lambda^{-k})+O(\hat{\lambda}^{2k}\lambda^{-2k})+h^2_{2}(0,\eta,0).
\end{equation}
This along with \eqref{etaeq} gives us the corresponding value of
\begin{equation}\label{eq:nontr2:7}
\mu_2=O(\lambda^{2k}+|\gamma|^{-k}).
\end{equation}
The lemma follows immediately.
\qed

\appendix
\renewcommand\thesection{Appendix}
\section{}
\setcounter{equation}{0}
\renewcommand{\theequation}{A\arabic{equation}}
\renewcommand{\thesubsection}{A.\arabic{subsection}}

Here we show that, for the map $T_0$ in the form \eqref{eq:setting:T_0}, there exists a coordinate transformation $\mathcal{T}$ such that after this transformation the map $T_0$ will satisfy identities \eqref{eq:nonlinearterms} and \eqref{eq:nonlinearterms1}, and also keep the symmetry $\mathcal{R}$.

This transformation is constructed as a composition of\\
\indent $\mathcal{T}_1$ which straightens the local stable and unstable manifolds of $O$, thus giving the first two identities in \eqref{eq:nonlinearterms};\\
\indent $\mathcal{T}_2$ which linearises both the restriction $T_0|_{W^u_{loc}}$ and the quotient of $T_0|_{W^s_{loc}}$ by the the strong-stable foliation -- after that the third and forth identities in \eqref{eq:nonlinearterms} become valid;\\
\indent $\mathcal{T}_3$ which gives the last two identities in \eqref{eq:nonlinearterms}; and\\
\indent $\mathcal{T}_4$ which straightens a certain local, $\mathcal{R}$-symmetric extended unstable manifold $W_{loc}^{uE}(O)$ along with the foliation $\mathcal{F}^{uE}$ on it -- this leads to identities \eqref{eq:nonlinearterms1}.

In what follows we discuss the transformations $\mathcal{T}_i$, ($i=1,2,3,4$) separately and show that they keep the system symmetric with respect to $\mathcal{R}$, i.e., they commute with $\mathcal{R}$.

\subsection{Transformation $\mathcal{T}_1$}\label{app:t1}
Let $(x=w_{ux}(y), z=w_{uz}(y))$ and $y=w_s(x,z)$ be the equations for the local unstable and stable invariant manifolds of $O$, respectively. The transformation $\mathcal{T}_1$ is defined as
\begin{equation}\label{eq:app:1}
(x^{new},z^{new}) = (x-w_{ux}(y), z-w_{uz}(y)), \qquad y^{new} = y-w_s(x,z).
\end{equation}
In the new coordinates, the manifolds $W^u_{loc}$ and $W^s_{loc}$ have equations $(x^{new},z^{new})=0$ and, respectively $y^{new}=0$. Thus,
the first two identities in \eqref{eq:nonlinearterms} follow immediately from the invariance of these manifolds with respect to $T_0$.

Let us show that the coordinate transformation given by \eqref{eq:app:1} commutes with $\mathcal{R}$. Consider first the transformation 
$\psi: (x,y,z)\mapsto (x,y-w_s(x,z),z)$. By uniqueness of the stable manifold, $\mathcal(R) W^s_{loc}=W^s_{loc}$. Therefore, for any $x,z$, the image by $\mathcal R$ of the point $(x,w_s(x,z),z)\in W^s_{loc}$ also lies in $W^s_{loc}$, i.e.,
\begin{equation}\label{eq:app:2}
w_s(x,\mathcal{S}z)=-w_s(x,z).
\end{equation}
(see formula \eqref{eq:symmetryR} for $\mathcal R$). Similarly, by the uniqueness of the unstable manifold,
\begin{equation}\label{eq:app:2a}
w_{ux}(-y)=w_{ux}(y), \qquad w_{uz}(-y)=\mathcal{S} w_{uz}(y).
\end{equation}
Now let $(x,y,z)$ be an arbitrary point in a neighbourhood of $O$. We have
$$\mathcal{R}\circ \mathcal{T}_1(x,y,z)=(x-w_{ux}(y),-y+w_s(x,z), \mathcal{S}z - \mathcal{S} w_{uz}(y)),$$
and
$$\mathcal{T}_1\circ \mathcal{R}(x,y,z)=(x-w_{ux}(-y),-y-w_s(x,\mathcal{S}z), \mathcal{S}z - w_{uz}(-y)).$$
By \eqref{eq:app:2},\eqref{eq:app:2}, this implies that $\mathcal{T}_1$ commutes with $\mathcal{R}$, as required.

\subsection{Transformations $\mathcal{T}_2$ and $\mathcal{T}_3$}\label{app:t2,3}
The construction of these two transformations is given in the proof of Lemma 6 in \citep{gst08}. Here we reconstruct them for 
our case and prove that they are $\mathcal{R}$-symmetric.

The transformation $\mathcal{T}_2$ in sought is in the form
\begin{equation}\label{eq:app:3}
x^{new}=x+h_1(x,z),\quad\quad y^{new}=y+h_2(y),\quad\quad z^{new}=z,
\end{equation}
where $h_{1}(0,0)=0,h_2(0)=0,\partial h_1(0,0)/\partial (x,z)=0$ and $\partial h_2(0)/\partial y=0$ (hence the first two identities in 
\eqref{eq:nonlinearterms} hold in the new coordinates). To obtain the identities
$$f_1(x,0,z)=0 \quad\mbox{and}\quad f_2(0,y,0)=0,$$
we must have $\bar{x}^{new}=\lambda x^{new}$ at $y=0$, and $\bar{y}^{new}=\gamma y^{new}$ at $(x,z)=0$, respectively. According to formula \eqref{eq:setting:T_0} for $T_0$, these conditions translate to
\begin{equation}\label{eq:app:4}
\begin{array}{rcl}
h_1(\bar{x},\bar{z})&=&\lambda h_1(x,z)-f(x,0,z),\\
h_2(\bar{y})&=&\gamma h_2(y)-f(0,y,0),
\end{array}
\end{equation}
where we denote here $\bar{x}=\lambda x +f_1(x,0,z)$, $\bar{y}=\gamma y +f_2(0,y,0)$, and $\bar{z}=A z + f_3(x,0,z)$.

It has been shown in \citep{gst08} that the above system has the following solution:
\begin{equation}\label{eq:app:6}
h_1(x,z)=\sum_{j=0}^{+\infty}\lambda^{-j-1}f_1(x_j,0,z_j)\quad\mbox{and}\quad h_2(y)=-\sum_{j=1}^{+\infty}\gamma^{j-1}f_2(0,y_j,0),
\end{equation}
where $\{(x_j,z_j)\}$ is the forward orbit of $(x,z)=:(x_0,z_0)$ under the restriction of the local map \eqref{eq:setting:T_0} to to $W^s(O)$, and $\{y_j\}$ is the backward orbit of $y=:y_0$ under the restriction of the local map to $W^u(O)$. The functions $h_1$ and $h_2$ given by \eqref{eq:app:6} are obviously $\mathcal{R}$-symmetric. We, therefore, proceed to the analysis of the transformation $\mathcal{T}_3$:
\begin{equation}\label{eq:app:13}
x^{new}=x+g_1(x,y),\quad\quad y^{new}=y+g_2(x,y,z),\quad\quad z^{new}=z+g_3(x,y),
\end{equation}
where $g_{1,3}$ vanish at $x=0$ and $y=0$ while $g_2$ equals to zero at $(x,z)=0$ and at $y=0$. These conditions ensure that $\mathcal{T}_3$ keeps the identities obtained previously. We need to achieve that
$$\dfrac{\partial f_1}{\partial x}(0,y,0)=0,$$
in the new coordinates, which is equivalent to
\begin{equation}\label{eq:app:14}
\dfrac{\partial (\bar{x}^{new}-\lambda x^{new})}{\partial x^{new}}(0,y^{new},0)=0.
\end{equation}
Since the first identity in \eqref{eq:nonlinearterms} ensures
$$\dfrac{\partial (\bar{x}^{new}-\lambda x^{new})}{\partial y^{new}}(0,y^{new},0)=0,$$
equation \eqref{eq:app:14} holds if and only if
$$\d (\bar{x}^{new}-\lambda x^{new})=0 \quad \mbox{when} \quad (x^{new},z^{new})=0 \;\;\mbox{and}\;\; \d z^{new}=0.$$
We have, from \eqref{eq:app:13}, that $(x,z)=0$ at $(x^{new},z^{new})=0$, and
$$\d z^{new}=\d z+ \dfrac{\partial g_3}{\partial x}(0,y) \d x,$$
so $\d z^{new}=0$ when
\begin{equation}\label{eq:app:14.2}
\d z=- \dfrac{\partial g_3}{\partial x}(0,y) \d x.
\end{equation}

Equations \eqref{eq:setting:T_0}, \eqref{eq:app:13}, and \eqref{eq:app:14.2} imply that, when $(x^{new},z^{new})=0$ and $\d z^{new}=0$, we have
\begin{equation}\label{eq:app:15}
\begin{array}{rcl}
&&\d (\bar{x}^{new}-\lambda x^{new})\\[5pt]
&=& \d (\bar{x}+g_1(0,\bar{y})-\lambda x -\lambda g_1(0,y) )\\[5pt]
&=& \d (f_1(0,y,0)+g_1(0,\bar{y}) -\lambda g_1(0,y) )\\[5pt]
&=& \dfrac{\partial f_1}{\partial x}(0,y,0) \d x -\dfrac{\partial f_1}{\partial z}(0,y,0)\dfrac{\partial g_3}{\partial x}(0,y)\d x\\[5pt]
&&+ \dfrac{\partial g_1}{\partial x}(0,\bar{y}) \left(\lambda\d x+\dfrac{\partial f_1}{\partial x}(0,y,0) \d x -\dfrac{\partial f_1}{\partial z}(0,y,0)\dfrac{\partial g_3}{\partial x}(0,y)\d x\right) -\lambda \dfrac{\partial g_1}{\partial x}(0,y)\d x.
\end{array}
\end{equation}
We need to find functions $g_1$ and $g_3$ such that the right-hand side of \eqref{eq:app:15} vanishes identically. Denote
\begin{equation}\label{eq:app:14.1}
\eta_1(y)=\dfrac{\partial g_1}{\partial x}(0,y) \quad\mbox{and}\quad \eta_3(y)=\dfrac{\partial g_3}{\partial x}(0,y),
\end{equation}
and equate the right-hand side of \eqref{eq:app:15} to zero. This gives the following condition:
\begin{equation}\label{eq:app:16}
\begin{array}{rcl}
\eta_1(\bar{y}) &=&\left(\lambda \eta_1(y)-\dfrac{\partial f_1}{\partial x}(0,y,0) + \dfrac{\partial f_1}{\partial z}(0,y,0)\eta_3(y) \right)\\[10pt]
&&\times \left(\lambda+\dfrac{\partial f_1}{\partial x}(0,y,0) -\dfrac{\partial f_1}{\partial z}(0,y,0)\eta_3(y)\right)^{-1},
\end{array}
\end{equation}
where we have used the fact $\bar{x}=0$ at $(x,z)=0$, and $\bar{y}=\gamma y+f_2(0,y,0)$.

Analogously, we will have identity
$$\dfrac{\partial f_3}{\partial z}(0,y,0)=0$$
satisfied in the new coordinates, if
\begin{equation}\label{eq:app:16.1}
\begin{array}{rcl}
\eta_3(\bar{y}) &=&\left(A \eta_3(y)-\dfrac{\partial f_3}{\partial x}(0,y,0) + \dfrac{\partial f_3}{\partial z}(0,y,0)\eta_3(y)\right)\\[10pt]
&&\times \left(\lambda+\dfrac{\partial f_1}{\partial x}(0,y,0) -\dfrac{\partial f_1}{\partial z}(0,y,0)\eta_3(y)\right)^{-1},
\end{array}
\end{equation}
Equations \eqref{eq:app:16} and \eqref{eq:app:16.1} are solved by noticing that they can be viewed as the conditions for the manifold
\begin{equation}\label{eq:app:16.2}
w_1:\{u_1=\eta_1(y), u_3=\eta_3(y)\}
\end{equation}
to be invariant under the map
\begin{equation}\label{eq:app:17a}
\begin{array}{rcl}
\bar{y}&=&\gamma y+f_2(0,y,0),\\[10pt]
\bar{u}_1 &=&\left(\lambda u_1-\dfrac{\partial f_1}{\partial x}(0,y,0) + \dfrac{\partial f_1}{\partial z}(0,y,0)u_3 \right) \left(\lambda+\dfrac{\partial f_1}{\partial x}(0,y,0) -\dfrac{\partial f_1}{\partial z}(0,y,0)u_3\right)^{-1},\\[10pt]
\bar{u}_3 &=&\left(A u_3-\dfrac{\partial f_3}{\partial x}(0,y,0) + \dfrac{\partial f_3}{\partial z}(0,y,0)u_3 \right) \left(\lambda+\dfrac{\partial f_1}{\partial x}(0,y,0) -\dfrac{\partial f_1}{\partial z}(0,y,0)u_3\right)^{-1}.
\end{array}
\end{equation}
Note that this map has a fixed point $(0,0,0)$. The multipliers of this point are the eigenvalues of the linearised map, which is given by
$$y\mapsto\gamma y, \quad\quad u_1\mapsto u_1 - \dfrac{\partial^2 f_1}{\partial x\partial y}(0,0,0)\lambda^{-1}y,\quad u_3\mapsto \lambda^{-1}A u_3 - \dfrac{\partial^2 f_3}{\partial x\partial y}(0,0,0)\lambda^{-1}y.$$
The spectrum of this map consists of the spectra of the following three operators: $y\mapsto \gamma y, u_1\mapsto u_1, u_3\mapsto \lambda^{-1}A u_3$. Therefore, the fixed point $(0,0,0)$ has one multiplier on the unit circle, one multiplier outside the unit circle and $(n-2)$ multipliers inside the unit circle. It has been known (see e.g. \citep{hps77,sstc1}) that such fixed point lies in a unique one-dimensional unstable manifold that is tangent at zero to the eigenspace corresponding to the multiplier outside the unit circle. It follows that such unique manifold in our case is the sought manifold $w_1$.

The map \eqref{eq:app:17a} is symmetric with respect to $(y,u_1,u_3)\mapsto (-y,u_1,\mathcal{S}u_3)$. Indeed, this follows immediately from the relations
\begin{equation*}
\begin{array}{l}
\dfrac{\partial f_1}{\partial x}(0,-y,0)=\dfrac{\partial f_1}{\partial x}(0,y,0),\quad
\dfrac{\partial f_1}{\partial z}(0,-y,0)\mathcal{S}=\dfrac{\partial f_1}{\partial z}(0,y,0),\\[10pt]
\dfrac{\partial f_3}{\partial x}(0,-y,0)=\mathcal{S}\dfrac{\partial f_3}{\partial x}(0,y,0),\quad
\dfrac{\partial f_3}{\partial z}(0,-y,0)\mathcal{S}=\mathcal{S}\dfrac{\partial f_3}{\partial z}(0,y,0),
\end{array}
\end{equation*}
which are, in turn, implied by the symmetry of $T_0$ with respect to $\mathcal{R}$. By uniqueness of $w_1$, it must be symmetric with respect to the transformation $(y,u_1,u_3)\mapsto (-y,u_1,\mathcal{S}u_3)$, which implies that $\eta_{1,3}$ are symmetric with respect to $y\mapsto -y$. Consequently, 
functions $g_{1,3}(x,z)$ can be any of those that vanish at $(x,z)=0$ and $y=0$ and satisfy \eqref{eq:app:14.1}. Due to the symmetry of $\eta_{1,2}$,
it is easy to show that $g_{1,3}$ can be chosen symmetric with respect to $(x,y,z)\mapsto (x,-y,\mathcal{S}z)$, as required.

The next identity to be satisfied in the new coordinates is
\begin{equation}\label{abovei}
\dfrac{\partial f_2}{\partial y} (x,0,z)=0.
\end{equation}
Similarly to the above, by letting
\begin{equation}\label{eq:app:17}
\eta_2(x,z)=\dfrac{\partial g_2}{\partial y}(x,0,z),
\end{equation}
the identity \eqref{abovei} is equivalent to
\begin{equation}\label{eq:app:18}
\eta_2(\bar{x},\bar{z}) =\left(\gamma \eta_2(x,z)-\dfrac{\partial f_2}{\partial y}(x,0,z) \d y \right)\left(\gamma+\dfrac{\partial f_2}{\partial y}(x,0,z) \d y\right)^{-1}.
\end{equation}
This is the condition for the manifold $w_2:v=\eta_2(x,z)$ to being invariant under the map
\begin{equation*}\label{eq:app:19}
\bar{x}=\lambda x+f_1(x,0,z),\quad
\bar{z}=Az + f_3(x,0,z)\quad
\bar{v} =\left(\gamma v-\dfrac{\partial f_2}{\partial y}(x,0,z) \d y \right)\left(\gamma+\dfrac{\partial f_2}{\partial y}(x,0,z) \d y\right)^{-1}.
\end{equation*}
This map is symmetric with respect to $(x,z,v)\mapsto (x,\mathcal{S}z,v)$, and has a unique $(n-1)$-dimensional stable invariant manifold. It follows that 
$\eta_2$ exists and is symmetric with respect to $(x,z)\mapsto (x,\mathcal{S}z)$. The function $g_2(x,y,z)$ can be any of those that vanish at 
$(x,z)=0$ and $y=0$ and satisfy \eqref{eq:app:17}. The symmetry of $\eta_2$ implies that $g_2$ can be chosen symmetric with respect to
$(x,y,z)\mapsto (x,-y,\mathcal{S}z)$, so we can now conclude that $\mathcal{T}_3$ is $\mathcal{R}$-symmetric.

\subsection{Transformation $\mathcal{T}_4$}\label{app:t4}
Recall that $\mathcal{T}_4$ is a transformation that straightens the extended-unstable invariant manifold $W^{uE}(O)$ of $O$ and the foliation on it. 
This manifold is not unique and we choose a special one as follows.

We consider the following map $G_0$: 
\begin{equation}\label{eq:smoothness:2}
\begin{array}{rcl}
\bar{x} &=& \lambda x +  f_1(x,y,z), \\
\bar{y} &=& \gamma y +f_2(x,y,z),\\
\bar{z} &=& Az + f_3(x,y,z),\\[5pt]
\bar{u} &=& \Big(\Big(\lambda+\dfrac{\partial f_1}{\partial x}\Big)u+\dfrac{\partial f_1}{\partial y}+\dfrac{\partial f_1}{\partial z}v\Big)
\Big(\gamma+\dfrac{\partial f_2}{\partial y}+\dfrac{\partial f_2}{\partial x}u+\dfrac{\partial f_2}{\partial z}v\Big)^{-1},\\[10pt]
\bar{v} &=& \Big(\Big(A+\dfrac{\partial f_3}{\partial z}\Big)v+\dfrac{\partial f_3}{\partial x}u+\dfrac{\partial f_3}{\partial y}\Big)
\Big(\gamma+\dfrac{\partial f_2}{\partial y}+\dfrac{\partial f_2}{\partial x}u+\dfrac{\partial f_2}{\partial z}v\Big)^{-1},
\end{array}
\end{equation}
where the first three lines give the map $T_0$, e.g. the functions $f_i$ satisfy all the identities in \eqref{eq:nonlinearterms}. 
Obviously, this map is $C^{r-1}$ smooth and $\mathcal{R}$-symmetric. Note that $G_0$ 
is defined in $V_0\times \mathbb{R}^{1+D}$, where $V_0$ is the domain of $T_0$. We now extend $G_0$ to the whole of $\mathbb{R}^{2D+3}$ by replacing the functions $f_i$ $(i=1,2,3)$ in \eqref{eq:smoothness:2} with $f_i(\xi(x,y,z))$, were $\xi$ is a $C^r$ function such that, for two small numbers 
$\delta_1,\delta_2>0$ with $\delta_1<\delta_2$, we have 
\begin{equation*}
\xi(x,y,z)=
\begin{cases}
(x,y,z) &\mathrm{if}\,\, \|(x,y,z)\|<\delta_1\\
0 &\mathrm{if}\,\, \|(x,y,z)\|>\delta_2
\end{cases}.
\end{equation*}
For simplicity we use the same notation for the new functions $f_i$ so that the extension map, denoted by $G$, assumes the same form as \eqref{eq:smoothness:2}. One can choose the function $\xi$ such that the map $G$ will be $\mathcal{R}$-symmetric.

It can be seen from \eqref{eq:smoothness:2} that $G$ has a fixed point at zero, and the corresponding multipliers are $\lambda$, $\gamma$
(these correspond to variables $x$ and $y$), the eigenvalues of $A$ (which corresponds to variables $z$) and also $\lambda/\gamma$ (corresponding to the 
variable $u$) and the eigenvalues of $A$ divided by $\gamma$ (corresponding to the variables $v$). Since $|\gamma|>1$, it follows that the eigenvalues corresponding to the variables $z,u,v$ are smaller in absolute value that $\max\{\tilde\lambda, \lambda/\gamma\}$ where
$\tilde{\lambda}>0$ is a value close to $|\lambda_1|$, the largest absolute value of the eigenvalues of $A$. Since $|\gamma|>1>|\lambda|$ and $|\lambda| > \tilde\lambda$, we have that there is a spectrum dichotomy between $x,y$ variables and $z,u,v$ variables. The other assumption $|\lambda\gamma|>1$ and 
$|\lambda_1|<\lambda^2$ (so, $\tilde{\lambda}<\lambda^2$) further implies
\begin{equation}\label{eq:smoothness:3}
\left|\dfrac{\ln|\lambda\gamma^{-1}|}{\ln|\lambda|}\right|>\left|\dfrac{\ln\lambda^2}{\ln|\lambda|}\right|=2
\quad\mbox{and}\quad
\left|\dfrac{\ln\tilde{\lambda}}{\ln|\lambda|}\right|>\left|\dfrac{\ln\lambda^2}{\ln|\lambda|}\right|=2.
\end{equation}
Thus, the spectrum gap $l$ between $(x,y)$ and $(z,u,v)$ is greater than 2. It follows that there exists a unique invariant $C^2$-manifold $W_G$ 
for the map $G$, and it attracts all the orbits near it (see e.g. Section 5 of \citep{sstc1} and \citep{hps77}). This manifold has the form
\begin{equation}\label{eq:smoothness:4}
\begin{array}{rcl}
z&=&\eta_{uE}(x,y),\\
u&=&\eta_1(x,y),\\
v&=&\eta_2(x,y).
\end{array}
\end{equation} 
We now take any surface $w$ of the form \eqref{eq:smoothness:4} such that it is $\mathcal{R}$-symmetric and satisfies
\begin{equation}\label{eq:smoothness:5}
\eta_2=\dfrac{\partial \eta_{uE}}{\partial x}\eta_1 + \dfrac{\partial \eta_{uE}}{\partial y}.
\end{equation}
This equation means that the line field given by $(\eta_1,1,\eta_2)$ belongs to the tangent space of the surface $z=\eta_{uE}(x,y)$. Since $W_G$ is attracting, the iterates $G^n(w)$ tend to $W_G$ as $n\to +\infty$. It is easy to check that each iteration will be again $\mathcal{R}$-symmetric and will satisfy
\eqref{eq:smoothness:5}. Therefore, the limit $W_G$ is $\mathcal{R}$-symmetric and satisfies \eqref{eq:smoothness:5}. By our construction, the extended-unstable manifold $W^{uE}(O)$ is given by the $C^2$ function $z=\eta_{uE}(x,y)$. A $C^2$ foliation $\mathcal{F}^{uE}$ on $W^{uE}(O)$ can be found by integrating the line field given by $(\eta_1,1,\eta_2)$. Namely, it consists of solutions to the system of differential equations 
$\dot{x}=\eta_1,\,\,\dot{y}=1,\,\,\dot{z}=\eta_2.$

We can now define $\mathcal{T}_4$ as the composition of two transformations which straighten the manifold $W^{uE}(O)$ and the leaves of $\mathcal{F}^{uE}$, respectively. The former can be obtained by the same way as we did for $\mathcal{T}_1$, and it will be $C^2$ and $\mathcal{R}$-symmetric. Regarding the latter, we explain as follows.

Parametrize the leaves by its intersection with $\{y=0\}$, which is denoted by $c$. Then, the leaf of $\mathcal{F}^{uE}$ that goes through the point $(c,0,0)$ is given by $(x,z)=h(y,c)=:(h_1(y,c), h_2(y,c))$, where $h_i$ are $C^2$ functions. The foliation $\mathcal{F}^{uE}$ also induces a $C^2$ function
\begin{equation*}
g:\mathbb{R}^2\to \mathbb{R},(x,y)\mapsto c,
\end{equation*}
where $c$ satisfies $x=h_1(y,c)$. In order to linearise the quotient map along the leaves of this foliation (i.e., to straighten the leaves), we use the following $C^2$ transformation:
\begin{equation}
 x^{new}= g(x,y),\quad y^{new}=y,\quad z^{new}=z.
\end{equation}
\par{}
Note that the foliation $\mathcal{F}^{uE}$ is $\mathcal{R}$-symmetric. This implies $h_1(y,c) =  h_1(-y,c)$ and $g(x,y) = g(x,-y)$. Consequently, the above transformation is $\mathcal{R}$-symmetric. Therefore, the transformation $\mathcal{T}_4$ is $C^2$ and $\mathcal{R}$-symmetric.
\begin{rem}
The transformation $\mathcal{T}_4$ is $C^1$-smooth in parameters. This can be seen by letting $\varepsilon$ be the vector of all parameters, and then adding $\bar{\varepsilon}=\varepsilon$ into system \eqref{eq:smoothness:2}. 
\end{rem}

\end{document}